\overfullrule=0pt
\let\printlabel=N
\tolerance=700
\def\fnt#1{\csname !f-#1\endcsname}
\def\deffnt#1{\expandafter\font\csname !f-#1\endcsname}
\def\scal#1{ scaled \magstep#1}
\deffnt{rm}=cmr10
\deffnt{rm.1}=cmr10 \scal1
\deffnt{rm.2}=cmr10 \scal2
\deffnt{rm9}=cmr9
\deffnt{rm8}=cmr8
\deffnt{rm7}=cmr7
\deffnt{rm5}=cmr5
\deffnt{bf}=cmbx10
\deffnt{bf.1}=cmbx10 \scal1
\deffnt{bf.2}=cmbx10 \scal2
\deffnt{bf9}=cmbx9
\deffnt{bf7}=cmbx7
\deffnt{bf5}=cmbx5
\deffnt{it}=cmti10
\deffnt{it.1}=cmti10 \scal1
\deffnt{it9}=cmti9
\deffnt{it7}=cmti7
\deffnt{mi}=cmmi10
\deffnt{mi.1}=cmmi10 \scal1
\deffnt{mi9}=cmmi9
\deffnt{mi7}=cmmi7
\deffnt{mi5}=cmmi5
\deffnt{sy}=cmsy10
\deffnt{sy9}=cmsy9
\deffnt{ex}=cmex10
\deffnt{ex9}=cmex9
\deffnt{sl}=cmsl10
\deffnt{sl8}=cmsl8
\deffnt{sf}=cmss10
\deffnt{sf.2}=cmss10 \scal2
\deffnt{sc}=cmcsc10
\deffnt{sc9}=cmcsc9
\deffnt{sc.1}=cmcsc10 \scal1
\deffnt{mb}=msbm10
\deffnt{mb7}=msbm7
\deffnt{mb5}=msbm5
\newfam\mbfam
\textfont\mbfam=\fnt{mb}
\scriptfont\mbfam=\fnt{mb7}
\scriptscriptfont\mbfam=\fnt{mb5}
\def\mb{\fam\mbfam\fnt{mb}}
\newfam\lgfam
\textfont\lgfam=\fnt{mi.1}
\scriptfont\lgfam=\fnt{mi9}
\scriptscriptfont\lgfam=\fnt{mi7}

\newfam\smfam
\textfont\smfam=\fnt{mi9}
\scriptfont\smfam=\fnt{mi7}
\scriptscriptfont\smfam=\fnt{mi5}
\deffnt{gt}=eufm10
\deffnt{gt5}=eufm5
\def\cmd#1{\csname #1\endcsname}
\def\defcmd#1{\expandafter\def\csname #1\endcsname}
\def\edefcmd#1{\expandafter\edef\csname #1\endcsname}
\def\ifundefined#1{\expandafter\ifx\csname #1\endcsname\relax}
\edef\warning#1{{\newlinechar=`|\message{|*********** #1}}}
\def\npar{\relax\par\noindent}
\newbox\xbox
\def\ifhempty#1#2#3{\setbox\xbox=\hbox{#1}\ifdim\wd\xbox=0pt\relax
                                                           #2\else #3\fi}
\def\lph#1{\hbox to 0pt{\hss #1}}
\def\rph#1{\hbox to 0pt{#1\hss}}
\def\ulph#1{\vbox to 0pt{\vss\hbox to 0pt{\hss #1}}}
\def\urph#1{\vbox to 0pt{\vss\hbox to 0pt{#1\hss}}}
\def\dlph#1{\vbox to 0pt{\hbox to 0pt{\hss #1}\vss}}
\def\drph#1{\vbox to 0pt{\hbox to 0pt{#1\hss}\vss}}
\def\nasp{\spacefactor=1000}
\def\small{\def\rm{\fnt{rm9}\fam\smfam}\def\it{\fnt{it9}}\def\bf{\fnt{bf9}}%
\textfont1=\fnt{mi9}\textfont2=\fnt{sy9}\textfont3=\fnt{ex9}%
\baselineskip=.8\baselineskip\parindent=.8\parindent%
\abovedisplayskip=.6\abovedisplayskip\belowdisplayskip=.6\belowdisplayskip%
\stateskip=.6\stateskip
\rm}
\def\nohyphen{\hyphenpenalty=5000 \exhyphenpenalty=5000
 \tolerance=5000 \pretolerance=5000}
\def\center{\leftskip=0pt plus 1fill\relax \rightskip=\leftskip\relax
 \parfillskip=0pt\relax \parindent=0pt\relax \nohyphen}
\def\\{\ifhmode\hfil\break\else\ifvmode\vskip\baselineskip\fi\fi}
\def\-{\hfill}
\def\vsp#1{\vadjust{\kern #1}}

\def\vbreak#1#2{\par\ifdim\lastskip<#2\removelastskip\penalty#1\vskip#2\fi}
\def\vbreakn#1#2{\vbreak{#1}{#2}\npar}

\newcount\multilabels \newcount\undeflabels
\newread\filelabelsold \newwrite\filelabels
\edef\filelabelsname{\jobname.lbl}
\def\getlabels{\openin\filelabelsold=\filelabelsname\relax
 \ifeof\filelabelsold 
  \warning{The file of labels \filelabelsname\ does not exist}
  \else\closein\filelabelsold
   {\catcode`/=0 \globaldefs=1 \input \filelabelsname\relax}\fi}
\def\openlabels{\immediate\openout\filelabels=\filelabelsname\relax}
\def\closelabels{\closeout\filelabels}

\def\label#1#2#3{
\ifx\printlabel Y\lnote{\fnt{rm7}#1}\fi%
\immediate\write\filelabels{%
/defcmd{!tl-#1}{#2}/defcmd{!nl-#1}{#3}}%
\ifundefined{!dl-#1}\defcmd{!dl-#1}{!earlier-defined!}%
\else\warning{Label #1 is multiply defined.}\advance\multilabels by 1\fi%
\edefcmd{!tl-#1}{#2}\edefcmd{!nl-#1}{#3}}
\def\xlabel#1#2#3#4{
\ifx\printlabel Y\ldnote{\fnt{rm7}#1}\fi%
\immediate\write\filelabels{%
/defcmd{!tl-#1}{#2}/defcmd{!nl-#1}{#3}/defcmd{!xl-#1}{#4}}%
\ifundefined{!dl-#1}\defcmd{!dl-#1}{!earlier-defined!}%
\else \warning{LABEL #1 is multiply defined.} \advance\multilabels by 1\fi%
\edefcmd{!tl-#1}{#2}\edefcmd{!nl-#1}{#3}\edefcmd{!xl-#1}{#4}}
\def\undelabel#1{\warning{Label #1 has not been defined.}%
\global\advance\undeflabels by 1{?#1?}}
\def\rfrn#1{\ifundefined{!nl-#1}\undelabel{#1}\else\cmd{!nl-#1}\fi}
\def\rfrt#1{\ifundefined{!tl-#1}\undelabel{#1}%
 \else\cmd{!tl-#1}\fi}
\def\rfrx#1{\ifundefined{!xl-#1}\undelabel{#1}%
 \else\cmd{!xl-#1}\fi}
\def\rfr#1{\ifundefined{!tl-#1}\undelabel{#1}%
 \else\cmd{!tl-#1}~\cmd{!nl-#1}\fi}

\def\title#1{\vbreak{-200}{1cm}{\center\fnt{bf.1}#1\par}}
\def\author#1{\vbreak{100}{5mm}{\center\fnt{rm}#1\par}}
\def\titledate{\vbreak{100}{5mm}\centerline{\fnt{it}\date}}

\newskip\absskip
\absskip=20mm
\long\def\abstract#1{\vbreak{50}{5mm}\centerline{\fnt{bf}Abstract}
 \kern3mm{\leftskip=\absskip\rightskip=\leftskip
 \small #1\vbreak{-9000}{4mm}}}

\def\support#1{\footnote{}{{\small#1\hfil}}}

\ifx\headdate Y{\headline={{\fnt{rm7}\date}\hfil}}\fi
\def\date{\ifcase\month%
\or January\or February\or March\or April\or May\or June\or July%
\or August\or September\or October\or November\or December\fi%
~\number\day,~\number\year}
\def\myaddress{\rm%
\setbox\xbox=\hbox{Department of mathematics}%
\vtop{\hsize=\wd\xbox\parindent=0pt
Department of Mathematics\\
The Ohio State University\\
Columbus, OH 43210, USA\\
\hbox to \wd\xbox{{\it e-mail}: leibman@math.ohio-state.edu\hss}}}

\def\acknowledgment#1{\vbreakn{-5000}{2mm}{\bf Acknowledgment.} #1}

\newcount\fnotenum
\def\fnote#1{\global\advance\fnotenum by 1
 {\parindent=10mm \parfillskip=0pt plus 1fill\relax
  \footnote{$^{(\the\fnotenum)}$}{{\small #1}}}}

\def\lnote#1{\ifvmode\ulph{#1 }\else\vadjust{\ulph{#1 }}\fi}
\def\ldnote#1{\ifvmode\ulph{#1 }\else\vadjust{\dlph{#1 }}\fi}
\newcount\sectionnum \sectionnum=-1 
\newcount\subsectionnum \subsectionnum=0 
\newcount\subsubsectionnum \subsubsectionnum=0
\newskip\sectionskip \sectionskip=6mm
\newskip\aftersectionskip \aftersectionskip=3mm
\newcount\sectionpenalty \sectionpenalty=-9000
\def\xsection#1{\global\advance\sectionnum by 1%
\global\subsectionnum=0\global\subsubsectionnum=0%
\ifx\nonewsecequation Y\else\global\equationnum=0\fi%
\ifx\nonewsecstate Y\else\global\statenum=0\fi%
\label{#1}{Section}{\the\sectionnum}}
\def\section#1#2{%
\vbreak{\sectionpenalty}{\sectionskip}\xsection{#1}%
\vbox{\ifx\nosectioncenter Y\npar\else\center\fi\fnt{bf}
\ifx\nosectionnum Y\else\rfrn{#1}.{\nasp} \fi#2\par}%
\vbreak{10000}{\aftersectionskip}}
\newskip\subsectionskip \subsectionskip=3mm
\newcount\subsectionpenalty \subsectionpenalty=-7000
\def\thesubsectionnum{\ifx\nosection Y\the\subsectionnum
\else\the\sectionnum.\the\subsectionnum\fi}
\def\xsubsection#1{\global\advance\subsectionnum by 1%
\global\subsubsectionnum=0%
\label{#1}{subsection}{\thesubsectionnum}}
\def\subsection#1{%
\vbreak{\subsectionpenalty}{\subsectionskip}
\ifx\nosubsectionnum Y\par\else
\noindent\xsubsection{#1}{\fnt{bf}\rfrn{#1}.{\nasp}}\fi}
\def\thesubsubsectionnum{\thesubsectionnum.\the\subsubsectionnum}
\def\xsubsubsection#1{\global\advance\subsubsectionnum by 1%
\label{#1}{subsection}{\thesubsubsectionnum}}
\def\subsubsection#1{%
\vbreak{-5000}{1.5mm}\noindent\xsubsubsection{#1}{\fnt{bf}\rfrn{#1}.{\nasp}}}
\newcount\equationnum   
\def\theequation{\ifx\nosection Y\the\equationnum
\else\the\sectionnum.\the\equationnum\fi}
\def\equn#1#2{\global\advance\equationnum by 1%
\label{#1}{Equation}{\theequation}%
$$\vcenter{\equalign{#2}}\eqno(\theequation)$$}
\def\equ#1{$$\vcenter{\equalign{#1}}$$}
\def\equalign#1{\let\\=\cr\let\-=\hfill
\ialign{&\hfil$\dsp ##$\hfil\cr#1\crcr}}
\def\lequ#1{$$\vcenter{\lequalign{#1}}$$}
\def\lequalign#1{\let\\=\cr\let\-=\hfill
\ialign{&\hbox to \hsize{$\dsp ##$\hfil}\cr#1\crcr}}
\def\matalign#1#2{{\let\\=\cr\let\-=\hfill%
\ialign{\hfil$##$\hfil&&#1\hfil$##$\hfil\cr#2\crcr}}}

\def\frfr#1{\hbox{(\rfrn{#1})}}
\newskip\stateskip\stateskip=2mm
\newskip\theoremskip\theoremskip=\stateskip
\newcount\statenum\statenum=0
\def\thestate{\ifx\afterstate Y\the\statenum\else\thesubsectionnum\fi}
\def\state#1#2#3#4{\global\advance\statenum by 1\xlabel{#2}{#1}{\thestate}{#3}%
{\bf #1\ifhempty{#3}{\ifx\afterstate Y\ \thestate\fi}{ \box\xbox}.}{\nasp} #4%
\begingroup}
\def\nstate#1#2{{\bf #1\ifhempty{#2}{}{ \box\xbox}.}{\nasp}\begingroup}
\def\endstate{\par\endgroup\vbreak{-7000}{\stateskip}}
\def\ltheorem#1#2#3{\state{Theorem}{#1}{#2}{#3}\it}
\def\theorem#1#2#3{\vbreakn{-7000}{\stateskip}\ltheorem{#1}{#2}{#3}}
\def\endtheorem{\endstate}
\def\llemma#1#2#3{\state{Lemma}{#1}{#2}{#3}\it}
\def\lemma#1#2#3{\vbreakn{-7000}{\stateskip}\llemma{#1}{#2}{#3}}
\def\endlemma{\endstate}
\def\lproposition#1#2#3{\state{Proposition}{#1}{#2}{#3}\it}
\def\proposition#1#2#3{\vbreakn{-7000}{\stateskip}\lproposition{#1}{#2}{#3}}
\def\endproposition{\endstate}



\def\lremark#1{\nstate{Remark}{#1}}
\def\remark#1{\vbreakn{-7000}{\stateskip}\lremark{#1}}
\def\endremark{\endstate}




\def\lproof#1{\nstate{Proof}{#1}}
\def\proof#1{\vbreakn{-7000}{\stateskip}\lproof{#1}}

\def\endproof{\endprr\endstate}
\def\enprule{\vrule height1mm depth1mm width2mm}
\def\endprr{\discretionary{}{\kern\hsize}{\kern 3ex}\llap{\enprule}}
\def\frgdsp{\par\penalty10000
 \vskip-\belowdisplayskip\kern-2mm\noindent\hbox to \hsize{\hfil}}
\newcount\biblnum \newcount\biblpenalty
\def\bibliography#1{\vbreak{-9000}{6mm}
\biblpenalty=10000\biblnum=0%
\line{\fnt{bf}\ifhempty{#1}{Bibliography}{\box\xbox}\hfil}\par\kern2mm
\begingroup\parskip=0pt\parindent=0pt\frenchspacing}
\def\endbibliography{\par\endgroup}
\newskip\biblleft \biblleft=10mm
\def\bibxitem#1/#2(#3) #4{\advance\biblnum by 1%
\vbreak{\biblpenalty}{.8mm}\biblpenalty=-9000%
\hangindent=\biblleft \hangafter=1
\noindent{\small\rlap{\brfr{#1}}\hskip\hangindent
\xlabel{#1}{#3}{\the\biblnum}{#2}#4}}
\def\bibitem#1/#2(#3) #4#5{\bibxitem#1/#2(#3) {{\frenchspacing #4}, #5}}
\def\bibx#1/#2 {\advance\biblnum by 1{#1/#2}
\xlabel{#1}{text}{\the\biblnum}{#2}}
\def\bibart#1/#2 a:#3 t:#4 j:#5 n:#6 y:#7 p:#8 *{%
\bibitem{#1}/{#2}(paper) {#3}{{#4}, {\it #5\/}
   {\bf #6} (#7), #8.}}
\def\bibartp#1/#2 a:#3 t:#4 *{%
\bibitem{#1}/{#2}(paper) {#3}{{#4}, {in preparation}.}}
\def\bibarts#1/#2 a:#3 t:#4 *{%
\bibitem{#1}/{#2}(paper) {#3}{{#4}, {submitted}.}}
\def\bibarta#1/#2 a:#3 t:#4 j:#5 *{%
\bibitem{#1}/{#2}(paper) {#3}{{#4}, {\rm to appear in \it #5}.}}
\def\bibartn#1/#2 a:#3 t:#4 j:#5 *{%
\bibitem{#1}/{#2}(paper) {#3}{{#4}, #5.}}
\def\bibartpr#1/#2 a:#3 t:#4 *{%
\bibitem{#1}/{#2}(paper) {#3}{{#4}, {preprint}.}}
\def\bibartx#1/#2 a:#3 t:#4 x:#5 *{%
\bibitem{#1}/{#2}(paper) {#3}{{#4}, #5.}}
\def\bibook#1/#2 a:#3 t:#4 i:#5 *{%
\bibitem{#1}/{#2}(book) {#3}{{\it #4\/}, #5.}}
\def\no#1{no.{\nasp}~{#1}}

\def\brfr#1{\hbox{\rm[\ifx\biblNUM Y\rfrn{#1}\else\rfrx{#1}\fi]}}
\def\start{
 \warning{|****** Start ******|}
 \getlabels\openlabels\multilabels=0\undeflabels=0
 \null}
\def\finish{\par\closelabels
 \warning{|****** Finish ******}
 \ifnum\multilabels>0\warning{\the\multilabels_multidefined labels!}\fi
 \ifnum\undeflabels>0\warning{\the\undeflabels_undefined labels!}
                     \warning{*************************}
                     \warning{Try to run TeX once again!}
                     \warning{*************************}\fi}
\output={\shipout\vbox{\makeheadline\pagebody\makefootline}
\ifx\doublepage Y\shipout\vbox{}\fi
 \advancepageno
 \ifnum\outputpenalty>-20000\else\dosupereject\fi}
\def\sdup#1#2{{\scr #1 \atop \scr #2}}

\def\Gal{\lnote{{\bf V\kern 1cm\relax}}}
\def\comment#1{\ifhmode\\\fi
\hbox to \hsize{\hss\vbox{\advance\hsize by 10mm
\baselineskip=.7\baselineskip\npar\ulph{\fnt{bf.2}V}\fnt{bf7}#1}}}

\def\rest#1{\raise-2pt\hbox{$|_{#1}$}}

\def\notdvd{\mathrel{\setbox\xbox=\hbox{$\big|$}%
\hbox to 0pt{\hbox to\wd\xbox{\hss/\hss}\hss}\big|}}

\def\frac#1#2{{#1\over #2}}
\def\matr#1{{\let\\=\cr\left(\matrix{#1\crcr}\right)}}
\def\vect#1{{\let\\=\cr\left(\matrix{#1\crcr}\right)}}
\def\smatr#1{{\baselineskip=2pt\lineskip=2pt
\left(\vcenter{\let\\=\cr\let\-=\hfill
\ialign{\hfil$\scr##$\hfil&&\hfil\kern2pt$\scr##$\hfil\cr#1\crcr}}\right)}}

\def\svd{\vbox to 2.4mm{\nolineskip
\kern.3pt\hbox{.}\vfil\hbox{.}\vfil\hbox{.}\kern.3pt}}
\def\vd{\vbox to 3.2mm{\nolineskip
\kern1pt\hbox{.}\vfil\hbox{.}\vfil\hbox{.}\kern1pt}}
\def\lvd{\vbox to 7mm{\nolineskip
\kern1pt\hbox{.}\vfil\hbox{.}\vfil\hbox{.}\vfil\hbox{.}
\vfil\hbox{.}\vfil\hbox{.}\kern1pt}}
\def\rvd{\vbox to 2.4mm{\nolineskip
\kern.3pt\hbox{.}\vfil\hbox{\kern3.5pt.}\vfil\hbox{\kern7pt.}\kern.3pt}}
\def\lld{\hbox to 5mm{\kern2pt.\hfil.\hfil.\kern2pt}}

\def\comp{\mathord{\hbox{$\scr\circ$}}}
\def\semprod{\mathrel{\times\kern-2pt
\vbox{\hrule width.5pt height4pt depth0pt\kern.5pt}\kern1pt}}

\def\dsc{\discretionary{}{}{}}
\long\def\omit#1\endomit{\par\vbox{%
\hrule\vskip2mm\hfil\vdots\hfil\vskip2mm\hrule}}
\long\def\ignore#1\endignore{}
\def\R{{\mb R}}

\def\Z{{\mb Z}}
\def\N{{\mb N}}
\def\Q{{\mb Q}}
\def\PP{\raise2.2pt\hbox{\fnt{mi.1}\char"7D}}

\def\const{\mathop{\hbox{\rm const}}}


\def\nolineskip{\baselineskip=0pt\lineskip=0pt}

\let\dsp=\displaystyle

\let\ovr=\overline
\let\und=\underline
\let\scr=\scriptstyle

\let\sln=\subset
\let\sle=\subseteq

\let\sm=\setminus

\let\col=\colon

\let\ld=\ldots
\let\vd=\vdots
\let\cd=\cdot
\let\pus=\emptyset
\let\ra=\longrightarrow
\let\ras=\rightarrow

\let\alf=\alpha
\let\bet=\beta

\let\del=\delta

\let\phi=\varphi
\let\sig=\sigma
\let\eps=\varepsilon
\let\lam=\lambda

\let\om=\omega

\let\afterstate=Y
\subsectionskip=5mm
\def\tsubsection#1#2{\subsection{#1} {\bf #2}\vbreak{9000}{1mm}}
\def\thestate{\the\sectionnum.\the\statenum}

\def\esslim{\mathop{\hbox{\rm ess-lim}}}
\def\D{\mathop{\hbox{\rm D}}}
\def\Du{\mathop{\und{\hbox{\rm D}}}}
\def\Do{\mathop{\ovr{\hbox{\rm D}}}}
\def\UD{\mathop{\hbox{\rm UD}}}
\def\UDu{\mathop{\und{\hbox{\rm UD}}}}
\def\UDo{\mathop{\ovr{\hbox{\rm UD}}}}
\def\ae{a.e.\nasp}
\def\Rp{\R_{+}}
\def\tri{\triangle}

\def\cH{{\cal H}}
\def\cT{{\cal T}}
\def\cP{{\cal P}}
\def\cG{{\cal G}}

\def\cL{{\cal L}}
\def\Sp{{\cal S}}
\let\Lam=\Lambda
\def\T{{\mb T}}
\def\tf{\tilde{f}}

\def\hf{\hat{f}}
\def\smatr#1{{\baselineskip=2pt\lineskip=2pt
\left(\vcenter{\let\\=\cr\let\-=\hfill
\ialign{\hfil$\scr##$\hfil&&\hfil\kern2pt$\scr##$\hfil\cr#1\crcr}}\right)}}

\def\svd{\vbox to 2.4mm{\nolineskip
\kern.3pt\hbox{.}\vfil\hbox{.}\vfil\hbox{.}\kern.3pt}}
\def\vd{\vbox to 3.2mm{\nolineskip
\kern1pt\hbox{.}\vfil\hbox{.}\vfil\hbox{.}\kern1pt}}
\def\lvd{\vbox to 7mm{\nolineskip
\kern1pt\hbox{.}\vfil\hbox{.}\vfil\hbox{.}\vfil\hbox{.}
\vfil\hbox{.}\vfil\hbox{.}\kern1pt}}
\def\rvd{\vbox to 2.4mm{\nolineskip
\kern.3pt\hbox{.}\vfil\hbox{\kern3.5pt.}\vfil\hbox{\kern7pt.}\kern.3pt}}
\def\lld{\hbox to 5mm{\kern2pt.\hfil.\hfil.\kern2pt}}
\start
\title{From discrete- to continuous-time ergodic theorems}
\kern3mm
\centerline{\it In memory of Dan Rudolph}
\author{V.~Bergelson, A.~Leibman, and C.~G.~Moreira}
\support{The first two authors are partially supported by NSF grant DMS-0901106.}
\titledate
\abstract{
We introduce methods that allow to derive 
continuous-time versions of various discrete-time ergodic theorems.
We then illustrate these methods
by giving simple proofs and refinements of some known results
as well as establishing new results of interest.}
\section{S-Int}{Introduction}

The goal of this paper is to introduce methods
that allow one to obtain continuous-time versions 
of various discrete-time ergodic results.
While the classical von~Neumann's and Birkhoff's ergodic theorems
were dealing with continuous families of invertible measure preserving transformations,
it was very soon observed that ergodic theorems for $\Z$-actions
hold true as well,
are somewhat easier to handle,
and, moreover, can be used as an auxiliary tool for the derivation
of the corresponding continuous-time results.
(See, for example, the formulation 
of the so-called Birkhoff's fundamental lemma in \brfr{BirKoo}.
See also \brfr{Kolmogorov} and \brfr{Hopf}, Section 8.)
Moreover, since not every measure preserving $\Z$-action
imbeds in a continuous measure preserving $\R$-flow,
and since there are various important classes 
of non-invertible measure-preserving transformations,
it became, over the years, more fashionable
to study ergodic theorems for measure preserving $\Z$- and $\N$-actions.
Numerous multiple recurrence and convergence results
obtained in the framework of the ergodic Ramsey theory also focused
(mainly due to combinatorial and number theoretical applications)
on $\Z$-actions and, more generally, actions of various discrete semigroups.

There are, however, questions in modern ergodic theory
pertaining to measure preserving $\R$-actions
that naturally present themselves and are connected with interesting applications
but do not seem to easily follow from the corresponding results for $\Z$-actions.
To better explain our point, let us consider some examples.
We start with the $\R$-version of the von~Neumann's ergodic theorem 
(\brfr{vonNeumann}):
if $T^{t}$, $t\in\R$, 
is an ergodic 1-parameter group of measure preserving transformations
of a probability measure space $(X,\mu)$,
then for any $f\in L^{2}(X)$, 
$\lim_{b-a\ras\infty}\frac{1}{b-a}\int_{a}^{b}T^{t}f\,dt=\int_{X}f\,d\mu$
\fnote{Here and below, $T^{t}f(\om)=f(T^{t}\om)$, $t\in\R$, $\om\in X$,
and the integral $\int_{a}^{b}T^{t}f\,dt$ is understood in the sense of Bochner.}
in $L^{2}(X)$.
An easy trick shows that this result immediately follows 
from the corresponding theorem for $\Z$-actions,
which says that if $T$ is an invertible ergodic measure preserving transformation
of a probability measure space $(X,\mu)$,
then for any $f\in L^{2}(X)$, 
$\lim_{N-M\ras\infty}\frac{1}{N-M}\sum_{n=M+1}^{N}T^{n}f=\int_{X}f\,d\mu$ in $L^{2}(X)$.
Indeed, all one has to do is to apply the $\Z$-version of von~Neumann's theorem
to the function $\tf=\int_{0}^{1}T^{t}f\,dt$ and the transformation $T^{1}$,
utilizing the fact that, for any $a,b\in\R$, 
$\int_{a}^{b}T^{t}f\,dt=\sum_{n=[a]}^{[b]-1}T^{n}\tf
-\int_{[a]}^{a}T^{t}f\,dt+\int_{[b]}^{b}T^{t}f\,dt$.
(The $\R$-version of Birkhoff's pointwise ergodic theorem
can be derived from its $\Z$-version in a similar way.)
This argument is no longer applicable to ``multiple ergodic averages''
\equn{f-Ilin}{
\frac{1}{b-a}\int_{a}^{b}T^{\alf_{1}t}f_{1}\cd\ld\cd T^{\alf_{r}t}f_{r}\,dt,
}
where $r\geq 2$, $\alf_{i}\in\R$, and $f_{i}\in L^{\infty}(X)$;
however, it can be modified so that one is still able to show 
that the averages \frfr{f-Ilin} converge in $L^{2}$-norm as $b-a\ra\infty$
as long as it is known that for arbitrarily small $u>0$ the averages
$\frac{1}{N-M}\sum_{n=M+1}^{N}T^{\alf_{1}un}f_{1}\cd\ld\cd T^{\alf_{r}un}f_{r}$
converge as $N-M\ra\infty$.
(See, for example, \brfr{Austin-l}.)
Indeed, given $\eps>0$,
find $\del>0$ such that $\|T^{\alf_{i}t}f-f\|<\eps$ for all $t\in(0,\del)$,
where $\|\cd\|=\|\cd\|_{L^{2}(X)}$;
then, assuming w.l.o.g.\nasp\ that $\sup|f_{i}|\leq 1$, $i=1,\ld,r$, 
we have, for any $t\in\R$, $\bigl\|T^{\alf_{i}t}f_{i}
-T^{\alf_{i}\del[t/\del]}f_{i}\bigr\|<\eps$, $i=1,\ld,r$,
and so $\bigl\|\prod_{i=1}^{r}T^{\alf_{i}t}f_{i}
-\prod_{i=1}^{r}T^{\alf_{i}\del[t/\del]}f_{i}\bigr\|<r\eps$.
Hence,
\lequ{
\limsup_{b-a\ras\infty}
\Bigl\|\frac{1}{b-a}\int_{a}^{b}\prod_{i=1}^{r}T^{\alf_{i}t}f_{i}\,dt-
\frac{1}{[b/\del]-[a/\del]}\sum_{n=[a/\del]}^{[b/\del]-1}
\prod_{i=1}^{r}T^{\alf_{i}\del n}f_{i}\Bigr\|
\-\\\-
\leq\limsup_{b-a\ras\infty}\Bigl(
\frac{1}{b-a}\sum_{n=[a/\del]}^{[b/\del]-1}\int_{n\del}^{(n+1)\del}
\Bigl\|\prod_{i=1}^{r}T^{\alf_{i}t}f_{i}-\prod_{i=1}^{r}T^{\alf_{i}\del n}f_{i}\Bigr\|\,dt
\-\\\-
+\frac{1}{b-a}\int_{[a/\del]\del}^{a}\Bigl\|\prod_{i=1}^{r}T^{\alf_{i}t}f_{i}\Bigr\|\,dt
+\frac{1}{b-a}\int_{[b/\del]\del}^{b}\Bigl\|\prod_{i=1}^{r}T^{\alf_{i}t}f_{i}\Bigr\|\,dt\Bigr)
\leq r\eps.
}
Since 
$\lim_{N-M\ras\infty}\frac{1}{N-M}\sum_{n=M+1}^{N}\prod_{i=1}^{r}T^{\alf_{i}\del n}f_{i}$
exists in $L^{2}(X)$, 
and since $\eps$ is arbitrary, we get that
$\lim_{b-a\ras\infty}\frac{1}{b-a}\int_{a}^{b}\prod_{i=1}^{r}T^{\alf_{i}t}f_{i}\,dt$
also exists.

However, even this argument stops working if we consider, say, 
the ``polynomial averages''
$\frac{1}{b-a}\int_{a}^{b}T^{p(t)}f\,dt$,
where $p$ is a polynomial,
or, more generally, the ``polynomial multiple averages''
\equn{f-mpe}{
\frac{1}{b-a}\int_{a}^{b}
T^{p_{1}(t)}f_{1}\cd\ld\cd T^{p_{r}(t)}f_{r}\,dt,
}
where $p_{i}$ are polynomials,
since in this case the function $\phi(t)=T^{p(t)}f$ from $\R$ to $L^{2}(X)$
is no longer uniformly continuous.
The convergence of the corresponding discrete-time averages
\equn{f-mped}{
\frac{1}{N-M}\sum_{n=M+1}^{N}
T^{p_{1}(n)}f_{1}\cd\ld\cd T^{p_{r}(n)}f_{r},
}
is known (see \brfr{HKp} and \brfr{cpm}),
but to establish the convergence in $L^{2}(X)$ of the averages \frfr{f-mpe}
one either has to go through all the main stages 
of the proof of the convegrence of averages \frfr{f-mped}
and verify the validity of the corresponding $\R$-statements
(see, for example, \brfr{Potts}, 
where the existence of the non-uniform limits
$\frac{1}{b}\int_{0}^{b}T^{p_{1}(t)}f_{1}\cd\ld\cd T^{p_{r}(t)}f_{r}dt$
is established),
or may try to find some alternative general method 
connecting the convergence of discrete- and of continuous-time averages.
(Yet another approach to proving convergence of multiple polynomial averages,
utilized in \brfr{Austin-cp},
is based on a ``change of variables'' trick and usage of equivalent methods of summation;
this method allows one to treat expressions like
$\frac{1}{b}\int_{0}^{b}\prod_{j=1}^{k}T_{j}^{p_{j,1}(t)}f_{1}\cd\ld
\cd\prod_{j=1}^{k}T_{j}^{p_{j,r}(t)}f_{r}\,dt$,
where $T_{j}$ are commuting measure preserving transformations.
However, this method gives no information about what the limits of such averages are,
and, also, it is not clear whether it can be extended 
to obtain convergence of uniform averages \frfr{f-mpe}.)

As another example where a passage from discrete to continuous setup 
is desirable but not apriori obvious
let us mention the problem of the study 
of the distribution of values of generalized polynomials.
A generalized polynomial is a function that is obtained from conventional polynomials
of one or several variables
by applying the operations of taking the integer part, addition and multiplication;
for example, if $p_{i}(x)$ are conventional polynomials,
then $u(x)=\bigl[[p_{1}(x)]p_{2}(x)+p_{3}(x)\bigr]p_{4}(x)
+\bigl[p_{5}(x)[p_{6}(x)]\bigr]^{2}p_{7}(x)$ 
is a generalized polynomial.
It was shown in \brfr{sko}
that the values of any bounded vector-valued generalized polynomial of integer argument
are well distributed on a piecewise polynomial surface,
with respect to a natural measure on this surface.
The proof was based on the theorem on well-distribution
of polynomial orbits on nilmanifolds;
since such a theorem for continuous polynomial flows on nilmanifolds 
was not known at the time of writing \brfr{sko},
we could not prove that bounded generalized polynomials of continuous argument
are well distributed on piecewise polynomial surfaces.
(See \brfr{sko}, Theorem~B$_{c}$.
As a matter of fact, 
the problem of extending the results from \brfr{sko} to the case of continuous parameter
served as an impetus for the present paper.)

In this paper we introduce two simple but quite general methods
that allow one to deduce continuous-time ergodic theorems
from their discrete-time couterparts.
To deliver the zest of these methods 
we will formulate now two easy to state theorems.
Let $F(t)$ be a bounded measurable function from $[0,\infty)$ to a Banach space.
(In our applications, 
$F$ will usually be ``an ergodic expression'' 
that depends on a continuous parameter $t$ and takes values in a functional space,
say $F(t)=T_{1}^{t}f_{1}\cd\ld\cd T_{k}^{t}f_{k}\in L^{1}(X)$, $t\in[0,\infty)$,
where $T_{i}$ are 1-parameter groups of measure preserving transformations 
of a measure space $X$
and $f_{i}\in L^{\infty}(X)$).

\proposition{P-addInt}{}{(Additive method)} 
If the limit $\lim_{N\ras\infty}\frac{1}{N}\sum_{n=0}^{N-1}F(t+n)=A_{t}$ exists
for a.e.\nasp\ $t\in[0,1)$
then the limit $\lim_{N\ras\infty}\frac{1}{b}\int_{0}^{b}F(t)\,dt$ also exists,
and equals $\int_{0}^{1}A_{t}\,dt$.
\endproposition

\proposition{P-mulInt}{}{(Multiplicative method)} 
If the limit $\lim_{N\ras\infty}\frac{1}{N}\sum_{n=0}^{N-1}F(nt)=L_{t}$ exists
for a.e.\nasp\ $t\in(0,1)$,
then the limit $L=\lim_{N\ras\infty}\frac{1}{b}\int_{0}^{b}F(t)\,dt$ also exists,
and, moreover, $L_{t}=L$ for \ae\ $t\in(0,1]$.
\endproposition

Each of these ``methods'' has its pros and cons.
The ``additive'' method is very easy to substantiate.
However, it has the disadvatage that, being non-homogeneous,
it ``desinchronizes'' the expression $F(t)$,
which may be an obstacle for certain applications.
Consider, for example, the expression $F(t)=T_{1}^{t}f\cd\ld\cd T_{k}^{t}f$,
appearing in the formulation of the $\R$-version of the ergodic Szemer\'{e}di theorem 
(see \rfr{s-szem} below).
In this case $F(t+n)=T_{1}^{n}f_{1}\cd\ld\cd T_{k}^{n}f_{k}$,
where $f_{i}=T_{i}^{t}f$, $i=1,\ld,k$ are, generally speaking, distinct functions,
which complicates application of the ``discrete'' ergodic Szemer\'{e}di theorem.
The ``multiplicative method'' is quite a bit harder to establish,
but it preserves the ``structure'' of $F(t)$:
for $F(t)=T_{1}^{t}f\cd\ld\cd T_{k}^{t}f$
we now have $F(nt)=(T_{1}^{t})^{n}f\cd\ld\cd(T_{k}^{t})^{n}f$.
An additional advantage of the multiplicative method
is that it guarantees the equality of almost all ``discrete'' limits $L_{t}$,
and therefore gives more information about the ``continuous'' limit $L$.
(See \rfr{P-poleqlim} below.)

When it comes to convergence on average, 
there are many types of it (uniform, strong Ces\`{a}ro, etc.)
which naturally appear in various situations 
in classical analysis, number theory, and ergodic theory,
and for each of them one can provide a statement
that connects discrete and continuous averages.
We therefore present several similar results;
their proofs are based on similar ideas,
but utilizing these ideas in diverse situations
we obtain a variety of useful theorems.
Here is the descriptive list of various kinds of averaging schemes 
we will be dealing with.
(In what follows $V$ stands for an abstract Banach space.)

\npar$\bullet$
{\it One-parameter standard Ces\`{a}ro limits:}
The Ces\`{a}ro limit of a sequence $(v_{n})$ in $V$ is
$\lim_{N\ras\infty}\frac{1}{N}\sum_{n=1}^{N}v_{n}$,
and for a measurable function $f\col[0,\infty)\ra V$ it is
$\lim_{b\ras\infty}\frac{1}{b}\int_{0}^{b}f(x)\,dx$.

\npar$\bullet$
{\it One-parameter uniform Ces\`{a}ro limits:}
The uniform Ces\`{a}ro limit for a sequence $(v_{n})$ in $V$ is
$\lim_{N-M\ras\infty}\frac{1}{N-M}\dsc\sum_{n=M+1}^{N}v_{n}$,
and for a measurable function $f\col[0,\infty)\ra V$ it is
$\lim_{b-a\ras\infty}\frac{1}{b-a}\int_{a}^{b}f(x)\,dx$.

\npar
(The ``one-parameter averaging schemes'' above
is, of course, a special case of the corresponding ``multiparameter schemes'' below,
but we start with the one-parameter case to make our proofs more transparent.)

An $\N^{d}$-sequence $(v_{n})$ in $V$
is a mapping $\N^{d}\ra V$, $n\mapsto v_{n}$.
For a parallelepiped $P=\prod_{i=1}^{d}[a_{i},b_{i}]\sln\R^{d}$
we define $l(P)=\min_{1\leq i\leq d}(b_{i}-a_{i})$
and $w(P)=\prod_{i=1}^{d}(b_{i}-a_{i})$.

\npar$\bullet$
{\it Multiparameter standard Ces\`{a}ro limits:}
The Ces\`{a}ro limit of an $\N^{d}$-sequence $(v_{n})$ in $V$ is
$\lim_{l(P)\ras\infty}\frac{1}{w(P)}\dsc\sum_{n\in\N^{d}\cap P}v_{n}$,
and for a measurable function $f\col[0,\infty)^{d}\ra V$ it is
$\lim_{l(P)\ras\infty}\frac{1}{w(P)}\int_{P}f(x)\,dx$,
where, in both cases, $P$ runs over the set of parallelepipeds
of the form $\prod_{i=1}^{d}[0,b_{i}]$ in $[0,\infty)^{d}$.

\npar$\bullet$
{\it Multiparameter uniform Ces\`{a}ro limits:}
The uniform Ces\`{a}ro limit 
of an $\N^{d}$-sequence $(v_{n})$ in $V$ is
$\lim_{l(P)\ras\infty}\frac{1}{w(P)}\dsc\sum_{n\in\Z^{d}\cap P}v_{n}$,
and for a measurable function $f\col[0,\infty)^{d}\ra V$ it is
$\lim_{l(P)\ras\infty}\frac{1}{w(P)}\int_{P}f(x)\,dx$,
where, in both cases, $P$ runs over the set of parallelepipeds
of the form $\prod_{i=1}^{d}[a_{i},b_{i}]$ in $[0,\infty)^{d}$.

\npar$\bullet$
{\it Two-sided {\rm(or, rather, {\it all sided})}
standard and uniform Ces\`{a}ro limits:}
Instead of $\N^{d}$-sequences, functions on $[0,\infty)^{d}$,
and parallelepipeds in $[0,\infty)^{d}$,
we deal with $\Z^{d}$-sequences, functions on $\R^{d}$, 
and parallelepipeds in $\R^{d}$.

\npar$\bullet$
{\it Limits of averages along general F{\o}lner sequences:}
Instead of averaging over parallelepipeds of the form $P=\prod_{i=1}^{d}[a_{i},b_{i}]$,
we consider averages over elements of a general F{\o}lner sequence
$(\Phi_{N})_{N=1}^{\infty}$ in $\R^{d}$,
\hbox{$\lim_{N\ras\infty}\frac{1}{w(\Phi_{N})}\int_{\Phi_{N}}f(x)\,dx$},
where $w$ stands for the Lebesgue measure on $\R^{d}$.

\npar$\bullet$
{\it Liminf and limsup versions for standard and uniform averages:}
When the limits above do not (or are not known to) exist,
but $(v_{n})$ is a real-valued sequence and $f$ is a real-valued function,
we consider the corresponding liminfs and limsups.

\npar$\bullet$
{\it Lim-limsup versions:}
If the limits $\lim_{N\ras\infty}\frac{1}{N}\sum_{n=0}^{N-1}F(nt)$ 
do not, or are not known to exist, 
it may still be possible that for some $L\in V$,
$\lim_{t\ras 0^{+}}\limsup_{N\ras\infty}\bigl\|\frac{1}{N}
\sum_{n=0}^{N-1}F(nt)-L\bigr\|=0$;
it turns out that this suffices for the multiplicative method to work.

After proving several versions of \rfrt{P-addInt}s~\rfrn{P-addInt} and \rfrn{P-mulInt}
corresponding to different averaging schemes,
we will apply them to re-prove some known 
and establish some new ergodic-theoretical results;
here is a list of the applications that we obtain in \rfr{S-Appl}:

\vbreakn{8000}{1mm}$\circ$
In section~\rfrn{s-charfac},
we show that characteristic factors for averages of the form
\equn{f-Impd}{
\frac{1}{w(\Phi_{N})}\int_{\Phi_{N}}T^{p_{1}(t)}f_{1}\cd\ld\cd T^{p_{r}(t)}f_{r}\,dt,
}
where $T^{t}$, $t\in\R$, is a continuous $1$-parameter group
of measure preserving transformations of a probability measure space $X$,
$p_{i}$ are polynomials $\R^{d}\ra\R$, $f_{i}\in L^{\infty}(X)$,
and $(\Phi_{N})$ is a F{\o}lner sequence in $\R^{d}$,
are Host-Kra-Ziegler factors of $X$.
(A non-uniform version of this result is obtained in \brfr{Potts}.)

\vbreakn{-9000}{1mm}$\circ$
In section~\rfrn{s-nilorb},
we prove that, for any $d\in\N$, 
a $d$-parameter polynomial flow on a nilmanifold $X$
is well distributed on a subnilmanifold of $X$.
(This result is new and refines the fact
that any such flow is uniformly distributed in a subnilmanifold of $X$.)

\vbreakn{-9000}{1mm}$\circ$
In section~\rfrn{s-pollim},
we prove the convergence of averages \frfr{f-Impd}.
(This result is new,
strengthening the results obtained in \brfr{Potts} 
and the one-parameter case of the results obtained in \brfr{Austin-cp}.) 
We also prove that the averages
$\frac{1}{b-a}\int_{a}^{b}T_{1}^{t}f_{1}\cd\ld\cd T_{r}^{t}\,dt$ converge,
where $T_{i}$ are pairwise commuting measure preserving transformations.
(This strengthens the linear case of the results obrained in \brfr{Austin-cp}.)

\vbreakn{-9000}{1mm}$\circ$
In section~\rfrn{s-szem},
we obtain a continuous-time version of the polynomial ergodic Szemer\'{e}di theorem.

\vbreakn{-9000}{1mm}$\circ$
In section~\rfrn{s-genpol},
we prove that the values of bounded vector-valued generalized polynomials 
are well-distributed on a piecewise polynomial surface.
This establishes the continuous version of the well-distribution result from \brfr{sko}
that we discussed above.

\vbreakn{-9000}{1mm}$\circ$
In section~\rfrn{s-hardywm},
we derive, from the corresponding discrete-time results in \brfr{Nikos2} and \brfr{Inger},
convergence of multiple averages \frfr{f-mpe}
with $p_{i}$ being functions of polynomial growth.

\vbreakn{-9000}{1mm}$\circ$
Finally, in section~\rfrn{s-pointwise},
we apply our methods to obtain continuous-time theorems
dealing with almost everywhere convergence of certain ergodic averages.
\acknowledgment
We thank the referee for many useful comments and corrections.
\section{S-dct}{A Fatou lemma and a dominated convergence theorem}
Throughout Sections~\rfrn{S-dct} -- \rfrn{S-Density},
$V$ stands for a separable Banach space.
We will repeatedly use the following Fatou-like lemma and its corollary:
\lemma{P-dctsup}{}{}
Let $(X,\mu)$ be a finite measure space
and let $(f_{n})$ be a sequence of uniformly bounded measurable functions 
from $X$ to $V$.
Then $\limsup_{n\ras\infty}\bigl\|\int_{X}f_{n}\,d\mu\bigr\|
\leq\int_{X}\limsup_{n\ras\infty}\|f_{n}\|\,d\mu$.
\endlemma
\proof{}
Let $M>0$ be such that $\|f_{n}(x)\|\leq M$ for all $x\in X$ and $n\in\N$,
and let $s(x)=\limsup_{n\ras\infty}\|f_{n}(x)\|$, $x\in X$.
Fix $\eps>0$.
For each $x\in X$ let $n(x)\in\N$ be such that 
$\|f_{n}(x)\|<s(x)+\eps$ for all $n\geq n(x)$.
For each $n\in\N$, let $A_{n}=\bigl\{x\in X:n(x)\leq n\bigr\}$.
Then $A_{1}\sle A_{2}\sle\ld$ and $\bigcup_{n=1}^{\infty}A_{n}=X$,
so $\lim_{n\ras\infty}\mu(X\sm A_{n})=0$.
Let $N$ be such that $\mu(X\sm A_{N})<\eps$.
Then for any $n\geq N$,
\lequ{
\Bigl\|\int_{X}f_{n}\,d\mu\Bigr\|
\leq\int_{X}\|f_{n}\|\,d\mu
=\int_{A_{N}}\|f_{n}\|\,d\mu+\int_{X\sm A_{N}}\kern-3mm\|f_{n}\|\,d\mu
\leq\int_{A_{N}}(s+\eps)\,d\mu+M\eps
\leq\int_{X}s\,d\mu+\eps(\mu(X)+M).
}
Since this is true for any positive $\eps$,
$\bigl\|\int_{X}f_{n}\,d\mu\bigr\|\leq\int_{X}s\,d\mu$.%
\endproof

As a corollary, we get:
\lemma{P-dct}{}{}
Let $(X,\mu)$ be a finite measure space.
If a sequence $(f_{n})$ of uniformly bounded measurable functions 
from $X$ to $V$
converges to a function $f\col X\ra V$ a.e.\nasp\ on $X$,
then $\int_{X}f_{n}\,d\mu\ra\int_{X}f\,d\mu$.
\endlemma

\remark{}
Of course, a more general dominated convergence theorem, 
where $\|f_{n}\|$ are not assumed to be bounded 
but only dominated by an integrable function,
like in the case of real-valued functions, also holds,
but we will only need its special case given by \rfr{P-dct}.
\endremark
\section{S-add}{Additive method}

When $a$ and $b$ are positive real numbers,
we define $\sum_{n>a}^{b}v_{n}
=\cases{\sum_{n\in(a,b]\cap\N}v_{n}\hbox{ if $a<b$}\cr 0\hbox{ if $a\geq b$.}\cr}$
\tsubsection{s-addC}{Standard Ces\`{a}ro limits}
\theorem{P-addC}{}{}
Let $f\col[0,\infty)\ra V$ be a bounded measurable function
such that the limit 
$$
A_{t}=\lim_{b\ras\infty}\frac{1}{b}\sum_{n>0}^{b}f(t+n)
$$ 
exists for \ae\ $t\in[0,1]$.
Then $\lim_{b\ras\infty}\frac{1}{b}\int_{0}^{b}f(x)\,dx$ also exists
and is equal to $\int_{0}^{1}A_{t}\,dt$.
\endtheorem
\proof{}
We may assume that the parameter $b$ is integer.
For any $b\in\N$ we have 
\equ{
\frac{1}{b}\int_{0}^{b}f(x)\,dx
=\frac{1}{b}\sum_{n=0}^{b-1}\int_{0}^{1}f(t+n)\,dt
=\int_{0}^{1}\frac{1}{b}\sum_{n=0}^{b-1}f(t+n)\,dt.
}
Since for \ae\ $t\in[0,1]$,
$\frac{1}{b}\sum_{n=0}^{b-1}f(t+n)\,dt\ra A_{t}$ as $b\ras\infty$,
by \rfr{P-dct},
\equ{
\lim_{b\ras\infty}\frac{1}{b}\int_{0}^{b}f(x)\,dx
=\int_{0}^{1}A_{t}\,dt.
}
\frgdsp\endproof

\remark{}
Of course, in the formulation of \rfr{P-addC} 
the interval $[0,1]$ and the expression $f(t+n)$
can be replaced by the interval $[0,\del]$ and the expression $f(t+n\del)$
for any positive $\del$.
\endremark
\tsubsection{s-addU}{Uniform Ces\`{a}ro limits}
\theorem{P-addU}{}{}
Let $f\col[0,\infty)\ra V$ be a bounded measurable function
such that the limit 
$$
A_{t}=\lim_{b-a\ras\infty}\frac{1}{b-a}\sum_{n>a}^{b}f(t+n)
$$ 
exists for \ae\ $t\in[0,1]$.
Then $\lim_{b-a\ras\infty}\frac{1}{b-a}\int_{a}^{b}f(x)\,dx$ also exists
and is equal to $\int_{0}^{1}A_{t}\,dt$.
\endtheorem
\proof{}
We may assume that the parameters $a$, $b$ are integer.
For any sequences $(a_{k})$, $(b_{k})$ of nonnegative integers with $b_{k}-a_{k}\ra+\infty$
we have
\equ{
\frac{1}{b_{k}-a_{k}}\int_{a_{k}}^{b_{k}}f(x)\,dx
=\frac{1}{b_{k}-a_{k}}\sum_{n=a_{k}}^{b_{k}-1}\int_{0}^{1}f(t+n)\,dt
=\int_{0}^{1}\frac{1}{b_{k}-a_{k}}\sum_{n=a_{k}}^{b_{k}-1}f(t+n)\,dt.
}
Since for \ae\ $t\in[0,1]$,
$\frac{1}{b_{k}-a_{k}}\sum_{n=a_{k}}^{b_{k}-1}f(t+n)\,dt\ra A_{t}$ as $k\ras\infty$,
by \rfr{P-dct},
\equ{
\lim_{k\ras\infty}\frac{1}{b_{k}-a_{k}}\int_{a_{k}}^{b_{k}}f(x)\,dx
=\int_{0}^{1}A_{t}\,dt.
}
\frgdsp\endproof
\tsubsection{s-addMC}{Multiparameter standard Ces\`{a}ro limits}
Let $d\in\N$.
We will call a mapping $\N^{d}\ra V$, $n\mapsto v_{n}$
{\it $\N^{d}$-sequence in $V$}.
We write $\Rp$ for $[0,\infty)$.
We will now introduce notation
that will allow us to formulate and prove the $d$-parameter versions of the above theorems 
in complete analogy with the case $d=1$.

For $a,b\in\Rp^{d}$, $a=(a_{1},\ld,a_{d})$, $b=(b_{1},\ld,b_{d})$, 
we write $a\leq b$ if $a_{i}\leq b_{i}$ for all $i=1,\ld,d$
and $a<b$ if $a_{i}<b_{i}$ for all $i$.
Under $\min(a,b)$ and $\max(a,b)$ 
we will understand $\bigl(\min(a_{1},b_{1}),\ld,\min(a_{d},b_{d})\bigr)$
and $\bigl(\max(a_{1},b_{1}),\ld,\max(a_{d},b_{d})\bigr)$ respectively.
For $a=(a_{1},\ld,a_{d})\in\Rp^{d}$ and $b=(b_{1},\ld,b_{d})\in\Rp^{d}$,
we define $ab=(a_{1}b_{1},\ld,a_{d}b_{d})$,
and if $b>0$, $a/b=(a_{1}/b_{1},\ld,a_{d}/b_{d})$,
and $b^{\alf}=(b_{1}^{\alf},\ld,b_{d}^{\alf})$, $\alf\in\R$.

For $a=(a_{1},\ld,a_{d})\in\Rp^{d}$ 
we define $w(a)=a_{1}\cdots a_{d}$ and $l(a)=\min\{a_{1},\ld,a_{d}\}$.
Note that if $a,b\in\Rp^{d}$ and $0<a\leq b$, then $w(a)/w(b)\leq l(a)/l(b)$.

For $a,b\in\Rp^{d}$, $a\leq b$, 
we define {\it intervals\/} $[a,b]=\bigl\{x\in\Rp^{d}:a\leq x\leq b\bigr\}$
and $(a,b]=\bigl\{x\in\Rp^{d}:a<x\leq b\bigr\}$.

For $a,b\in\Rp^{d}$,
under $\sum_{n=a}^{b}v_{n}$ we will understand
$\sum_{n\in\N^{d}\cap[a,b]}v_{n}$ if $a\leq b$ and 0 otherwise,
under $\sum_{n>a}^{b}v_{n}$ we will understand
$\sum_{n\in\N^{d}\cap(a,b]}v_{n}$ if $a\leq b$ and 0 otherwise,
and under $\int_{a}^{b}v(x)\,dx$ we will understand $\int_{[a,b]}v(x)\,dx$.

Finally, for $c\in\Rp$, by $\bar{c}$ we will denote $(c,\ld,c)\in\Rp^{d}$.

\theorem{P-addMC}{}{}
Let $f\col\Rp^{d}\ra V$ be a bounded measurable function
such that the limit 
$$
A_{t}=\lim_{l(b)\ras\infty}\frac{1}{w(b)}\sum_{n>0}^{b}f(t+n)
$$ 
exists for \ae\ $t\in[0,1]^{d}$.
Then $\lim_{l(b)\ras\infty}\frac{1}{w(b)}\int_{0}^{b}f(x)\,dx$ also exists
and is equal to $\int_{[0,1]^{d}}A_{t}\,dt$.
\endtheorem
\proof{}
We may assume that $b\in\N^{d}$.
Let $(b_{k})$ be a sequence in $\N^{d}$ with $l(b_{k})\ras\infty$ as $k\ras\infty$.
For any $k\in\N$ we have
\equ{
\frac{1}{w(b_{k})}\int_{0}^{b_{k}}f(x)\,dx
=\frac{1}{w(b_{k})}\sum_{n=0}^{b_{k}-\bar{1}}\int_{[0,1]^{d}}f(t+n)\,dt
=\int_{[0,1]^{d}}\frac{1}{w(b_{k})}\sum_{n=0}^{b_{k}-\bar{1}}f(t+n)\,dt.
}
Since for \ae\ $t\in[0,1]^{d}$,
$\frac{1}{w(b_{k})}\sum_{n=0}^{b_{k}-\bar{1}}f(t+n)\,dt\ra A_{t}$ as $k\ras\infty$,
by \rfr{P-dct},
\equ{
\lim_{k\ras\infty}\frac{1}{w(b_{k})}\int_{0}^{b_{k}}f(x)\,dx
=\int_{[0,1]^{d}}A_{t}\,dt.
}
\frgdsp\endproof
\tsubsection{s-addMU}{Multiparameter uniform Ces\`{a}ro limits}
\theorem{P-addMU}{}{}
Let $f\col\Rp^{d}\ra V$ be a bounded measurable function
such that the limit 
$$
A_{t}=\lim_{l(b-a)\ras\infty}\frac{1}{w(b-a)}\sum_{n>a}^{b}f(t+n)
$$ 
exists for \ae\ $t\in[0,1]^{d}$.
Then $\lim_{l(b-a)\ras\infty}\frac{1}{w(b-a)}\int_{a}^{b}f(x)\,dx$ also exists
and is equal to $\int_{[0,1]^{d}}A_{t}\,dt$.
\endtheorem
\proof{}
We may assume that $a,b\in\N^{d}$.
Let $(a_{k})$, $(b_{k})$ be sequences in $\N^{d}$ 
with $a_{k}<b_{k}$ and $l(b_{k}-a_{k})\ras\infty$ as $k\ras\infty$.
For any $k\in\N$ we have
\equ{
\frac{1}{w(b_{k}-a_{k})}\int_{a_{k}}^{b_{k}}f(x)\,dx
=\frac{1}{w(b_{k}-a_{k})}\sum_{n=a_{k}}^{b_{k}-\bar{1}}\int_{[0,1]^{d}}f(t+n)\,dt
=\int_{[0,1]^{d}}\frac{1}{w(b_{k}-a_{k})}\sum_{n=a_{k}}^{b_{k}-\bar{1}}f(t+n)\,dt.
}
Since for \ae\ $t\in[0,1]^{d}$,
$\frac{1}{w(b_{k}-a_{k})}\sum_{n=0}^{b_{k}-\bar{1}}f(t+n)\,dt\ra A_{t}$ as $k\ras\infty$,
by \rfr{P-dct},
\equ{
\lim_{k\ras\infty}\frac{1}{w(b_{k}-a_{k})}\int_{a_{k}}^{b_{k}}f(x)\,dx
=\int_{[0,1]^{d}}A_{t}\,dt.
}
\frgdsp\endproof
\tsubsection{s-addMSC}{Liminf and limsup versions}
In the case $f$ is a real-valued function
we may obtain similar results involving liminfs of limsups,
even if the limits $A_{t}$ do not exist:
\theorem{P-addMSC}{}{}
If $f\col\Rp^{d}\ra\R$ is a bounded measurable function,
then
$$
\liminf_{l(b)\ras\infty}\frac{1}{w(b)}\int_{0}^{b}f(x)\,dx
\geq\int_{[0,1]^{d}}\liminf_{l(b)\ras\infty}\frac{1}{w(b)}\sum_{n>0}^{b}f(t+n)\,dt.
$$
and
$$
\limsup_{l(b)\ras\infty}\frac{1}{w(b)}\int_{0}^{b}f(x)\,dx
\leq\int_{[0,1]^{d}}\limsup_{l(b)\ras\infty}\frac{1}{w(b)}\sum_{n>0}^{b}f(t+n)\,dt.
$$ 
\endtheorem
\proof{}
After adding a constant to $f$ we may assume that $f\geq 0$.
We may also assume that $b\in\N^{d}$.
Let $(b_{k})$ be a sequence in $\N^{d}$ with $l(b_{k})\ras\infty$ as $k\ras\infty$.
For any $k\in\N$ we have
\equ{
\frac{1}{w(b_{k})}\int_{0}^{b_{k}}f(x)\,dx
=\frac{1}{w(b_{k})}\sum_{n=0}^{b_{k}-\bar{1}}\int_{[0,1]^{d}}f(t+n)\,dt
=\int_{[0,1]^{d}}\frac{1}{w(b_{k})}\sum_{n=0}^{b_{k}-\bar{1}}f(t+n)\,dt.
}
By (the classical, real-valued) Fatou's theorem,
\equ{
\liminf_{k\ras\infty}\frac{1}{w(b_{k})}\int_{0}^{b_{k}}f(x)\,dx
\geq\int_{[0,1]^{d}}\liminf_{k\ras\infty}
\frac{1}{w(b_{k})}\sum_{n=0}^{b_{k}-\bar{1}}f(t+n)\,dt
=\int_{[0,1]^{d}}\liminf_{k\ras\infty}\frac{1}{w(b_{k})}\sum_{n>0}^{b_{k}}f(t+n)\,dt.
}
\frgdsp\endproof

And similarly,
\theorem{P-addMSU}{}{}
If $f\col\Rp^{d}\ra\R$ is a bounded measurable function,
then
$$
\liminf_{l(b-a)\ras\infty}\frac{1}{w(b-a)}\int_{a}^{b}f(x)\,dx
\geq\int_{[0,1]^{d}}\liminf_{l(b-a)\ras\infty}\frac{1}{w(b-a)}\sum_{n>a}^{b}f(t+n)\,dt.
$$
and
$$
\limsup_{l(b-a)\ras\infty}\frac{1}{w(b-a)}\int_{a}^{b}f(x)\,dx
\leq\int_{[0,1]^{d}}\limsup_{l(b-a)\ras\infty}\frac{1}{w(b-a)}\sum_{n>a}^{b}f(t+n)\,dt.
$$ 
\endtheorem
\section{S-mul}{Multiplicative method -- the one-parameter case}
\tsubsection{s-mulC}{Standard Ces\`{a}ro limits}
\theorem{P-mulC}{}{}
Let a bounded measurable function $f\col[0,\infty)\ra V$ be such that for some $c>0$,
the limit $L_{t}=\lim_{b\ras\infty}\frac{1}{b}\sum_{n>0}^{b}f(nt)$ exists for \ae\ $t\in[0,c]$.
Then $L_{t}$ is \ae\ constant, $L_{t}=L\in V$ \ae\ on $[0,c]$,
and $\lim_{b\ras\infty}\frac{1}{b}\int_{0}^{b}f(x)\,dx$ exists and equals $L$.
\endtheorem

Following the referee's suggestion we will derive \rfr{P-mulC} from the following classical fact.
The proof of this result which we provide for reader's convenience
has an advantage of being easily extendible to multiparameter case
(\rfr{P-mulTaub} below).
\lemma{P-Taub}{}{}
Let $(v_{n})$ be a sequence in $V$
such that $\|v_{n+1}-v_{n}\|=O(1/n)$ and the Ces\'{a}ro limit
$L=\lim_{N\ras\infty}\frac{1}{N}\sum_{n=1}^{N}v_{n}$ exists.
Then $\lim_{n\ras\infty}v_{n}=L$.
\endlemma
\proof{}
We may assume that $L=0$, that is, $\lim_{N\ras\infty}\frac{1}{N}\sum_{n=1}^{N}v_{n}=0$.
Assume that $v_{n}\not\ras 0$;
let $\eps>0$ be such that for any $N\in\N$ there exists $n>N$ such that $\|v_{n}\|>\eps$.
Let $\alf>0$ be such that $\|v_{n+1}-v_{n}\|<\alf/n$ for all $n$;
put $\del=\frac{\eps^{2}}{16(1+\eps/2\alf)\alf}$.
Find $N\in\N$ such that $\bigl\|\frac{1}{M}\sum_{n=1}^{M}v_{n}\bigr\|<\del$ for all $M>N$.
Find $M>N$ such that $\|v_{M}\|>\eps$ and $1/M<\eps/4\alf$.
Let, by the Hahn-Banach theorem, $\phi\in V^{*}$ 
be such that $|\phi(v)|\leq\|v\|$ for all $v\in V$ and $\phi(v_{M})=\|v_{M}\|$.
Then for any $n>M$,
$$
\phi(v_{n})
=\phi(v_{M})+\sum_{m=M}^{n-1}\phi(v_{m+1}-v_{m})
\geq\|v_{M}\|-\sum_{m=M}^{n-1}\|v_{m+1}-v_{m}\|
>\eps-(n-M)\frac{\alf}{M},
$$
which is $\geq\eps/2$ when $\frac{(n-M)\alf}{M}\leq\frac{\eps}{2}$,
that is, when $M<n\leq M+\eps M/2\alf$.
Put $K=M+\lfloor\eps M/2\alf\rfloor$,
then $\phi(v_{n})>\frac{\eps}{2}$ for $n=M+1,\ld,K$.
Thus
\equ{
\Bigl\|\frac{1}{K}\sum_{n=1}^{K}v_{n}\Bigr\|
\geq\frac{1}{K}\Bigl(\Bigl\|\sum_{n=M+1}^{K}v_{n}\Bigr\|
-\Bigl\|\sum_{n=1}^{M}v_{n}\Bigr\|\Bigr)
\geq\frac{1}{K}\Bigl(\phi\Bigl(\sum_{n=M+1}^{K}v_{n}\Bigr)-\del M\Bigr)
=\frac{1}{K}\Bigl(\sum_{n=M+1}^{K}\phi(v_{n})-\del M\Bigr)
\-\\
>\frac{1}{K}\Bigl((K-M)\frac{\eps}{2}-\del M\Bigr)
\geq\frac{1}{M(1+\eps/2\alf)}
\Bigl(\Bigl(\frac{\eps M}{2\alf}-1\Bigr)\frac{\eps}{2}-\del M\Bigr)
\geq\frac{1}{1+\eps/2\alf}
\Bigl(\Bigl(\frac{\eps}{2\alf}-\frac{1}{M}\Bigr)\frac{\eps}{2}-\del\Bigr)
\\\-
\geq\frac{\eps^{2}}{(1+\eps/2\alf)8\alf}-\del=\del,
}
which contradicts the choice of $N$.
\endproof

\proof{of \rfr{P-mulC}}
Let $v_{n}=\frac{1}{nc}\int_{0}^{nc}f(x)\,dx$, $n\in\N$.
Then for any $n$,
$$
\|v_{n+1}-v_{n}\|
=\frac{1}{c}\Bigl\|\Bigl(\frac{1}{n+1}-\frac{1}{n}\Bigr)\int_{0}^{nc}f(x)\,dx
+\frac{1}{n+1}\int_{nc}^{(n+1)c}f(x)\,dx\Bigr\|
\leq\frac{cn\sup\|f\|}{cn(n+1)}+\frac{c\sup\|f\|}{c(n+1)}=O(1/n).
$$
Also we have, for any $N\in\N$,
$$
\frac{1}{N}\sum_{n=0}^{N-1}v_{n}
=\frac{1}{N}\sum_{n=0}^{N-1}\frac{1}{nc}\int_{0}^{nc}f(x)\,dx
=\frac{1}{N}\sum_{n=0}^{N-1}\frac{1}{c}\int_{0}^{c}f(nt)\,dt
=\frac{1}{c}\int_{0}^{c}\frac{1}{N}\sum_{n=0}^{N-1}f(nt)\,dt.
$$
Since $\frac{1}{N}\sum_{n=0}^{N-1}f(nt)\ra L_{t}$ as $N\ra\infty$ for \ae\ $t\in[0,c]$,
by \rfr{P-dct}, $\frac{1}{N}\sum_{n=0}^{N-1}v_{n}\ra\frac{1}{c}\int_{0}^{c}L_{t}\,dt$.
By \rfr{P-Taub}, $v_{n}\ra\frac{1}{c}\int_{0}^{c}L_{t}\,dt$.
On the other hand, $\lim_{n\ras\infty}v_{n}=\lim_{b\ras\infty}\frac{1}{b}\int_{0}^{b}f(x)\,dx$.
So, $L=\lim_{b\ras\infty}\frac{1}{b}\int_{0}^{b}f(x)\,dx$ exists 
and equals $\frac{1}{c}\int_{0}^{c}L_{t}\,dt$.

Next, for any $z\in(0,c)$ we also have $L=\frac{1}{z}\int_{0}^{z}L_{t}\,dt$.
So, for any $z\in[0,c]$, $\int_{0}^{z}L_{t}\,dt=zL$,
which implies that $L_{t}=L$ a.e. on $[0,c]$.
\endproof
\tsubsection{s-mulU}{Uniform Ces\`{a}ro limits}
\theorem{P-mulU}{}{}
Let a bounded measurable function $f\col[0,\infty)\ra V$
be such that for some $c>0$, for \ae\ $t\in(0,c]$
the limit $L_{t}=\lim_{b-a\ras\infty}\frac{1}{b-a}\sum_{n>a}^{b}f(nt)$ exists.
Then $L_{t}$ is constant \ae\ on $(0,c]$, 
$L_{t}=L\in V$ for \ae\ $t\in(0,c]$,
and $\lim_{b-a\ras\infty}\frac{1}{b-a}\int_{a}^{b}f(x)\,dx$ exists and equals $L$.
\endtheorem
\lemma{P-asU}{}{}
Let $(v_{n})$ be a bounded sequence in $V$ 
such that $\lim_{b-a\ras\infty}\frac{1}{b-a}\sum_{n>a}^{b}v_{n}=0$,
let $(b_{k})$ be a sequence of positive real numbers with $b_{k}\ras\infty$,
and let $(\alf_{k})$, $(\bet_{k})$ be sequences of real numbers
such that $0<\bet_{k}-\alf_{k}\leq b_{k}$ for all $k$.
Then $\lim_{k\ras\infty}\frac{1}{b_{k}}\sum_{n>\alf_{k}}^{\bet_{k}}v_{n}=0$.
\endlemma
\proof{}
Assume that $\|v_{n}\|\leq 1$ for all $n$.
Let $\eps>0$.
Let $B>1$ be such that
$\bigl\|\frac{1}{b-a}\sum_{n>a}^{b}v_{n}\bigr\|<\eps$
whenever $b-a>B$.
Let $K$ be such that $b_{k}>(B+1)/\eps$ for all $k>K$.
Then for any $k>K$,
if $\bet_{k}-\alf_{k}>B$,
then 
$\bigl\|\frac{1}{b_{k}}\sum_{n>\alf_{k}}^{\bet_{k}}v_{n}\bigr\|
\leq\bigl\|\frac{1}{\bet_{k}-\alf_{k}}\sum_{n>\alf_{k}}^{\bet_{k}}v_{n}\bigr\|<\eps$,
and if $\bet_{k}-\alf_{k}\leq B$,
then also $\bigl\|\frac{1}{b_{k}}\sum_{n>\alf_{k}}^{\bet_{k}}v_{n}\bigr\|
\leq\frac{\bet_{k}-\alf_{k}+1}{b_{k}}<\eps$.%
\endproof
\proof{of \rfr{P-mulU}}
Since, in particular, 
$L_{t}=\lim_{b\ras\infty}\frac{1}{b}\sum_{n>0}^{b}f(nt)$ for \ae\ $t\in(0,c]$,
we have $L_{t}=\const=L$ for \ae\ $t\in(0,c]$ by \rfr{P-mulC}.
Replacing $f$ by $f-L$, we may assume that $L=0$.
After replacing $f(x)$ by $f(cx/2)$, we assume that $c=2$.
Let $a\geq 0$, $b\geq a+1$.
For any $n\in\N$,
$\frac{1}{n}\int_{a}^{b}f(x)\,dx=\int_{a/n}^{b/n}f(nt)\,dt$.
Adding these equalities for all $n\in(b/2,b]$,
and taking into account that $b/n<2$ for $n>b/2$,
we get $\lam\int_{a}^{b}f(x)\,dx=\int_{0}^{2}\sum_{n>\alf(a,b,t)}^{\bet(a,b,t)}f(nt)\,dt$,
where $\lam=\sum_{n>b/2}^{b}\frac{1}{n}\geq\frac{1}{2}$,
and for every $t\in(0,2]$, $\alf(a,b,t)=\max\{b/2,a/t\}$ and $\bet(a,b,t)=\min\{b,b/t\}$.
Thus, 
\equn{f-dmnunf}{
\Bigl\|\frac{1}{b-a}\int_{a}^{b}f(x)\,dx\Bigr\|
\leq 2\Bigl\|\int_{0}^{2}f_{a,b}(t)\,dt\Bigr\|
}
where $f_{a,b}(t)=\frac{1}{b-a}\sum_{n>\alf(a,b,t)}^{\bet(a,b,t)}f(nt)$, $t\in(0,2]$.

We will now show that the functions $f_{a,b}$, for $a\geq 0$, $b\geq a+1$,
are uniformly bounded.
Let us assume that $\sup_{x\in(0,\infty)}\|f(x)\|\leq 1$.
If $a\leq b/2$, then $b-a\geq b/2$,
and since $\bet(a,b,t)-\alf(a,b,t)\leq b/2$,
for any $t\in(0,2]$ we have $\|f_{a,b}(t)\|\leq\frac{1}{b/2}(b/2+1)\leq 3$.
If $a>b/2$, then for any $t\in(0,1/2]$ 
we have $\alf(a,b,t)\geq a/t\geq 2a>b\geq\bet(a,b,t)$,
so $f_{a,b}(t)=0$;
and since, for any $t\in(0,2]$, $\bet(a,b,t)-\alf(a,b,t)\leq(b-a)/t$,
we have $\|f_{a,b}(t)\|\leq\frac{1}{b-a}((b-a)/t+1)\leq 3$ for $t\in[1/2,2]$.

For \ae\ $t\in(0,2]$, 
since $\bet(a,b,t)-\alf(a,b,t)\leq(b-a)/t$ for all $a\geq 0$, $b\geq a+1$, 
by \rfr{P-asU}, 
$$
\lim_{b-a\ras\infty}f_{a,b}(t)
=\frac{1}{t}\lim_{b-a\ras\infty}\frac{t}{b-a}\sum_{n>\alf(a,b,t)}^{\bet(a,b,t)}f(nt)=0.
$$
Hence, by \rfr{P-dct},
$\int_{0}^{2}f_{a,b}(t)\,dt\ra 0$ as $b-a\ras\infty$.
So, by \frfr{f-dmnunf}, 
$\frac{1}{b-a}\int_{a}^{b}f(x)\,dx\ra 0$ as $b-a\ras\infty$.%
\endproof
\tsubsection{s-mulSU}{Liminf and limsup versions for uniform averages}
In the case $f$ is a real-valued function and the limits $L_{t}$ do not exist
we have liminf/limsup versions of the above theorems.
We start with the uniform case:
\theorem{P-mulSU}{}{}
For any bounded measurable function $f\col[0,\infty)\ra\R$ and any $c>0$,
$$
\liminf_{b-a\ras\infty}\frac{1}{b-a}\int_{a}^{b}f(x)\,dx
\geq\frac{1}{c}\int_{0}^{c}\liminf_{b-a\ras\infty}\frac{1}{b-a}\sum_{n>a}^{b}f(nt)\,dt
$$
and
$$
\limsup_{b-a\ras\infty}\frac{1}{b-a}\int_{a}^{b}f(x)\,dx
\leq\frac{1}{c}\int_{0}^{c}\limsup_{b-a\ras\infty}\frac{1}{b-a}\sum_{n>a}^{b}f(nt)\,dt.
$$
\endtheorem
\proof{}
We will only prove the first inequality.
We may assume that $f\geq 0$.
Let $L=\liminf\limits_{b-a\ras\infty}\frac{1}{b-a}\int_{a}^{b}f(x)\,dx$;
find a sequence of intervals $[a_{k},b_{k}]$ with $b_{k}-a_{k}\ra\infty$
such that $L=\lim_{k\ras\infty}\frac{1}{b_{k}-a_{k}}\int_{a_{k}}^{b_{k}}f(x)\,dx$.
We may assume that $a_{k}\ra\infty$ 
(after replacing each $a_{k}$ by $\max\{a_{k},\sqrt{b_{k}}\}$),
and that $b_{k}/a_{k}\ra 1$
(after replacing each interval $[a_{k},b_{k}]$ by a suitable subinterval).

For any $t>0$ we have
$$
\liminf_{b-a\ras\infty}\frac{1}{b-a}\sum_{n>a}^{b}f(nt)
\leq
\liminf_{k\ras\infty}\frac{1}{b_{k}/t-a_{k}/t}\sum_{n>a_{k}/t}^{b_{k}/t}f(nt),
$$
so
$$
\frac{1}{c}\int_{0}^{c}\liminf_{b-a\ras\infty}\frac{1}{b-a}\sum_{n>a}^{b}f(nt)\,dt
\leq
\frac{1}{c}\int_{0}^{c}\liminf_{k\ras\infty}\frac{1}{b_{k}/t-a_{k}/t}
\sum_{n>a_{k}/t}^{b_{k}/t}f(nt)\,dt.
$$
By (the classical) Fatou's lemma we have
$$
\frac{1}{c}\int_{0}^{c}\liminf_{k\ras\infty}\frac{1}{b_{k}/t-a_{k}/t}
\sum_{n>a_{k}/t}^{b_{k}/t}f(nt)\,dt
\leq
\frac{1}{c}\liminf_{k\ras\infty}\int_{0}^{c}\frac{1}{b_{k}/t-a_{k}/t}
\sum_{n>a_{k}/t}^{b_{k}/t}f(nt)\,dt,
$$
and
\lequ{
\frac{1}{c}\liminf_{k\ras\infty}\int_{0}^{c}\frac{1}{b_{k}/t-a_{k}/t}
\sum_{n>a_{k}/t}^{b_{k}/t}f(nt)\,dt
=
\frac{1}{c}\liminf_{k\ras\infty}\frac{1}{b_{k}-a_{k}}\int_{0}^{c}
\sum_{n>a_{k}/t}^{b_{k}/t}tf(nt)\,dt
\-\\\-
\leq
\frac{1}{c}\liminf_{k\ras\infty}\frac{1}{b_{k}-a_{k}}\sum_{n>a_{k}/c}
\int_{a_{k}/n}^{b_{k}/n}tf(nt)\,dt.
}
For every $k\in\N$, for each $n$ we have
$$
I_{n}=\int_{a_{k}/n}^{b_{k}/n}tf(nt)\,dt
=\frac{1}{n^{2}}\int_{a_{k}}^{b_{k}}xf(x)\,dx
=\frac{1}{n^{2}}\alf_{k}\int_{a_{k}}^{b_{k}}f(x)\,dx
$$
with $\alf_{k}\in[a_{k},b_{k}]$, so
$$
\sum_{n>a_{k}/c}I_{n}
=\alf_{k}\int_{a_{k}}^{b_{k}}f(x)\,dx\sum_{n>a_{k}/c}\frac{1}{n^{2}}
=\alf_{k}s_{k}\int_{a_{k}}^{b_{k}}f(x)\,dx,
$$
where $s_{k}=\sum_{n>a_{k}/c}\frac{1}{n^{2}}$ satisfies $s_{k}a_{k}/c\ra 1$ as $k\ra\infty$.
Since, by our assumption, also $\alf_{k}/a_{k}\ra 1$, we get
$$
\frac{1}{c}\lim_{k\ras\infty}\frac{1}{b_{k}-a_{k}}\sum_{n>a_{k}/c}
\int_{a_{k}/n}^{b_{k}/n}tf(nt)\,dt
=\lim_{k\ras\infty}\frac{\alf_{k}s_{k}}{c}\cd\frac{1}{b_{k}-a_{k}}\int_{a_{k}}^{b_{k}}f(x)\,dx
=L.
$$
So, $\frac{1}{c}\int_{0}^{c}
\liminf\limits_{b-a\ras\infty}\frac{1}{b-a}\sum_{n>a}^{b}f(nt)\,dt\leq L$.
\endproof
\tsubsection{s-mulSC}{Liminf and limsup versions for standard averages}
\theorem{P-mulSC}{}{}
For any bounded measurable function $f\col[0,\infty)\ra\R$ and any $c>0$,
$$
\liminf_{b\ras\infty}\frac{1}{b}\int_{0}^{b}f(x)\,dx
\geq\frac{1}{c}\int_{0}^{c}\liminf_{b\ras\infty}\frac{1}{b}\sum_{n>0}^{b}f(nt)\,dt
$$
and
$$
\limsup_{b\ras\infty}\frac{1}{b}\int_{0}^{b}f(x)\,dx
\leq\frac{1}{c}\int_{0}^{c}\limsup_{b\ras\infty}\frac{1}{b}\sum_{n>0}^{b}f(nt)\,dt.
$$
\endtheorem
\proof{}
We may assume that $f\geq 0$.
Let $L=\liminf_{b\ras\infty}\frac{1}{b}\int_{0}^{b}f(x)\,dx$;
choose a sequence $(b_{k})$, with $b_{k}\ra\infty$ as $k\ra\infty$,
such that $\lim_{k\ras\infty}\frac{1}{b_{k}}\int_{0}^{b_{k}}f(x)\,dx=L$.
Then also $\lim_{k\ras\infty}\frac{1}{b_{k}-a_{k}}\int_{a_{k}}^{b_{k}}f(x)\,dx=L$,
where $a_{k}=\sqrt{b_{k}}$, $k\in\N$.
For all $t>0$ we have 
$$
\liminf_{b\ras\infty}\frac{1}{b}\sum_{n>0}^{b}f(nt)
=\liminf_{b\ras\infty}\frac{1}{b-\sqrt{b}}\sum_{n>\sqrt{b}}^{b}f(nt)
\leq\liminf_{k\ras\infty}\frac{1}{b_{k}-a_{k}}\sum_{n>a_{k}}^{b_{k}}f(nt),
$$
and as in the proof of \rfr{P-mulSU}, by Fatou's lemma,
$$
\frac{1}{c}\int_{0}^{c}\liminf_{k\ras\infty}
\frac{1}{b_{k}-a_{k}}\sum_{n>a_{k}}^{b_{k}}f(nt)\,dt
\leq\frac{1}{c}\liminf_{k\ras\infty}\frac{1}{b_{k}-a_{k}}
\int_{0}^{c}\sum_{n>a_{k}/t}^{b_{k}/t}tf(nt)\,dt,
$$
so it suffices to show that this last expression is $\leq L$.

For every $k\in\N$
put $M_{k}=\lfloor b_{k}^{2/3}\rfloor$
and subdivide the interval $\bigl[a_{k},b_{k}\bigr]$ into $M_{k}$ equal parts:
put $b_{k,j}=a_{k}+j(b_{k}-a_{k})/M_{k}$, $j=0,\ld,M_{k}$.
As in the proof of \rfr{P-mulSU},
for any $k$ and $j$ we have
\equn{f-msc}{
\frac{1}{b_{k,j}-b_{k,j-1}}\int_{0}^{c}\sum_{n>b_{k,j-1}/t}^{b_{k,j}/t}tf(nt)\,dt
\leq\frac{\alf_{k,j}s_{k,j}}{b_{k,j}-b_{k,j-1}}\int_{b_{k,j-1}}^{b_{k,j}}f(x)\,dx,
}
where $\alf_{k,j}\in[b_{k,j-1},b_{k,j}]$ and $s_{k,j}=\sum_{n>b_{k,j-1}/c}\frac{1}{n^{2}}$.
Since the function $\phi(t)=\frac{t+\del}{t-1}$ with $\del>0$ is decreasing,
for any $k$ and any $j$ we have
$$
\alf_{k,j}s_{k,j}<\frac{b_{k,j}}{(b_{k,j-1}-1)/c}
\leq c\frac{a_{k}+b_{k}/M_{k}}{a_{k}-1}
=c\frac{b_{k}^{1/2}+b_{k}/\lfloor b_{k}^{2/3}\rfloor}{b_{k}^{1/2}-1}=:r_{k},
$$
which tends to $c$ as $k\ra\infty$.
Replacing $\alf_{k,j}s_{k,j}$ by $r_{k}$
and taking the average of both sides of the inequality \frfr{f-msc} 
for a fixed $k$ and $j=1,\ld,M_{k}$ we get
$$
\frac{1}{b_{k}-a_{k}}\int_{0}^{c}\sum_{n>a_{k}/t}^{b_{k}/t}tf(nt)\,dt
<\frac{r_{k}}{b_{k}-a_{k}}\int_{a_{k}}^{b_{k}}f(x)\,dx,
$$
so
$$
\frac{1}{c}\liminf_{k\ras\infty}\frac{1}{b_{k}-a_{k}}
\int_{0}^{c}\sum_{n>a_{k}/t}^{b_{k}/t}tf(nt)\,dt
\leq\lim_{k\ras\infty}\frac{r_{k}}{c(b_{k}-a_{k})}\int_{a_{k}}^{b_{k}}f(x)\,dx=L.
$$
\frgdsp\endproof
\section{S-mulM}{Multiplicative method -- the multiparameter case}
\subsectionpenalty=9000
\tsubsection{s-mulMC}{Standard Ces\`{a}ro limits}
We will use the notation introduced in \rfr{s-addMC}.
\theorem{P-mulMC}{}{}
Let a bounded measurable function $f\col\Rp^{d}\ra V$ 
be such that for some $c\in\Rp^{d}$, $c>0$,
the limit $L_{t}=\lim_{l(b)\ras\infty}\frac{1}{w(b)}\sum_{n>0}^{b}f(nt)$ exists 
for \ae\ $t\in[0,c]$.
Then $L_{t}$ is \ae\ constant, $L_{t}=L\in V$ \ae\ on $[0,c]$,
and $\lim_{l(b)\ras\infty}\frac{1}{w(b)}\int_{0}^{b}f(x)\,dx$ exists and equals $L$.
\endtheorem
Let $e_{1}=(1,0,\ld,0,0)$, $\ld$, $e_{d}=(0,0,\ld,0,1)$.
\lemma{P-mulTaub}{}{}
Let $(v_{n})$ be an $\N^{d}$-sequence in $V$
such that the limit
$v=\lim_{l(N)\ras\infty}\frac{1}{w(N)}\sum_{n\leq N}v_{n}$ exists
and for some $\alf>0$, 
for any $n=(n_{1},\ld,n_{d})\in\N^{d}$ and for any $i$ one has $\|v_{n+e_{i}}-v_{n}\|<\alf n_{i}$.
Then $\lim_{l(n)\ras\infty}v_{n}=v$.
\endlemma
\proof{}
We may assume that $v=0$, that is, $\lim_{l(N)\ras\infty}\frac{1}{w(N)}\sum_{n\leq N}v_{n}=0$.
Assume that $v_{n}\not\ras 0$ as $l(n)\ras\infty$;
let $\eps>0$ be such that for any $N\in\N^{d}$ there exists $n>N$ such that $\|v_{n}\|>\eps$.
Put $\del=\frac{\eps(\eps/4\alf d)^{d}}{2(2^{d}+1)(1+\eps/2\alf d)^{d}}$.
Find $N\in\N^{d}$ such that $\bigl\|\frac{1}{w(M)}\sum_{n\leq M}v_{n}\bigr\|<\del$ for all $M>N$.
Find $M=(M_{1},\ld,M_{d})>N$, such that $\|v_{M}\|>\eps$ and $1/M_{i}<\eps/4\alf d$ for all $i$.
Let, by the Hahn-Banach theorem, $\phi\in V^{*}$ 
be such that $|\phi(v)|\leq\|v\|$ for all $v\in V$ and $\phi(v_{M})=\|v_{M}\|$.
Then for any $n=(n_{1},\ld,n_{d})>M$,
\equ{
\phi(v_{n})
=\phi(v_{M})+\sum_{i=1}^{d}\sum_{m=M_{i}}^{n_{i}-1}
\phi\bigl(v_{(n_{1},\ld,n_{i-1},m+1,M_{i+1},\ld,M_{d})}
-v_{(n_{1},\ld,n_{i-1},m,M_{i+1},\ld,M_{d})}\bigr)
\-\\
\geq
\bigl\|v_{M}\|-\sum_{m=M}^{n-1}
\bigl\|v_{(n_{1},\ld,n_{i-1},m+1,M_{i+1},\ld,M_{d})}
-v_{(n_{1},\ld,n_{i-1},m,M_{i+1},\ld,M_{d})}\bigr\|
\\\-
>\eps-\sum_{i=1}^{d}(n_{i}-M_{i})\frac{\alf}{M_{i}},
}
which is $\geq\eps/2$ when $\frac{(n_{i}-M_{i})\alf}{M_{i}}\leq\frac{\eps}{2d}$ for all $i$,
that is, when $M_{i}<n_{i}\leq M_{i}+\eps M_{i}/2\alf d$ for all $i$.
Put $K_{i}=M_{i}+\lfloor\eps M_{i}/2\alf d\rfloor$ and $K=(K_{1},\ld,K_{d})$,
then $\phi(v_{n})>\frac{\eps}{2}$ for $M\leq n\leq K$
and $w(K)\leq w(M)(1+\eps/2\alf d)^{d}$, 
$w(K-M)\geq w(M)\prod_{i=1}^{d}(\eps/2\alf d-1/M_{i})\geq w(M)(\eps/4\alf d)^{d}$.
Now, we can represent $\sum_{n\leq K}v_{n}$ as an alternating sum
$$
\sum_{n\leq K}v_{n}=\sum_{j=1}^{2^{d}-1}\Bigl(\pm\sum_{n\leq R_{j}}v_{n}\Bigr)+\sum_{n>M}^{K}v_{n},
$$
where for each $j$, for every $i$,
the $i$th entry of $R_{j}$ is either equal to $M_{i}$ or to $K_{i}$.
(For $d=2$, for instance, the formula is
$\sum_{n\leq K}v_{n}=\sum_{n\leq(M_{1},K_{2})}v_{n}+\sum_{n\leq(K_{1},M_{2})}v_{n}
-\sum_{n\leq(M_{1},M_{2})}v_{n}+\sum_{n>M}^{K}v_{n}$.)
For each $j$, $\bigl\|\sum_{n\leq R_{j}}v_{n}\bigr\|<w(R_{j})\del\leq w(K)\del$,
thus
\equ{
\Bigl\|\frac{1}{w(K)}\sum_{n\leq K}v_{n}\Bigr\|
\geq\frac{1}{w(K)}\Bigl(\Bigl\|\sum_{n>M}^{K}v_{n}\Bigr\|-2^{d}w(K)\del\Bigr)
\geq\frac{1}{w(K)}\,\phi\Bigl(\sum_{n=M+1}^{K}v_{n}\Bigr)-2^{d}\del
=\frac{1}{w(K)}\sum_{n=M+1}^{K}\phi(v_{n})-2^{d}\del
\-\\
>\frac{w(K-M)}{w(K)}\cd\frac{\eps}{2}-2^{d}\del
\geq\frac{w(M)(\eps/4\alf d)^{d}}{w(M)(1+\eps/2\alf d)^{d}}\cd\frac{\eps}{2}-2^{d}\del
=\frac{\eps(\eps/4\alf d)^{d}}{2(1+\eps/2\alf d)^{d}}-2^{d}\del=\del,
}
which contradicts the choice of $N$.
\endproof
\proof{of \rfr{P-mulMC}}
Let $v_{n}=\frac{1}{w(nc)}\int_{0}^{nc}f(x)\,dx$, $n\in\N^{d}$.
Then for any $n=(n_{1},\ld,n_{d})\in\N^{d}$ and any $i\in\{1,\ld,d\}$,
\lequ{
\|v_{n+e_{i}}-v_{n}\|
=\frac{1}{w(c)}\Bigl\|\Bigl(\frac{1}{w(n+e_{i})}-\frac{1}{w(n)}\Bigr)\int_{0}^{nc}f(x)\,dx
+\frac{1}{w(n+e_{i})}\int_{nc}^{(n+e_{i})c}f(x)\,dx\Bigr\|
\-\\\-
\leq\frac{w(cn)\sup\|f\|}{w(c)w(n)w(n+e_{i})}+\frac{w(c)\sup\|f\|}{w(c)w(n+e_{i}))}
=2\sup\|f\|/(n_{i}+1).
}
Also we have, for any $N\in\N^{d}$,
\equ{
\frac{1}{w(N)}\sum_{n\leq N}v_{n}
=\frac{1}{w(N)}\sum_{n\leq N}\frac{1}{w(nc)}\int_{0}^{nc}f(x)\,dx
=\frac{1}{w(N)}\sum_{n\leq N}\frac{1}{w(c)}\int_{0}^{c}f(nt)\,dt
=\frac{1}{w(c)}\int_{0}^{c}\frac{1}{w(N)}\sum_{n\leq N}f(nt)\,dt.
}
Since $\frac{1}{w(N)}\sum_{n\leq N}f(nt)\ra L_{t}$ as $l(N)\ra\infty$ for \ae\ $t\in[0,c]$,
by \rfr{P-dct}, $\frac{1}{w(N)}\sum_{n\leq N}v_{n}\ra\frac{1}{w(c)}\int_{0}^{c}L_{t}\,dt$.
By \rfr{P-mulTaub}, $v_{n}\ra\frac{1}{w(c)}\int_{0}^{c}L_{t}\,dt$.
On the other hand, $\lim_{l(n)\ras\infty}v_{n}
=\lim_{l(b)\ras\infty}\frac{1}{w(b)}\int_{0}^{b}f(x)\,dx$.
So, $L=\lim_{l(b)\ras\infty}\frac{1}{w(b)}\int_{0}^{b}f(x)\,dx$ exists
and equals $\frac{1}{w(c)}\int_{0}^{c}L_{t}\,dt$.

Next, for any $z\in(0,c)$, we also have $L=\frac{1}{w(z)}\int_{0}^{z}L_{t}\,dt$.
So, for any $z\in(0,c]$, $\int_{0}^{z}L_{t}\,dt=w(z)L$,
which implies that $L_{t}=L$ a.e. on $[0,c]$.
\endproof
\tsubsection{s-mulMU}{Uniform Ces\`{a}ro limits}
\theorem{P-mulMU}{}{}
Let a bounded measurable function $f\col\Rp^{d}\ra V$
be such that for some $c\in\Rp^{d}$, $c>0$, for \ae\ $t\in(0,c]$
the limit $L_{t}=\lim_{b-a\ras\infty}\frac{1}{b-a}\sum_{n>a}^{b}f(nt)$ exists.
Then $L_{t}$ is constant \ae\ on $(0,c]$, 
$L_{t}=L\in V$ for \ae\ $t\in(0,c]$,
and $\lim_{b-a\ras\infty}\frac{1}{b-a}\int_{a}^{b}f(x)\,dx$ exists and equals $L$.
\endtheorem
\lemma{P-asMU}{}{}
Let $(v_{n})$ be a bounded $\N^{d}$-sequence in $V$
such that $\lim_{l(b-a)\ras\infty}\frac{1}{w(b-a)}\sum_{n>a}^{b}v_{n}=0$,
and let $(b_{k})$, $(\alf_{k})$ and $(\bet_{k})$ be sequences in $\Rp^{d}$
such that $0<\bet_{k}-\alf_{k}\leq b_{k}$ for all $k$ and $l(b_{k})\ras\infty$.
Then $\lim_{k\ras\infty}\frac{1}{w(b_{k})}\sum_{n>\alf_{k}}^{\bet_{k}}v_{n}=0$.
\endlemma
\proof{}
Assume that $\|v_{n}\|\leq 1$ for all $n$.
Let $\eps>0$.
Let $B>1$ be such that
$\bigl\|\frac{1}{w(b-a)}\sum_{n>a}^{b}v_{n}\bigr\|<\eps$
whenever $l(b-a)>B$.
Let $K$ be such that $l(b_{k})>(B+1)/\eps$ for all $k>K$.
Then for any $k>K$,
if $l(\bet_{k}-\alf_{k})>B$,
then $\bigl\|\frac{1}{w(b_{k})}\sum_{n>\alf_{k}}^{\bet_{k}}v_{n}\bigr\|
\leq\bigl\|\frac{1}{w(\bet_{k}-\alf_{k})}\sum_{n>\alf_{k}}^{\bet_{k}}v_{n}\bigr\|<\eps$,
and if $l(\bet_{k}-\alf_{k})\leq B$,
then also $\bigl\|\frac{1}{w(b_{k})}\sum_{n>\alf_{k}}^{\bet_{k}}v_{n}\bigr\|
\leq\frac{l(\bet_{k}-\alf_{k})+1}{l(b_{k})}<\eps$.%
\endproof
\proof{of \rfr{P-mulMU}}
Since, in particular, 
$L_{t}=\lim_{l(b)\ras\infty}\frac{1}{w(b)}\sum_{n>0}^{b}f(nt)$ for \ae\ $t\in(0,c]$,
we have $L_{t}=\const=L$ for \ae\ $t\in(0,c]$ by \rfr{P-mulMC}.
Replacing $f$ by $f-L$, we may assume that $L=0$.
After replacing $f(x)$ by $f(cx/2)$, we assume that $c=\bar{2}$.
Let $a\geq 0$, $b\geq a+1$.
For any $n\in\N$,
$\frac{1}{w(n)}\int_{a}^{b}f(x)\,dx=\int_{a/n}^{b/n}f(nt)\,dt$.
Adding these equalities for all $n\in(b/2,b]$,
and taking into account that $b/n<\bar{2}$ for $n>b/2$,
we get $\lam\int_{a}^{b}f(x)\,dx=\int_{0}^{\bar{2}}\sum_{n>\alf(a,b,t)}^{\bet(a,b,t)}f(nt)\,dt$,
where $\lam=\sum_{n>b/2}^{b}\frac{1}{w(n)}\geq\frac{1}{2^{d}}$,
and for every $t\in(0,2]^{d}$, $\alf(a,b,t)=\max\{b/2,a/t\}$ and $\bet(a,b,t)=\min\{b,b/t\}$.
Thus, 
\equn{f-muldmnunf}{
\Bigl\|\frac{1}{w(b-a)}\int_{a}^{b}f(x)\,dx\Bigr\|
\leq 2\Bigl\|\int_{0}^{\bar{2}}f_{a,b}(t)\,dt\Bigr\|
}
where $f_{a,b}(t)=\frac{1}{w(b-a)}\sum_{n>\alf(a,b,t)}^{\bet(a,b,t)}f(nt)$, $t\in(0,2]^{d}$.

We will now show that the functions $f_{a,b}$, for $a\geq 0$, $b\geq a+\bar{1}$,
are uniformly bounded.
Let us assume that $\sup_{x\in(0,\infty)}\|f(x)\|\leq 1$.
If $a\leq b/2$, then $b-a\geq b/2$,
and since $\bet(a,b,t)-\alf(a,b,t)\leq b/2$,
for any $t\in(0,2]^{d}$ we have $\|f_{a,b}(t)\|\leq\frac{1}{w(b/2)}w(b/2+\bar{1})\leq 3^{d}$.
If $a>b/2$, then for any $t\in(0,2]^{d}$ with $t_{i}<1/2$ for some $i$ 
we have $\alf(a,b,t)_{i}\geq a_{i}/t_{i}\geq 2a_{i}>b_{i}\geq\bet(a,b,t)_{i}$,
so $f_{a,b}(t)=0$;
and since, for any $t\in(0,2]^{d}$, $\bet(a,b,t)-\alf(a,b,t)\leq(b-a)/t$,
we have $\|f_{a,b}(t)\|\leq\frac{1}{w(b-a)}w((b-a)/t+\bar{1})\leq 3^{d}$ for all $t\in[1/2,2]^{d}$.

For \ae\ $t\in(0,2]^{d}$, 
since $\bet(a,b,t)-\alf(a,b,t)\leq(b-a)/t$ for all $a,b\in\Rp^{d}$, $b\geq a+\bar{1}$, 
by \rfr{P-asMU}, 
$$
\lim_{l(b-a)\ras\infty}f_{a,b}(t)
=\frac{1}{w(t)}\lim_{l(b-a)\ras\infty}\frac{1}{w((b-a)/t)}\sum_{n>\alf(a,b,t)}^{\bet(a,b,t)}f(nt)=0.
$$
Hence, by \rfr{P-dct},
$\int_{0}^{\bar{2}}f_{a,b}(t)\,dt\ra 0$ as $l(b-a)\ras\infty$.
So, by \frfr{f-muldmnunf},
$\frac{1}{w(b-a)}\int_{a}^{b}f(x)\,dx\ra 0$ as $l(b-a)\ras\infty$.%
\endproof
\tsubsection{s-mulMSU}{Liminf and limsup versions for uniform limits}
\theorem{P-mulMSU}{}{}
For any bounded measurable function $f\col\Rp^{d}\ra\R$ and any $c\in\Rp^{d}$, $c>0$,
one has
$$
\liminf_{l(b-a)\ras\infty}\frac{1}{w(b-a)}\int_{a}^{b}f(x)\,dx
\geq\frac{1}{w(c)}\int_{0}^{c}\liminf_{l(b-a)\ras\infty}\frac{1}{w(b-a)}
\sum_{n>a}^{b}f(nt)\,dt
$$
and
$$
\limsup_{l(b-a)\ras\infty}\frac{1}{w(b-a)}\int_{a}^{b}f(x)\,dx
\leq\frac{1}{w(c)}\int_{0}^{c}\limsup_{l(b-a)\ras\infty}\frac{1}{w(b-a)}
\sum_{n>a}^{b}f(nt)\,dt.
$$
\endtheorem
\proof{}
We will only prove the first inequality.
We may assume that $f\geq 0$.
Let $L=\liminf\limits_{l(b-a)\ras\infty}\frac{1}{w(b-a)}\int_{a}^{b}f(x)\,dx$;
find a sequence of intervals $[a_{k},b_{k}]\sln\Rp^{d}$ with $l(b_{k}-a_{k})\ra\infty$
such that $L=\lim_{k\ras\infty}\frac{1}{w(b_{k}-a_{k})}\int_{a_{k}}^{b_{k}}f(x)\,dx$.
We may assume that $l(a_{k})\ra\infty$
(after replacing each $a_{k}$ by $\max\{a_{k},\sqrt{b_{k}}\}$),
and that $w(b_{k})/w(a_{k})\ra 1$
(after replacing each interval $[a_{k},b_{k}]$ by a suitable subinterval).

For any $t>0$ we have
$$
\liminf_{l(b-a)\ras\infty}\frac{1}{w(b-a)}\sum_{n>a}^{b}f(nt)
\leq
\liminf_{k\ras\infty}\frac{1}{w(b_{k}/t-a_{k}/t)}\sum_{n>a_{k}/t}^{b_{k}/t}f(nt),
$$
so
$$
\frac{1}{w(c)}\int_{0}^{c}\liminf_{l(b-a)\ras\infty}\frac{1}{w(b-a)}\sum_{n>a}^{b}f(nt)\,dt
\leq
\frac{1}{w(c)}\int_{0}^{c}\liminf_{k\ras\infty}\frac{1}{w(b_{k}/t-a_{k}/t)}
\sum_{n>a_{k}/t}^{b_{k}/t}f(nt)\,dt.
$$
By (the classical) Fatou's lemma we have
$$
\frac{1}{w(c)}\int_{0}^{c}\liminf_{k\ras\infty}\frac{1}{w(b_{k}/t-a_{k}/t)}
\sum_{n>a_{k}/t}^{b_{k}/t}f(nt)\,dt
\leq
\frac{1}{w(c)}\liminf_{k\ras\infty}\int_{0}^{c}\frac{1}{w(b_{k}/t-a_{k}/t)}
\sum_{n>a_{k}/t}^{b_{k}/t}f(nt)\,dt,
$$
and
\lequ{
\frac{1}{w(c)}\liminf_{k\ras\infty}\int_{0}^{c}\frac{1}{w(b_{k}/t-a_{k}/t)}
\sum_{n>a_{k}/t}^{b_{k}/t}f(nt)\,dt
=
\frac{1}{w(c)}\liminf_{k\ras\infty}\frac{1}{w(b_{k}-a_{k})}\int_{0}^{c}
\sum_{n>a_{k}/t}^{b_{k}/t}tf(nt)\,dt
\-\\\-
\leq
\frac{1}{w(c)}\liminf_{k\ras\infty}\frac{1}{w(b_{k}-a_{k})}\sum_{n>a_{k}/c}
\int_{a_{k}/n}^{b_{k}/n}tf(nt)\,dt.
}
For every $k\in\N$, for each $n$ we have
$$
I_{n}=\int_{a_{k}/n}^{b_{k}/n}tf(nt)\,dt
=\frac{1}{w(n^{2})}\int_{a_{k}}^{b_{k}}xf(x)\,dx
=\frac{1}{w(n^{2})}w(\alf_{k})\int_{a_{k}}^{b_{k}}f(x)\,dx
$$
with $\alf_{k}\in[a_{k},b_{k}]$, so
$$
\sum_{n>a_{k}/c}I_{n}
=w(\alf_{k})\int_{a_{k}}^{b_{k}}f(x)\,dx\sum_{n>a_{k}/c}\frac{1}{w(n^{2})}
=w(\alf_{k})s_{k}\int_{a_{k}}^{b_{k}}f(x)\,dx,
$$
where $s_{k}=\sum_{n>a_{k}/c}\frac{1}{w(n^{2})}$ 
satisfies $s_{k}w(a_{k}/c)\ra 1$ as $k\ra\infty$.
Since, by our assumption, also $w(\alf_{k})/w(a_{k})\ra 1$, we get
$$
\frac{1}{w(c)}\lim_{k\ras\infty}\frac{1}{w(b_{k}-a_{k})}\sum_{n>a_{k}/c}
\int_{a_{k}/n}^{b_{k}/n}tf(nt)\,dt
=\lim_{k\ras\infty}\frac{w(\alf_{k})s_{k}}{w(c)}\cd\frac{1}{w(b_{k}-a_{k})}
\int_{a_{k}}^{b_{k}}f(x)\,dx=L.
$$
So, $\frac{1}{w(c)}\int_{0}^{c}
\liminf\limits_{l(b-a)\ras\infty}\frac{1}{w(b-a)}\sum_{n>a}^{b}f(nt)\,dt\leq L$.
\endproof
\tsubsection{s-mulMSC}{Liminf and limsup versions for standard averages}
\theorem{P-mulMSC}{}{}
For any bounded measurable function $f\col\Rp^{d}\ra\R$ and any $c\in\Rp^{d}$, $c>0$,
one has
$$
\liminf_{l(b)\ras\infty}\frac{1}{w(b)}\int_{0}^{b}f(x)\,dx
\geq\frac{1}{w(c)}\int_{0}^{c}\liminf_{l(b)\ras\infty}\frac{1}{w(b)}\sum_{n>0}^{b}f(nt)\,dt
$$
and
$$
\limsup_{l(b)\ras\infty}\frac{1}{w(b)}\int_{0}^{b}f(x)\,dx
\leq\frac{1}{w(c)}\int_{0}^{c}\limsup_{l(b)\ras\infty}\frac{1}{w(b)}\sum_{n>0}^{b}f(nt)\,dt.
$$
\endtheorem
\proof{}
We may assume that $f\geq 0$.
Let $L=\liminf_{b\ras\infty}\frac{1}{b}\int_{0}^{b}f(x)\,dx$;
choose a sequence $(b_{k})$ in $\Rp^{d}$, with $l(b_{k})\ra\infty$ as $k\ra\infty$,
such that $\lim_{k\ras\infty}\frac{1}{w(b_{k})}\int_{0}^{b_{k}}f(x)\,dx=L$.
Then also $\lim_{k\ras\infty}\frac{1}{w(b_{k}-a_{k})}\int_{a_{k}}^{b_{k}}f(x)\,dx=L$,
where $a_{k}=\sqrt{b_{k}}$, $k\in\N$.
For all $t>0$ we have
$$
\liminf_{l(b)\ras\infty}\frac{1}{w(b)}\sum_{n>0}^{b}f(nt)
=\liminf_{l(b)\ras\infty}\frac{1}{w(b-\sqrt{b})}\sum_{n>\sqrt{b}}^{b}f(nt)
\leq\liminf_{k\ras\infty}\frac{1}{w(b_{k}-a_{k})}\sum_{n>a_{k}}^{b_{k}}f(nt),
$$
and as in the proof of \rfr{P-mulMSU}, by Fatou's lemma,
$$
\frac{1}{w(c)}\int_{0}^{c}\liminf_{k\ras\infty}
\frac{1}{w(b_{k}-a_{k})}\sum_{n>a_{k}}^{b_{k}}f(nt)\,dt
\leq\frac{1}{w(c)}\liminf_{k\ras\infty}\frac{1}{w(b_{k}-a_{k})}
\int_{0}^{c}\sum_{n>a_{k}/t}^{b_{k}/t}tf(nt)\,dt,
$$
so it suffices to show that this last expression is $\leq L$.

For every $k\in\N$
put $M_{k}=\lfloor b_{k}^{2/3}\rfloor$
and subdivide the interval $\bigl[a_{k},b_{k}\bigr]$ into $w(M_{k})$ equal parts:
put $b_{k,j}=a_{k}+j(b_{k}-a_{k})/M_{k}$, $j\in(\{0\}\cup\N)^{d}\cap[0,M_{k}]$.
As in the proof of \rfr{P-mulMSU},
for any $k$ and $j$ we have
\equn{f-mmsc}{
\frac{1}{w(b_{k,j}-b_{k,j-1})}\int_{0}^{c}\sum_{n>b_{k,j-1}/t}^{b_{k,j}/t}tf(nt)\,dt
\leq\frac{w(\alf_{k,j})s_{k,j}}{w(b_{k,j}-b_{k,j-1})}\int_{b_{k,j-1}}^{b_{k,j}}f(x)\,dx,
}
where $\alf_{k,j}\in[b_{k,j-1},b_{k,j}]$ and $s_{k,j}=\sum_{n>b_{k,j-1}/c}\frac{1}{w(n^{2})}$.
For any $k$ and any $j$,
$$
w(\alf_{k,j})s_{k,j}<\frac{w(b_{k,j})}{w(b_{k,j-1}-\bar{1})/w(c)}
\leq w(c)\frac{w(a_{k}+b_{k}/M_{k})}{w(a_{k}-\bar{1})}
=w(c)\frac{w(b_{k}^{1/2}+b_{k}/\lfloor b_{k}^{2/3}\rfloor)}{w(b_{k}^{1/2}-\bar{1})}=:r_{k},
$$
which tends to $w(c)$ as $k\ra\infty$.
Replacing $w(\alf_{k,j})s_{k,j}$ by $r_{k}$
and taking the average of both sides of the inequality \frfr{f-mmsc} 
for a fixed $k$ and $j=1,\ld,M_{k}$ we get
$$
\frac{1}{w(b_{k}-a_{k})}\int_{0}^{c}\sum_{n>a_{k}/t}^{b_{k}/t}tf(nt)\,dt
<\frac{r_{k}}{w(b_{k}-a_{k})}\int_{a_{k}}^{b_{k}}f(x)\,dx,
$$
so
$$
\frac{1}{w(c)}\liminf_{k\ras\infty}\frac{1}{w(b_{k}-a_{k})}
\int_{0}^{c}\sum_{n>a_{k}/t}^{b_{k}/t}tf(nt)\,dt
\leq\lim_{k\ras\infty}\frac{r_{k}}{w(c)w(b_{k}-a_{k})}\int_{a_{k}}^{b_{k}}f(x)\,dx=L.
$$
\frgdsp\endproof
\tsubsection{s-mulMCsup}{A lim-limsup version for standard averages}
It turns out that if the limits 
$\lim_{l(b)\ras\infty}\frac{1}{w(b)}\sum_{n>0}^{b}f(nt)$, $t>0$,
do not exist, 
but, for some $L\in V$, one has
$\lim_{t\ras0^{+}}\limsup_{l(b)\ras\infty}\Bigl\|\frac{1}{w(b)}\sum_{n>0}^{b}f(nt)-L\Bigr\|=0$,
we still have the result.
For a function $h\col(0,r)^{d}\ra V$, $r>0$,
we write $\esslim_{t\ras 0^{+}}h(t)=h_{0}$ 
if for any $\eps>0$ there exists $\del\in\Rp^{d}$ 
such that $\|h(t)-h_{0}\|<\eps$ for \ae\ $t\in(0,\del]$.
\theorem{P-mulMCsup}{}{}
Let $f\col\Rp^{d}\ra V$ be a bounded measurable function satisfying, for some $L\in V$,
$$
\esslim_{t\ras 0^{+}}\limsup_{l(b)\ras\infty}
\Bigl\|\frac{1}{w(b)}\sum_{n>0}^{b}f(nt)-L\Bigr\|=0.
$$
Then $\lim_{l(b)\ras\infty}\frac{1}{w(b)}\int_{0}^{b}f(x)\,dx=L$.
\endtheorem
\lemma{P-asMCsup}{}{}
Let $(v_{n})$ be a bounded $\N^{d}$-sequence in $V$
and let $(b_{k})$, $(\bet_{k})$ be sequences in $\Rp^{d}$
with $0<\bet_{k}\leq b_{k}$ for all $k$ and $l(b_{k})\ras\infty$.
Then $\limsup_{k\ras\infty}\bigl\|\frac{1}{w(b_{k})}\sum_{n>0}^{\bet_{k}}v_{n}\bigr\|
\leq\limsup_{l(b)\ras\infty}\bigl\|\frac{1}{w(b)}\sum_{n>0}^{b}v_{n}\bigr\|$,
and $\limsup_{k\ras\infty}
\bigl\|\frac{1}{w(b_{k})}\sum_{n>b_{k}/2}^{\bet_{k}}v_{n}\bigr\|
\leq2^{d}\limsup_{l(b)\ras\infty}\bigl\|\frac{1}{w(b)}\sum_{n>0}^{b}v_{n}\bigr\|$.
\endlemma
\proof{}
Assume that $\|v_{n}\|\leq 1$ for all $n$.
Let $s=\limsup_{l(b)\ras\infty}\bigl\|\frac{1}{w(b)}\sum_{n>0}^{b}v_{n}\bigr\|$.
Let $\eps>0$,
and let $B>1$ be such that 
$\bigl\|\frac{1}{w(b)}\sum_{n>0}^{b}v_{n}\bigr\|<s+\eps$
whenever $l(b)>B$.
Let $K$ be such that $l(b_{k})>B/\eps$ for all $k>K$.
Then for any $k>K$,
if $l(\bet_{k})>B$,
then $\bigl\|\frac{1}{w(b_{k})}\sum_{n>0}^{\bet_{k}}v_{n}\bigr\|
\leq\bigl\|\frac{1}{w(\bet_{k})}\sum_{n>0}^{\bet_{k}}v_{n}\bigr\|<s+\eps$,
and if $l(\bet_{k})\leq B$,
then also $\bigl\|\frac{1}{w(b_{k})}\sum_{n>0}^{\bet_{k}}v_{n}\bigr\|
\leq\frac{w(\bet_{k})}{w(b_{k})}\leq\frac{l(\bet_{k})}{l(b_{k})}<\eps\leq s+\eps$.
So, $\limsup_{k\ras\infty}
\bigl\|\frac{1}{w(b_{k})}\sum_{n>0}^{\bet_{k}}v_{n}\bigr\|\leq s$.

For any $S\sle\{1,\ld,d\}$ and $a=(a_{1},\ld,a_{d}),\ b=(b_{1},\ld,b_{d})\in\R^{d}$,
let $\sig_{S}(a,b)=(c_{1},\ld,c_{d})$ where for each $i$,
$c_{i}=a_{i}$ if $i\in S$ and $c_{i}=b_{i}$ if $i\not\in S$.
Then for any $a,b\in\Rp^{d}$ with $a\leq b$
we have $\sum_{n>a}^{b}v_{n}
=\sum_{S\sle\{1,\ld,d\}}(-1)^{|S|}\sum_{n>0}^{\sig_{S}(a,b)}v_{n}$.
Since, for any $S\sle\{1,\ld,d\}$,
$\limsup_{k\ras\infty}
\bigl\|\frac{1}{w(b_{k})}\sum_{n>0}^{\sig_{S}(b_{k}/2,\bet_{k})}v_{n}\bigr\|\leq s$,
we also get that
$\limsup_{k\ras\infty}\bigl\|\frac{1}{w(b_{k})}\sum_{n>b_{k}/2}^{\bet_{k}}v_{n}\bigr\|
\leq 2^{d}s$.%
\endproof
\proof{of \rfr{P-mulMCsup}}
We may assume that $L=0$.
Fix $\eps>0$.
Find $\del\in\Rp^{d}$, $\del>0$, such that 
$\limsup_{b\ras\infty}\bigl\|\frac{1}{b}\sum_{n>0}^{b}f(nt)\bigr\|<\eps$
for \ae\ $t\in[0,\del]$,
and define $g(x)=f\bigl(\frac{\del}{2}x)$, $x\in\Rp^{d}$.
Then $\limsup_{b\ras\infty}\bigl\|\frac{1}{b}\sum_{n>0}^{b}g(nt)\bigr\|<\eps$
for \ae\ $t\in[0,2]^{d}$.

Let $b\in\Rp^{d}$, $b\geq\bar{1}$.
For any $n\in\N^{d}$,
$\frac{1}{w(n)}\int_{0}^{b}g(x)\,dx=\int_{0}^{b/n}g(nt)\,dt$.
Adding these equalities for all $n\in(b/2,b]$,
and taking into account that $b/n<\bar{2}$ for $n>b/2$,
we get $\lam\int_{0}^{b}g(x)\,dx=\int_{0}^{2}\sum_{n>b/2}^{\bet(b,t)}g(nt)\,dt$,
where $\lam=\sum_{n>b/2}^{b}\frac{1}{w(n)}\geq\frac{1}{2^{d}}$,
and for every $t\in(0,2]^{d}$, $\bet(b,t)=\min\{b,b/t\}$.
Thus, 
\equn{f-dmnmordsup}{
\Bigl\|\frac{1}{w(b)}\int_{0}^{b}g(x)\,dx\Bigr\|
\leq 2^{d}\Bigl\|\int_{0}^{\bar{2}}g_{b}(t)\,dt\Bigr\|,
}
where $g_{b}(t)=\frac{1}{w(b)}\sum_{n>b/2}^{\bet(b,t)}g(nt)$, $t\in(0,2]^{d}$.

Let us assume that $\sup_{x\in\Rp^{d}}\|g(x)\|=\sup_{x\in\Rp^{d}}\|f(x)\|\leq 1$.
Then for any $b\geq\bar{1}$,
for every $t\in(0,2]^{d}$ we have 
$\|g_{b}(t)\|\leq\hbox{$\frac{1}{w(b)}w(\bet(b,t))$}\leq 1$,
so the functions $g_{b}$ are uniformly bounded.
For \ae\ $t\in(0,2]^{d}$, since $\bet(b,t)\leq b$,
by \rfr{P-asMCsup} we have 
$\limsup_{l(b)\ras\infty}\|g_{b}(t)\|\leq2^{d}\eps$.
Hence, by \rfr{P-dctsup},
$\limsup_{b\ras\infty}\bigl\|\int_{0}^{\bar{2}}g_{b}(t)\,dt\bigr\|\leq 2^{2d}\eps$.
So, by \frfr{f-dmnmordsup},
$\limsup_{b\ras\infty}\bigl\|\frac{1}{b}\int_{0}^{b}g(x)\,dx\bigr\|\leq 2^{3d}\eps$.
Since for any $b\in\Rp$, $b>0$, we have
$\frac{1}{w(b)}\int_{0}^{b}g(x)\,dx
=\frac{1}{w(b\del/2)}\int_{0}^{b\del/2}f(x)\,dx$,
we get $\limsup_{l(b)\ras\infty}
\bigl\|\frac{1}{w(b)}\int_{0}^{b}f(x)\,dx\bigr\|\leq2^{3d}\eps$.
Since this is true for any positive $\eps$,
$\lim_{l(b)\ras\infty}\frac{1}{w(b)}\int_{0}^{b}f(x)\,dx=0$.%
\endproof
\tsubsection{s-mulMUsup}{A lim-limsup version for uniform averages}
\theorem{P-mulMUsup}{}{}
Let $f\col\Rp^{d}\ra V$ be a bounded measurable function satisfying, for some $L\in V$,
\equ{
\esslim_{t\ras 0^{+}}\limsup_{l(b-a)\ras\infty}
\Bigl\|\frac{1}{w(b-a)}\sum_{n>a}^{b}f(nt)-L\Bigr\|=0.
}
Then $\lim_{l(b-a)\ras\infty}\frac{1}{w(b-a)}\int_{a}^{b}f(x)\,dx=L$.
\endtheorem
\lemma{P-asMUsup}{}{}
Let $(v_{n})$ be a bounded $\N^{d}$-sequence in $V$
and let $(b_{k})$, $(\alf_{k})$, and $(\bet_{k})$ be sequences in $\Rp^{d}$
such that $0<\bet_{k}-\alf_{k}\leq b_{k}$ for all $k$ and $l(b_{k})\ras\infty$.
Then $\limsup_{k\ras\infty}\bigl\|\frac{1}{w(b_{k})}
\sum_{n>\alf_{k}}^{\bet_{k}}v_{n}\bigr\|
\leq\limsup_{l(b-a)\ras\infty}\bigl\|\frac{1}{w(b-a)}\sum_{n>a}^{b}v_{n}\bigr\|$.
\endlemma
\proof{}
Assume that $\|v_{n}\|\leq 1$ for all $n$.
Let $s=\limsup_{l(b-a)\ras\infty}\bigl\|\frac{1}{w(b-a)}\sum_{n>a}^{b}v_{n}\bigr\|$.
Let $\eps>0$.
Let $B>1$ be such that
$\bigl\|\frac{1}{w(b-a)}\sum_{n>a}^{b}v_{n}\bigr\|<s+\eps$
whenever $l(b-a)>B$.
Let $K$ be such that $l(b_{k})>(B+1)/\eps$ for all $k>K$.
Then for any $k>K$,
if $l(\bet_{k}-\alf_{k})>B$,
then 
$\bigl\|\frac{1}{w(b_{k})}\sum_{n>\alf_{k}}^{\bet_{k}}v_{n}\bigr\|
\leq\bigl\|\frac{1}{w(\bet_{k}-\alf_{k})}\sum_{n>\alf_{k}}^{\bet_{k}}v_{n}\bigr\|
<s+\eps$,
and if $l(\bet_{k}-\alf_{k})\leq B$,
then also $\bigl\|\frac{1}{w(b_{k})}\sum_{n>\alf_{k}}^{\bet_{k}}v_{n}\bigr\|
\leq\frac{l(\bet_{k}-\alf_{k})+1}{l(b_{k})}<\eps\leq s+\eps$.
So, $\limsup_{k\ras\infty}\bigl\|\frac{1}{w(b_{k})}
\sum_{n>\alf_{k}}^{\bet_{k}}v_{n}\bigr\|\leq s$.%
\endproof
\proof{of \rfr{P-mulMUsup}}
We may assume that $L=0$.
Fix $\eps>0$.
Find $\del\in\Rp^{d}$, $\del>0$ such that
$\limsup_{l(b-a)\ras\infty}\bigl\|\frac{1}{w(b-a)}\sum_{n>a}^{b}f(nt)\bigr\|<\eps$
for \ae\ $t\in[0,\del]$,
and define $g(x)=f\bigl(\frac{\del}{2}x\bigr)$, $x\in\Rp^{d}$.
Then $\limsup_{l(b-a)\ras\infty}\bigl\|\frac{1}{w(b-a)}\sum_{n>a}^{b}g(nt)\bigr\|<\eps$
for \ae\ $t\in[0,2]^{d}$.

Let $a,b\in\Rp^{d}$, $a\geq 0$, $b\geq a+\bar{1}$.
For any $n\in\N^{d}$,
$\frac{1}{w(n)}\int_{a}^{b}g(x)\,dx=\int_{a/n}^{b/n}g(nt)\,dt$.
Adding these equalities for all $n\in(b/2,b]$,
and taking into account that $b/n<\bar{2}$ for $n>b/2$,
we get $\lam\int_{a}^{b}g(x)\,dx
=\int_{0}^{\bar{2}}\sum_{n>\alf(a,b,t)}^{\bet(a,b,t)}g(nt)\,dt$,
where $\lam=\sum_{n>b/2}^{b}\frac{1}{w(n)}\geq\frac{1}{2^{d}}$,
and for every $t\in(0,2]^{d}$, 
$\alf(a,b,t)=\max\{b/2,a/t\}$ and $\bet(a,b,t)=\min\{b,b/t\}$.
Thus, 
\equn{f-muldmnunfsup}{
\Bigl\|\frac{1}{w(b-a)}\int_{a}^{b}g(x)\,dx\Bigr\|
\leq 2^{d}\Bigl\|\int_{0}^{\bar{2}}g_{a,b}(t)\,dt\Bigr\|
}
where $g_{a,b}(t)=\frac{1}{w(b-a)}\sum_{n>\alf(a,b,t)}^{\bet(a,b,t)}g(nt)$, 
$t\in(0,2]^{d}$.

We will now show that the functions $g_{a,b}$, 
for $a,b\in\Rp^{d}$, $b\geq a+1$,
are uniformly bounded.
Let us assume that $\sup_{x\in\Rp^{d}}\|g(x)\|=\sup_{x\in\Rp^{d}}\|f(x)\|\leq 1$.
If $a\leq b/2$, then $b-a\geq b/2$,
and since $\bet(a,b,t)-\alf(a,b,t)\leq b/2$,
for any $t\in(0,2]^{d}$ we have $\|g_{a,b}(t)\|\leq\frac{1}{w(b/2)}w(b/2+\bar{1})\leq 3^{d}$.
If $a>b/2$, then for any $t=(t_{1},\ld,t_{d})\in(0,2]^{d}$ 
with $t_{i}<1/2$ for some $i$ 
we have $\alf(a,b,t)_{i}\geq a_{i}/t_{i}\geq 2a_{i}>b_{i}\geq\bet(a,b,t)_{i}$,
so $f_{a,b}(t)=0$;
and since, for any $t\in(0,2]^{d}$, $\bet(a,b,t)-\alf(a,b,t)\leq(b-a)/t$,
we have $\|f_{a,b}(t)\|\leq\frac{1}{w(b-a)}w((b-a)/t+\bar{1})\leq 3^{d}$ 
for all $t\in[1/2,2]^{d}$.

For \ae\ $t\in(0,2]^{d}$,
since $\bet(a,b,t)-\alf(a,b,t)\leq(b-a)/t$ for all $a,b\in\Rp^{d}$, $b\geq a+\bar{1}$, 
by \rfr{P-asMUsup},
\equ{
\limsup_{l(b-a)\ras\infty}\|g_{a,b}(t)\|
\leq\frac{1}{w(t)}\limsup_{l(b-a)\ras\infty}
\Bigl\|\frac{1}{w((b-a)/t)}\sum_{n>\alf(a,b,t)}^{\bet(a,b,t)}g(nt)\Bigr\|
<\frac{\eps}{w(t)}.
}
Since also $\limsup_{l(b-a)\ras\infty}\|g_{a,b}(t)\|\leq 3^{d}$,
we obtain that 
\equ{
\int_{0}^{\bar{2}}\limsup_{l(b-a)\ras\infty}\|g_{a,b}(t)\|\,dt
\leq\int_{[0,2]^{d}\sm[\eps,2]^{d}}3^{d}\,dt
+\int_{[\eps,2]^{d}}\frac{\eps}{w(t)}\,dt=c_{\eps},
}
where $c_{\eps}\leq d\eps2^{d-1}3^{d}+\eps(\log(2/\eps))^{d}$.
Hence, by \rfr{P-dctsup},
$\limsup_{l(b-a)\ras\infty}\bigl\|\int_{0}^{\bar{2}}g_{a,b}(t)\,dt\bigr\|\leq c_{\eps}$.
So, by \frfr{f-muldmnunfsup}, 
$\limsup_{l(b-a)\ras\infty}
\bigl\|\frac{1}{w(b-a)}\int_{a}^{b}g(x)\,dx\bigr\|\leq 2^{d}c_{\eps}$.
Since for any $0<a<b$,
$\frac{1}{w(b-a)}\int_{a}^{b}g(x)\,dx
=\frac{1}{w(b\del/2-a\del/2)}\int_{a\del/2}^{b\del/2}f(x)\,dx$,
we get 
$\limsup_{l(b-a)\ras\infty}\bigl\|\frac{1}{w(b-a)}\int_{a}^{b}f(x)\,dx\bigr\|
\leq 2^{d}c_{\eps}$.
Since this is true for any positive $\eps$ and $c_{\eps}\ras 0$ as $\eps\ras 0^{+}$,
we obtain that $\lim_{l(b-a)\ras\infty}\frac{1}{w(b-a)}\int_{a}^{b}f(x)\,dx=0$.%
\endproof
\section{S-Folner}{Two-sided limits and limits with respect to F{\o}lner sequences}
\tsubsection{s-mulZ}{Two-sided multiparameter limits}
We will now pass from the $(\N^{d},\R_{+}^{d})$ setup
to the $(\Z^{d},\R^{d})$ setup.
We adapt the notation introduced above to this new situation:
for $a,b\in\R^{d}$, $a=(a_{1},\ld,a_{d})$, $b=(b_{1},\ld,b_{d})$, 
we write $a\leq b$ if $a_{i}\leq b_{i}$ for all $i=1,\ld,d$,
and $a<b$ if $a_{i}<b_{i}$ for all $i$.
When writing $l(b)$ or $w(b)$, we will always assume that $b>0$.
As before, under $\sum_{n>a}^{b}v_{n}$ we understand
$\sum_{\sdup{n\in\Z^{d}}{a<n\leq b}}v_{n}$,
and under $\int_{a}^{b}v(x)\,dx$ we understand $\int_{a\leq x\leq b}v(x)\,dx$.

\rfr{P-addMC} clearly implies:
\theorem{P-addMCZ}{}{}
Let $f\col\R^{d}\ra V$ be a bounded measurable function
such that the limit 
$$
A_{t}=\lim_{l(b)\ras\infty}\frac{1}{w(2b)}\sum_{n>-b}^{b}f(t+n)
$$
exists for \ae\ $t\in[0,1]^{d}$.
Then $\lim_{l(b)\ras\infty}\frac{1}{w(2b)}\int_{-b}^{b}f(x)\,dx$ also exists
and is equal to $\int_{[0,1]^{d}}A_{t}\,dt$.
\endtheorem

\rfr{P-addMU} can also be easily adapted to the $\R^{d}$ case:
\theorem{P-addMUZ}{}{}
Let $f\col\R^{d}\ra V$ be a bounded measurable function
such that the limit 
$$
A_{t}=\lim_{l(b-a)\ras\infty}\frac{1}{w(b-a)}\sum_{n>a}^{b}f(t+n)
$$
exists for \ae\ $t\in[0,1]^{d}$.
Then $\lim_{l(b-a)\ras\infty}\frac{1}{w(b-a)}\int_{a}^{b}f(x)\,dx$ also exists
and is equal to $\int_{[0,1]^{d}}A_{t}\,dt$.
\endtheorem
\npar
The derivation of \rfr{P-addMUZ} from \rfr{P-addMU} is based on the following fact:
\lemma{P-UtoZ}{}{}
For any $s=(s_{1},\ld,s_{d})\in\{+,-\}^{d}=S$,
let $\R^{d}_{s}=\R_{s_{1}}\times\ld\times\R_{s_{d}}$.
Let $f\col\R^{d}\ra V$ be a bounded function and let $L$ be an element of $V$
such that for any $s\in S$,
$\lim_{\sdup{a,b\in\R_{s}}{l(b-a)\ras\infty}}\frac{1}{w(b-a)}\int_{a}^{b}f(x)\,dx=L$.
Then $\lim_{l(b-a)\ras\infty}\frac{1}{w(b-a)}\int_{a}^{b}f(x)\,dx=L$.
\endlemma
\proof{}
We will assume that $\sup|f|\leq 1$ and that $L=0$.
Given $\eps>0$, find $l\in\R$ such that 
$\bigl\|\frac{1}{w(b-a)}\int_{a}^{b}f(x)\,dx\bigr\|<\eps$
whenever $a,b\in\R^{d}_{s}$ for some $s\in S$ and $l(b-a)\geq l$.

Now let $a,b\in\R^{d}$, $b>a$ and $l(b-a)>l/\eps$.
Let $a=(a_{1},\ld,a_{d})$ and $b=(b_{1},\ld,b_{d})$.
For each $i$ such that $a_{i}<0<b_{i}$,
partition the interval $[a_{i},b_{i}]$ into subintervals $[a_{i},0]$ and $[0,b_{i}]$,
and thus partition the $d$-dimensional interval $[a,b]=\{x:a\leq x\leq b\}$
into $\leq 2^{d}$ $d$-dimensional subintervals $[p_{j},q_{j}]$ 
such that for each $j$, $[p_{j},q_{j}]\sle\R_{s}^{d}$ for some $s\in S$.
Then, for each $j$, 
if $l(q_{j}-p_{j})\geq l$,
then $\bigl\|\frac{1}{w(b-a)}\int_{p_{j}}^{q_{j}}f(x)\,dx\bigr\|
\leq\bigl\|\frac{1}{w(q_{j}-p_{j})}\int_{p_{j}}^{q_{j}}f(x)\,dx\bigr\|<\eps$,
and if $l(q_{j}-p_{j})<l$,
then $\bigl\|\frac{1}{w(b-a)}\int_{p_{j}}^{q_{j}}f(x)\,dx\bigr\|
\leq\frac{w(q_{j}-p_{j})}{w(b-a)}\leq\frac{l(q_{j}-p_{j})}{l(b-a)}<\frac{l}{l/\eps}=\eps$;
so, $\bigl\|\frac{1}{w(b-a)}\int_{a}^{b}f(x)\,dx\bigr\|
=\sum_{j}\bigl\|\frac{1}{w(b-a)}\int_{p_{j}}^{q_{j}}f(x)\,dx\bigr\|<2^{d}\eps$.%
\endproof

In the case $f$ is a real-valued function,
the same proof gives a stronger result:
\lemma{P-UtoSZ}{}{}
For any bounded measurable function $f\col\R^{d}\ra\R$,
$$
\liminf_{l(b-a)\ras\infty}\frac{1}{w(b-a)}\int_{a}^{b}f(x)\,dx
=\min_{s\in S}\left\{
\liminf_{\sdup{a,b\in\R_{s}}{l(b-a)\ras\infty}}\frac{1}{w(b-a)}\int_{a}^{b}f(x)\,dx\right\}
$$
and
$$
\limsup_{l(b-a)\ras\infty}\frac{1}{w(b-a)}\int_{a}^{b}f(x)\,dx
=\max_{s\in S}\left\{
\limsup_{\sdup{a,b\in\R_{s}}{l(b-a)\ras\infty}}\frac{1}{w(b-a)}\int_{a}^{b}f(x)\,dx\right\}.
$$
\endlemma

\rfr{P-UtoSZ} allows us to derive the ``two-sided'' version 
of Theorems~\rfrn{P-addMSC} and \rfrn{P-addMSU}:
\theorem{P-addMSCZ}{}{}
If $f\col\R^{d}\ra\R$ is a bounded measurable function, then
$$
\liminf_{l(b)\ras\infty}\frac{1}{w(2b)}\int_{-b}^{b}f(x)\,dx
\geq\int_{[0,1]^{d}}\liminf_{l(b)\ras\infty}\frac{1}{w(2b)}\sum_{n>-b}^{b}f(t+n)\,dt
$$
and
$$
\limsup_{l(b-a)\ras\infty}\frac{1}{w(2b)}\int_{-b}^{b}f(x)\,dx
\leq\int_{[0,1]^{d}}\limsup_{l(b)\ras\infty}\frac{1}{w(2b)}\sum_{n>-b}^{b}f(t+n)\,dt.
$$ 
\endtheorem
\theorem{P-addMSUZ}{}{}
If $f\col\R^{d}\ra\R$ is a bounded measurable function, then
$$
\liminf_{l(b-a)\ras\infty}\frac{1}{w(b-a)}\int_{a}^{b}f(x)\,dx
\geq\int_{[0,1]^{d}}\liminf_{l(b-a)\ras\infty}\frac{1}{w(b-a)}\sum_{n>a}^{b}f(t+n)\,dt
$$
and
$$
\limsup_{l(b-a)\ras\infty}\frac{1}{w(b-a)}\int_{a}^{b}f(x)\,dx
\leq\int_{[0,1]^{d}}\limsup_{l(b-a)\ras\infty}\frac{1}{w(b-a)}\sum_{n>a}^{b}f(t+n)\,dt.
$$ 
\endtheorem

For $a,b\in\R^{d}$ with $a<0<b$,
let us call ``the interval'' $(a,b)=\{t\in\R^{d}:a<t<b\}$
{\it a $P$-neighborhood of $0$ in $\R^{d}$},
and for $b\in\Rp^{d}$ with $b>0$
let us call ``the interval'' $[0,b)=\{t\in\Rp^{d}:t<b\}$
{\it a $P$-neighborhood of $0$ in $\Rp^{d}$}.

The ``multiplicative'' theorems for $\Z^{d}$ and $\R^{d}$-actions take the following form:
\theorem{P-mulMCZ}{}{}
Let a bounded measurable function $f\col\R^{d}\ra V$ 
be such that for some P-neighborhood of $0$ in $\R^{d}$, for \ae\ $t\in P$
the limit $L_{t}=\lim_{l(b)\ras\infty}\frac{1}{w(2b)}\sum_{n>-b}^{b}f(nt)$ exists.
Then $L_{t}=\const=L$ \ae\ on $P$
and $\lim_{l(b)\ras\infty}\frac{1}{w(2b)}\int_{-b}^{b}f(x)\,dx=L$.
\endtheorem
\theorem{P-mulMUZ}{}{}
Let a bounded measurable function $f\col\R^{d}\ra V$ 
be such that for some P-neighborhood of $0$ in $\R^{d}$, for \ae\ $t\in P$
the limit $L_{t}=\lim_{l(b-a)\ras\infty}\frac{1}{w(b-a)}\sum_{n>a}^{b}f(nt)$ exists.
Then $L_{t}=\const=L$ \ae\ on $P$
and $\lim_{l(b-a)\ras\infty}\frac{1}{w(b-a)}\int_{a}^{b}f(x)\,dx=L$.
\endtheorem
\theorem{P-mulMSCZ}{}{}
For any bounded measurable function $f\col\R^{d}\ra\R$ and any $c\in\Rp^{d}$, $c>0$,
$$
\liminf_{l(b)\ras\infty}\frac{1}{w(2b)}\int_{-b}^{b}f(x)\,dx
\geq\frac{1}{w(c)}\int_{0}^{c}\liminf_{l(b)\ras\infty}\frac{1}{w(2b)}
\sum_{n>-b}^{b}f(nt)\,dt
$$
and
$$
\limsup_{l(b)\ras\infty}\frac{1}{w(2b)}\int_{-b}^{b}f(x)\,dx
\leq\frac{1}{w(c)}\int_{0}^{c}\limsup_{l(b)\ras\infty}\frac{1}{w(2b)}
\sum_{n>-b}^{b}f(nt)\,dt.
$$
\endtheorem
\theorem{P-mulMSUZ}{}{}
For any bounded measurable function $f\col\R^{d}\ra\R$ and any $c\in\Rp^{d}$, $c>0$,
$$
\liminf_{l(b-a)\ras\infty}\frac{1}{w(b-a)}\int_{a}^{b}f(x)\,dx
\geq\frac{1}{w(c)}\int_{0}^{c}\liminf_{l(b-a)\ras\infty}\frac{1}{w(b-a)}
\sum_{n>a}^{b}f(nt)\,dt
$$
and
$$
\limsup_{l(b-a)\ras\infty}\frac{1}{w(b-a)}\int_{a}^{b}f(x)\,dx
\leq\frac{1}{w(c)}\int_{0}^{c}\limsup_{l(b-a)\ras\infty}\frac{1}{w(b-a)}
\sum_{n>a}^{b}f(nt)\,dt.
$$
\endtheorem
\theorem{P-mulMCZsup}{}{}
Let $f\col\R^{d}\ra V$ be a bounded measurable function satisfying
\equ{
\esslim_{t\ras 0^{+}}\limsup_{l(b)\ras\infty}
\Bigl\|\frac{1}{w(2b)}\sum_{n>-b}^{b}f(nt)\Bigr\|=0.
}
Then $\lim_{l(b)\ras\infty}\frac{1}{w(2b)}\int_{-b}^{b}f(x)\,dx=0$.
\endtheorem
\theorem{P-mulMUZsup}{}{}
Let $f\col\R^{d}\ra V$ be a bounded measurable function satisfying, for some $L\in V$,
\equ{
\esslim_{t\ras 0^{+}}\limsup_{l(b-a)\ras\infty}
\Bigl\|\frac{1}{w(b-a)}\sum_{n>a}^{b}f(nt)-L\Bigr\|=0.
}
Then $\lim_{l(b-a)\ras\infty}\frac{1}{w(b-a)}\int_{a}^{b}f(x)\,dx=L$.
\endtheorem
\tsubsection{s-Folner}{Limits with respect to an arbitrary F{\o}lner sequence}
Let us denote by $w$ the standard Lebesgue measure on $\R^{d}$
(this agrees with the notation used in the previous sections).
{\it A F{\o}lner sequence in $\R^{d}$}
is a sequence $(\Phi_{N})_{N=1}^{\infty}$ of subsets of finite mesure
such that for any $y\in\R^{d}$,
$\frac{w(\Phi_{N}\tri(\Phi_{N}+y))}{w(\Phi_{N})}\ra 0$ as $N\ras\infty$.
\lemma{P-Fol}{}{}
Let $f\col\R^{d}\ra V$ be a bounded measurable function
with the property that
$\lim_{l(b-a)\ras\infty}\dsc\frac{1}{w(b-a)}\int_{a}^{b}f(x)\,dx=L\in V$.
Then for any F{\o}lner sequence $(\Phi_{N})$ in $\R^{d}$,
$\lim_{N\ras\infty}\frac{1}{w(\Phi_{N})}\int_{\Phi_{N}}f(x)\,dx=L$.
\endlemma
\proof{}
We will assume that $L=0$ and that $\sup|f|\leq 1$.
Let $\eps>0$,
and let $Q$ be a $d$-dimensional interval $\bigl\{x\in\R^{d}:0\leq x\leq c\bigr\}$
with $l(c)$ large enough
so that $\bigl\|\frac{1}{w(Q)}\int_{Q+y}f(x)\,dx\bigr\|<\eps$
for any $y\in\R^{d}$.
Let $(\Phi_{N})$ be a F{\o}lner sequence in $\R^{d}$.
For any $y\in Q$ we have
\lequ{
2\geq
\frac{1}{w(\Phi_{N})}\Bigl\|\int_{\Phi_{N}}f(x+y)\,dx-\int_{\Phi_{N}}f(x)\,dx\Bigr\|
\leq
\frac{1}{w(\Phi_{N})}\Bigl(\Bigl\|\int_{(\Phi_{N}+y)\sm\Phi_N}f(x)\,dx\Bigr\|
+\Bigl\|\int_{\Phi_{N}\sm(\Phi_N+y)}f(x)\,dx\Bigr\|\Bigr)
\-\\\-
\leq\frac{w(\Phi_{N}\tri(\Phi_{N}+y))}{w(\Phi_{N})}\ra 0
\hbox{ as $N\ras\infty$}.
}
So, by \rfr{P-dct},
$$
\frac{1}{w(Q)w(\Phi_{N})}\int_{Q}
\Bigl(\int_{\Phi_{N}}f(x+y)\,dx-\int_{\Phi_{N}}f(x)\,dx\Bigr)\,dy\dsc\ra 0 
\hbox{ as $N\ras\infty$.}
$$
But 
$\frac{1}{w(Q)w(\Phi_{N})}\int_{Q}\int_{\Phi_{N}}f(x)\,dx\,dy
=\frac{1}{w(\Phi_{N})}\int_{\Phi_{N}}f(x)\,dx$ for all $N$,
whereas
\lequ{
\Bigl\|\frac{1}{w(Q)w(\Phi_{N})}\int_{Q}\int_{\Phi_{N}}f(x+y)\,dx\,dy\Bigr\|
\leq\frac{1}{w(\Phi_{N})}\int_{\Phi_{N}}
\Bigl\|\frac{1}{w(Q)}\int_{Q}f(x+y)\,dy\Bigr\|\,dx
\-\\\-
=\frac{1}{w(\Phi_{N})}\int_{\Phi_{N}}
\Bigl\|\frac{1}{w(Q)}\int_{Q-x}f(y)\,dy\Bigr\|\,dx
\leq\frac{1}{w(\Phi_{N})}\int_{\Phi_{N}}\eps\,dx=\eps
}
for all $N$.
Hence, $\limsup_{N\ras\infty}
\bigl\|\frac{1}{w(\Phi_{N})}\int_{\Phi_{N}}f(x)\,dx\bigr\|<\eps$.%
\endproof
\rfr{P-Fol} allows to strengthen Theorems~\rfrn{P-addMUZ},
\rfrn{P-mulMUZ}, and \rfrn{P-mulMUZsup}:
\theorem{P-addF}{}{}
Let $f\col\R^{d}\ra V$ be a bounded measurable function
such that the limit 
$$
A_{t}=\lim_{l(b-a)\ras\infty}\frac{1}{w(b-a)}\sum_{n>a}^{b}f(t+n)
$$
exists for \ae\ $t\in[0,1]^{d}$.
Then for any F{\o}lner sequence $(\Phi_{N})$ in $\R^{d}$,
$\lim_{N\ras\infty}\frac{1}{w(\Phi_{N})}\int_{\Phi_{N}}f(x)\,dx=\int_{[0,1]^{d}}A_{t}\,dt$.
\endtheorem
\theorem{P-mulF}{}{}
Let $f\col\R^{d}\ra V$ be a bounded measurable function
such that for some $P$ neighborhood $P$ of $0$ in $\R^{d}$, 
for \ae\ $t\in P$
the limit $L_{t}=\lim_{l(b-a)\ras\infty}\frac{1}{w(b-a)}\sum_{n>a}^{b}f(nt)$ exists.
Then $L_{t}=\const=L$ \ae\ on $P$
and for any F{\o}lner sequence $(\Phi_{N})$ in $\R^{d}$,
$\lim_{N\ras\infty}\frac{1}{w(\Phi_{N})}\int_{\Phi_{N}}f(x)\,dx=L$.
\endtheorem
\theorem{P-mulFsup}{}{}
Let $f\col\R^{d}\ra V$ be a bounded measurable function satisfying, for some $L\in V$,
\equ{
\esslim_{t\ras 0^{+}}\limsup_{l(b-a)\ras\infty}
\Bigl\|\frac{1}{w(b-a)}\sum_{n>a}^{b}f(nt)-L\Bigr\|=0.
}
Then for any F{\o}lner sequence $(\Phi_{N})$ in $\R^{d}$,
$\lim_{N\ras\infty}\frac{1}{w(\Phi_{N})}\int_{\Phi_{N}}f(x)\,dx=L$.
\endtheorem
In the case $f$ is a real-valued function we can get the following version of \rfr{P-Fol}:
\lemma{P-SFol}{}{}
For any bounded measurable function $f\col\R^{d}\ra\R$,
for any F{\o}lner sequence $(\Phi_{N})$ in $\R^{d}$,
$$
\liminf_{N\ras\infty}\frac{1}{w(\Phi_{N})}\int_{\Phi_{N}}f(x)\,dx
\geq\liminf_{l(b-a)\ras\infty}\dsc\frac{1}{w(b-a)}\int_{a}^{b}f(x)\,dx
$$
and
$$
\limsup_{N\ras\infty}\frac{1}{w(\Phi_{N})}\int_{\Phi_{N}}f(x)\,dx
\leq\limsup_{l(b-a)\ras\infty}\dsc\frac{1}{w(b-a)}\int_{a}^{b}f(x)\,dx.
$$
\endlemma
Using \rfr{P-SFol}, we may also strengthen Theorems~\rfrn{P-addMSUZ} and \rfrn{P-mulMSUZ}:
\theorem{P-addSF}{}{}
For any bounded measurable function $f\col\R^{d}\ra\R$
and any F{\o}lner sequence $(\Phi_{N})$ in $\R^{d}$,
$$
\liminf_{N\ras\infty}\frac{1}{w(\Phi_{N})}\int_{\Phi_{N}}f(x)\,dx
\geq\int_{[0,1]^{d}}\liminf_{l(b-a)\ras\infty}\frac{1}{w(b-a)}\sum_{n>a}^{b}f(t+n)\,dt.
$$
and
$$
\limsup_{N\ras\infty}\frac{1}{w(\Phi_{N})}\int_{\Phi_{N}}f(x)\,dx
\leq\int_{[0,1]^{d}}\limsup_{l(b-a)\ras\infty}\frac{1}{w(b-a)}\sum_{n>a}^{b}f(t+n)\,dt.
$$
\endtheorem
\theorem{P-mulSF}{}{}
For any bounded measurable function $f\col\R^{d}\ra\R$, any $c\in\Rp^{d}$, $c>0$,
and any F{\o}lner sequence $(\Phi_{N})$ in $\R^{d}$,
$$
\liminf_{N\ras\infty}\frac{1}{w(\Phi_{N})}\int_{\Phi_{N}}f(x)\,dx
\geq\frac{1}{w(c)}\int_{0}^{c}\liminf_{l(b-a)\ras\infty}\frac{1}{w(b-a)}
\sum_{n>a}^{b}f(nt)\,dt
$$
and
$$
\limsup_{N\ras\infty}\frac{1}{w(\Phi_{N})}\int_{\Phi_{N}}f(x)\,dx
\leq\frac{1}{w(c)}\int_{0}^{c}\limsup_{l(b-a)\ras\infty}\frac{1}{w(b-a)}
\sum_{n>a}^{b}f(nt)\,dt.
$$
\endtheorem
\section{S-Density}{Density of sets and convergence in density}

We will now formulate some special cases of the theorems above.
For a set $S\sle\N^{d}$, 
{\it the density\/} of $S$ is
$\D(S)=\lim_{l(b)\ras\infty}\frac{1}{w(b)}\bigl|S\cap(0,b]\bigr|$,
if it exists;
for a measurable set $S\sle\Rp^{d}$, 
{\it the density\/} of $S$ is
$\D(S)=\lim_{l(b)\ras\infty}\frac{1}{w(b)}w(S\cap[0,b])$, if it exists.
(As before, $w$ stands for the standard Lebesgue measure on $\R^{d}$).
{\it The lower density\/} $\Du(S)$ and {\it the upper density\/} $\Do(S)$
of a set $S\sle\N^{d}$ or $S\sle\Rp^{d}$
are defined as the $\liminf$ and, respectively, the $\limsup$ of the above expressions.

Taking $f=1_{S}$ in Theorems~\rfrn{P-addMSC}, \rfrn{P-addMC} 
and in Theorems~\rfrn{P-mulMSC}, \rfrn{P-mulMC},
we get, respectively:

\theorem{P-addDC}{}{}
Let $S$ be a measurable subset of $\Rp^{d}$,
and for each $t\in[0,1]^{d}$ let $S_{t}=\bigl\{n\in\N^{d}:t+n\in S\bigr\}$.
Then $\Du(S)\geq\int_{[0,1]^{d}}\Du(S_{t})\,dt$ 
and $\Do(S)\leq\int_{[0,1]^{d}}\Do(S_{t})\,dt$.
If $\D(S_{t})$ exists for \ae\ $t\in[0,1]^{d}$,
then $\D(S)$ also exists and equals $\int_{[0,1]^{d}}\D(S_{t})\,dt$.
\endtheorem

\theorem{P-mulDC}{}{}
Let $S$ be a measurable subset of $\Rp^{d}$,
and for each $t\in\Rp^{d}$ let $S_{t}=\bigl\{n\in\N^{d}:nt\in S\bigr\}$.
Then for any $c\in\Rp^{d}$, $c>0$,
one has $\Du(S)\geq\int_{[0,c]}\Du(S_{t})\,dt$ 
and $\Do(S)\leq\int_{[0,c]}\Do(S_{t})\,dt$.
If $\D(S_{t})$ exists for \ae\ $t$ in a P-neighborhood $P$ of $0$ in $\Rp^{d}$,
then $\D(S(t))=\const=D$ for \ae\ $t\in P$ and $\D(S)=D$.
\endtheorem

{\it The uniform {\rm(or {\it Banach})} density\/} of a set $S\sle\N^{d}$ is
$\UD(S)=\lim_{l(b-a)\ras\infty}\frac{1}{w(b-a)}\bigl|S\cap(a,b]\bigr|$,
if it exists;
for a measurable set $S\sle\Rp^{d}$ {\it the uniform density\/} of $S$ is
$\UD(S)=\lim_{l(b-a)\ras\infty}\frac{1}{w(b-a)}w(S\cap[a,b])$, if it exists.
(And it follows from (an $\Rp^{d}$-version) of \rfr{P-Fol} 
that for $S\sle\Rp^{d}$,
if $\UD(S)$ exists, then for any F{\o}lner sequence $(\Phi_{N})$ in $\Rp^{d}$,
$\lim_{N\ras\infty}\frac{1}{|\Phi_{N}|}w(S\cap\Phi_{N})=\UD(S)$.)
{\it The lower uniform density\/} $\UDu(S)$
and {\it the upper uniform density\/} $\UDo(S)$
of a set $S\sle\N^{d}$ or $S\sle\Rp^{d}$
is the $\liminf$ and, respectively, the $\limsup$ of the above expressions.

From Theorems~\rfrn{P-addMSU}, \rfrn{P-addMU}, \rfrn{P-mulMSU}, and \rfr{P-mulMU}
we get, respectively:

\theorem{P-addDU}{}{}
Let $S$ be a measurable subset of $\Rp^{d}$,
and for each $t\in[0,1]^{d}$ let $S_{t}=\bigl\{n\in\N^{d}:n+t\in S\bigr\}$.
Then $\UDu(S)\geq\int_{[0,1]^{d}}\UDu(S_{t})\,dt$ 
and $\UDo(S)\leq\int_{[0,1]^{d}}\UDo(S_{t})\,dt$.
If $\UD(S_{t})$ exists for \ae\ $t\in[0,1]^{d}$,
then $\UD(S)$ also exists and equals $\int_{[0,1]^{d}}\UD(S_{t})\,dt$.
\endtheorem

\theorem{P-mulDU}{}{}
Let $S$ be a measurable subset of $\Rp^{d}$,
and for each $t\in\Rp^{d}$ let $S_{t}=\bigl\{n\in\N^{d}:nt\in S\bigr\}$.
Then, for any $c\in\Rp^{d}$, $c>0$, 
one has $\Du(S)\geq\frac{1}{w(c)}\int_{[0,c]}\Du(S_{t})\,dt$ 
and $\Do(S)\leq\frac{1}{w(c)}\int_{[0,c]}\Do(S_{t})\,dt$.
If $\UD(S_{t})$ exists for \ae\ $t$ in a P-neighborhood $P$ of $0$ in $\Rp^{d}$,
then $\UD(S(t))=\const=D$ in $P$ and $\UD(S)=D$.
\endtheorem


\npar
Of course, the ``two-sided'' versions of Theorems~\rfrn{P-addDC} -- \rfrn{P-mulDU},
where one deals with $\Z^{d}$-sequences and functions on $\R^{d}$
instead of $\N^{d}$-sequences and functions on $\Rp^{d}$,
are also true.

\vbreak{-9000}{2mm}
We will now bring two theorems 
that deal with limits in density instead of Ces\`{a}ro limits.
We say that an $\N^{d}$-sequence $(v_{n})$ in $V$
{\it converges in density to $L\in V$} if for any $\eps>0$, 
the set $S_{\eps}=\bigl\{n\in\N^{d}:\|v_{n}-L\|>\eps\bigr\}$ has zero density, 
$\D(S_{\eps})=0$,
and {\it converges to $L$ in uniform density\/}
if for any $\eps>0$, $\UD(S_{\eps})=0$. 
We say that a (measurable) function $f\col\Rp^{d}\ra V$
converges to $L\in V$ in density
if for any $\eps>0$,
the set $S_{\eps}=\bigl\{x\in\Rp^{d}:\|f(x)-L\|>\eps\bigr\}$ has zero density, 
$\D(S_{\eps})=0$,
and converges to $L$ in uniform density
if for any $\eps>0$, $\UD(S_{\eps})=0$. 
Applying Theorems~\rfrn{P-addDC} -- \rfrn{P-mulDU}
to the real-valued function $\|f(x)-L\|$
we obtain:
\theorem{P-addCD}{}{}
Let $f\col\Rp^{d}\ra V$ be a bounded measurable function
such that for some $L\in V$,
for \ae\ $t\in[0,1]^{d}$
the $\N^{d}$-sequence $f(n+t)$, $n\in\N^{d}$, converges to $L$ in density 
(respectively, in uniform density).
Then $f$ converges to $L$ in density
(respectively, in uniform density).
\endtheorem

\theorem{P-mulCD}{}{}
Let $f\col\Rp^{d}\ra V$ be a bounded measurable function
such that for some $L\in V$,
for \ae\ $t$ in a P-neighborhood of $0$ in $\Rp^{d}$
the $\N^{d}$-sequence $f(nt)$, $n\in\N^{d}$, converges to $L$ in density 
(respectively, in uniform density).
Then $f$ converges to $L$ in density
(respectively, in uniform density).
\endtheorem


\npar
Of course, the two-sided versions 
of Theorems~\rfrn{P-addCD} and \rfrn{P-mulCD} also hold.
\section{S-Appl}{Applications}
\tsubsection{s-charfac}{Characteristic factors for multiple averages along polynomials}

Let $X$ be a probability measure space;
we will always assume that $X$ is sufficiently regular so that $L^{1}(X)$ is separable.

Let $G$ be a group of measure preserving transformations of $X$
and let $g_{1}(n),\ld,g_{r}(n)$, $n\in\Z^{d}$, 
be ($d$-parameter) sequences of elements of $G$.
A factor $Z$ of the system $(X,G)$ is said to be {\it characteristic\/}
for the averages 
$\frac{1}{|\Psi_{N}|}\sum_{n\in\Psi_{N}}g_{1}(n)f_{1}\cd\ld\cd g_{r}(n)f_{r}$,
where $(\Psi_{N})$ is a F{\o}lner sequence in $\Z^{d}$,
if for any $f_{1},\ld,f_{r}\in L^{\infty}(X)$,
$$
\lim_{N\ras\infty}\frac{1}{|\Psi_{N}|}\sum_{n\in\Psi_{N}}
\Bigl(g_{1}(n)f_{1}\cd\ld\cd g_{r}(n)f_{r}
-g_{1}(n)E(f_{1}|Z)\cd\ld\cd g_{r}(n)E(f_{r}|Z)\Bigr)=0
$$
in $L^{1}(X)$
(where $E(f|Z)$ stands for the conditional expectation of $f$ with respect to $Z$).
An analogous notion can be introduced for averages
$\frac{1}{w(\Phi_{N})}\int_{\Phi_{N}}g_{1}(x)f_{1}\cd\ld\cd g_{r}(x)f_{r}\,dx$,
where $g_{1},\ld,g_{r}$ are functions $\R^{d}\ra G$
and $(\Phi_{N})$ is a F{\o}lner sequence in $\R^{d}$.

Let $T$ be an ergodic invertible measure preserving transformation of $X$.
{\it The $k$-th Host-Kra-Ziegler factor\/} $Z_{k}(T)$ of $(X,T)$
is the minimal characteristic factor for the averages 
$\frac{1}{|\Psi_{N}|}\sum_{n\in\Psi_{N}}
\prod_{\pus\neq\sig\sle\{0,\ld,k\}}T^{n_{\sig}}f_{\sig}$,
where $n_{\sig}=\sum_{i\in\sig}n_{i}$,
and $(\Psi_{N})$ are F{\o}lner sequences in $\Z^{k+1}$.
$Z_{k}(T)$ is the maximal factor of $(X,T)$
isomorphic to a $k$-step pro-nilmanifold
(an inverse limit of compact $k$-step nilmanifolds)
on which $T$ acts as a translation.
(See \brfr{HKo} and \brfr{Zo}.)
The factors $Z_{k}(T)$ turn out to be characteristic 
for any system of polynomial powers of $T$:
\theorem{P-chardis}{}{(\brfr{cpm})}
For any system of polynomials $p_{1},\ld,p_{r}\col\Z^{d}\ra\Z$
there exists $k\in\N$
such that for any measure preserving transformation of a probability measure space $X$,
$Z_{k}(T)$ is a characteristic factor for the averages
$\frac{1}{|\Phi_{N}|}\sum_{n\in\Phi_{N}}T^{p_{1}(n)}f_{1}\cd\ld\cd T^{p_{r}(n)}f_{r}$.
\endtheorem

It is easy to see (see, for example, \brfr{FrKr})
that if $S$ is another ergodic transformation of $X$ commuting with $T$
then for all $k$, $Z_{k}(S)=Z_{k}(T)$.
Thus, if $T$ is a family of pairwise commuting ergodic transformations of $X$,
we may denote by $Z_{k}(T)$ the $k$-th Host-Kra-Ziegler factor 
of any (and so, of every) element of $T$.
This allows one to generalize \rfr{P-chardis} in the following way:

\theorem{P-charmdis}{}{(\brfr{Johnson})}
For any finite system of polynomials $p_{i}\col\Z^{d}\ra\Z^{c}$, $i=1,\ld,r$,
there exists $k\in\N$ such that, 
given any totally ergodic%
\fnote{A group $G$ of measure preserving transformations of a measure space
is {\it totally ergodic\/}
if every nonidentical element of $G$ is totally ergodic}
discrete $c$-parameter commutative group $T^{m}$, $m\in\Z^{c}$,
of measure preserving transformations $T$ of a probability measure space $X$,
the factor $Z_{k}(T)$ is characteristic for the averages
$\frac{1}{|\Psi_{N}|}\sum_{n\in\Psi_{N}}
T^{p_{1}(n)}f_{1}\cd\ld\cd T^{p_{r}(n)}f_{r}$,
where $(\Psi_{N})$ are F{\o}lner sequences in $\Z^{d}$.
\endtheorem

Now let $T^{t}$, $t\in\R$, be a continuous $1$-parameter group 
of measure preserving transformations of $X$
and assume that it is ergodic on $X$.
Then for almost all (actually, for all but countably many) $t\in\R$ 
the transformation $T^{t}$ is ergodic,
so for any $k$, $Z_{k}(T^{t})$ coinside for \ae\ $t$;
we will denote this factor by $Z_{k}(T)$.
We can now prove the following fact
(obtained in \brfr{Potts} for non-uniform averages).

\theorem{P-charfac}{}{}
For any system of polynomials $p_{1},\ld,p_{r}\col\R^{d}\ra\R$
there exists $k\in\N$
such that for any continuous $1$-parameter group $T^{t}$, $t\in\R$,
of measure preserving transformations of a probability measure space $X$,
$Z_{k}(T)$ is a characteristic factor for the averages
$\frac{1}{w(\Phi_{N})}\int_{\Phi_{N}}T^{p_{1}(x)}f_{1}\cd\ld\cd T^{p_{r}(x)}f_{r}\,dx$.
\endtheorem

\proof{}
Given polynomials $p_{1},\ld,p_{r}$ on $\R^{d}$,
find monomials $q_{\lam}(x)=c_{\lam}x^{\alf_{\lam}}$, $\lam=1,\ld,\Lam$,
where $c_{\lam}\in\R$ and $\alf_{\lam}$ are multi-indices,
that are $\Q$-linearly independent
and such that each of the polynomials $p_{i}$ 
is a sum of the monomials $q_{\lam}$ with integer coefficients,
$p_{i}=\sum_{\lam=1}^{\Lam}b_{i,\lam}q_{\lam}$, $b_{i,\lam}\in\Z$.
Then for any $x\in\R^{d}$, any $n\in\Z^{d}$, and any $i$,
$T^{p_{i}(nx)}=\prod_{\lam=1}^{\Lam}T_{x,\lam}^{b_{i,\lam}n^{\alf_{\lam}}}$
where $T_{x,\lam}=T^{c_{\lam}x^{\alf_{\lam}}}$,
and since $T^{t}$ is ergodic for \ae\ $t\in\R$,
the $\Lam$-parameter group generated by the transformations $T_{x,\lam}$, $\lam=1,\ld,\Lam$,
satisfies the assumptions of \rfr{P-charmdis} for \ae\ $x\in\R^{d}$.
Find $k$ which, by \rfr{P-charmdis}, corresponds to the polynomials
$b_{i,\lam}n^{\alf_{\lam}}$, $i=1,\ld,r$, $\lam=1,\ld,\Lam$,
so that for \ae\ $x\in\R^{d}$,
$$
\lim_{N\ras\infty}\frac{1}{|\Psi_{N}|}\sum_{n\in\Psi_{N}}
\Bigl(T^{p_{1}(nx)}f_{1}\cd\ld\cd T^{p_{r}(nx)}f_{r}
-T^{p_{1}(nx)}E(f_{1}|Z_{k}(T))\cd\ld\cd T^{p_{r}(nx)}E(f_{r}|Z_{k}(T))\Bigr)=0
$$
for any F{\o}lner sequence $(\Psi_{N})$ in $\Z^{d}$.
Then by \rfr{P-mulF},
$$
\lim_{N\ras\infty}\frac{1}{w(\Phi_{N})}\int_{\Phi_{N}}
\Bigl(T^{p_{1}(x)}f_{1}\cd\ld\cd T^{p_{r}(x)}f_{r}
-T^{p_{1}(x)}E(f_{1}|Z_{k}(T))\cd\ld\cd T^{p_{r}(x)}E(f_{r}|Z_{k}(T))\Bigr)\,dx=0
$$
for any F{\o}lner sequence $(\Phi_{N})$ in $\R^{d}$,
which proves \rfr{P-charfac}.
\endproof
\tsubsection{s-nilorb}{Polynomial orbits on nilmanifolds}

Let $X$ be a topological space with a probability Borel measure $\mu$.
We say that a $d$-parameter sequence $g(n)$, $n\in\Z^{d}$,
{\it is well distributed with respect to $\mu$}
if for any $h\in C(X)$
and any F{\o}lner sequence $(\Psi_{N})$ in $\Z^{d}$
one has
$\lim_{N\ras\infty}\frac{1}{|\Psi_{N}|}\sum_{n\in\Psi_{N}}h(g(n))
=\int_{X}h\,d\mu$.
We also say that a measurable function $g(t)$, $t\in\R^{d}$, in $X$
is well distributed with respect to $\mu$
if for any $h\in C(X)$
and any F{\o}lner sequence $(\Phi_{N})$ in $\R^{d}$,
$\lim_{N\ras\infty}\frac{1}{w(\Phi_{N})}\int_{\Phi_{N}}h(g(t))\,dt
=\int_{X}h\,d\mu$.

The following proposition is an immediate corollary of \rfr{P-addF},
applied to the functions $h\comp g$, $h\in C(X)$.
\proposition{P-addWdis}{}{}
Let $X$ be a topological space
and let $g\col\R^{d}\ra X$ be a function
such that for \ae\ $t\in[0,1]^{d}$ the sequence $g(n+t)$, $n\in\Z^{d}$,
is well distributed in $X$ with respect to a probability Borel measure $\mu_{t}$.
Then $g$ is well distributed with respect to the measure
$\mu=\int_{[0,1]^{d}}\mu_{t}\,dt$.
\endproposition

From \rfr{P-mulF} we get:
\proposition{P-mulWdis}{}{}
Let $X$ be a compact Hausdorff space for which $C(X)$ is separable
and let $g\col\R^{d}\ra X$ be a function
such that for \ae\ $t$ in a P-neighborhood $P$ of $0$ in $\Rp^{d}$ 
the sequence $g(nt)$, $n\in\Z^{d}$,
is well distributed in $X$ with respect to a probability Radon measure $\mu_{t}$.
Then $\mu_{t}=\const=\mu$ for \ae\ $t\in P$
and $g$ is well distributed with respect to the measure $\mu$.
\endproposition
\proof{}
By \rfr{P-mulF}, apllied to the function $h\comp g$,
for any $h\in C(X)$ we have
$\mu_{t}(h)=\const=\mu(h)$ for \ae\ $t\in P$
and $\lim_{N\ras\infty}\frac{1}{w(\Phi_{N})}\int_{\Phi_{N}}h(g(t))\,dt=\mu(h)$
for any F{\o}lner sequence $(\Phi_{N})$ in $\R^{d}$.
Excluding those $t$ for which $\mu_{t}(h)\neq\mu(h)$ 
for all functions $h$ from a fixed countable subset of $C(X)$,
we obtain that $\mu_{t}=\const=\mu$ for \ae\ $t\in P$
and $g$ is well distributed with respect to $\mu$.
(The assumption that $\mu_{t}$ are Radon measures
allows us to identify them with continuous linear functionals on $C(X)$.)
\endproof

We will apply these propositions in the following situation.
Let $X$ be a compact {\it nilmanifold},
that is, a homogeneous space of a nilpotent Lie group $G$,
and let $g\col\R^{d}\ra X$ be a {\it polynomial mapping\/},
that is, $g(t)=a_{1}^{p_{1}(t)}\ld a_{k}^{p_{k}(t)}\om$, $t\in\R^{d}$,
where $a_{1},\ld,a_{k}\in G$, $p_{1},\ld,p_{k}$ are polynomials $\R^{d}\ra\R$,
and $\om\in X$.
Let $Y=\ovr{\{g(t),\ t\in\R^{d}\}}$.
It follows from a general result obtained in \brfr{Shah}
that $Y$ is a connected sub-nilmanifold of $X$
(that is, a closed subset of $X$ of the form $H\om$ 
where $H$ is a connected closed subgroup of $G$ and $\om\in X$),
and $g$ is uniformly distributed in $Y$ in the following sense:
for any $h\in C(Y)$,
$\lim_{R\ras\infty}\frac{1}{w(B_{R})}\int_{B_{R}}h(g(t))\,dw(t)
=\int_{Y}h\,d\mu$,
where $w$ is the Lebesgue measure on $\R^{d}$,
$B_{R}$, $R>0$, is the ball $\{t\in\R^{d}:|t|\leq R\}$,
and $\mu$ is the Haar measure on $Y$.
We would like to have a stronger result
which states that $g$ is not only uniformly distributed,
but is well distributed in $Y$.
A discrete analogue of this fact, which we will presently formulate,
was obtained in \brfr{mpn} and \brfr{ran},
but before formulating it 
we need to introduce some terminology.
We call a finite disjoint union of connected subnilmanifolds of $X$
{\it a FU~subnilmanifold}.
We say that an element $\om'$ of $X$ is {\it rational\/} 
with respect to an element $\om\in X$
if $\om'=a\om$ for some $a\in G$ such that $a^{m}\om=\om$ for some $m\in\N$.
We say that a subnilmanifold $Y$ of $X$
is {\it rational with respect to $\om$} 
if $Y$ contains an element $\om'$ rational with respect to $\om$.
(Then such elements $\om'$ are dense in $Y$.)
Finally, we say that a FU~subnilmanifold of $X$ is rational with respect to $\om$
if all connected components of $Y$ are subnilmanifolds rational with respect to $\om$.
\proposition{P-disorb}{}{(See \brfr{mpn} and \brfr{ran}.)}
Let $g$ be a $d$-parameter polynomial sequence in $X$,
that is, $g(n)=a_{1}^{p_{1}(n)}\ld a_{k}^{p_{k}(n)}\om$
where $a_{1},\ld,a_{k}\in G$, $p_{1},\ld,p_{k}$ are polynomials $\Z^{d}\ra\R$,
and $\om\in X$.
Then the closure $Y=\ovr{\{g(n),\ n\in\Z^{d}\}}$ of $g$
is a FU~subnilmanifold of $X$
rational with respect to the point $g(0)$.
If $Y$ is connected,
then the sequence $g(n)$, $n\in\Z^{d}$, is well distributed in $Y$
(with respect to the Haar measure on $Y$).
\endproposition

We may now use \rfr{P-mulF} or \rfr{P-addF}
to deduce from \rfr{P-disorb} its continuous analogue.
We will also need the following fact: 
\proposition{P-genorb}{}{(\brfr{orb}, Theorem~2.1)}
Let $M$ be a set
and let $\phi\col\R^{d}\times M\ra X$ be a mapping
such that for every $m\in M$,
$\phi(\cd,m)$ is a polynomial mapping $\R^{d}\ra X$,
and there exists $\om\in X$
such that for each $t\in\R^{d}$
the set $Y_{t}=\ovr{\phi(t,M)}$ 
is a FU~subnilmanifold of $X$ rational with respect to $\om$.
Then there exists a FU~subnilmanifold\/ $Y$ of\/ $X$ 
such that $Y_{t}\sle Y$ for all $t\in\R^{d}$
and $Y_{t}=Y$ for \ae\ $t\in\R^{d}$.
\endproposition

Now let $g\col\R^{d}\ra X$ be a polynomial mapping.
By \rfr{P-disorb},
the mapping $\phi\col\R^{d}\times\Z^{d}\ra X$
defined by $\phi(t,n)=g(nt)$
satisfies the assumptions of \rfr{P-genorb} (with $\om=g(0)$),
thus there exists a FU~subnilmanifold $Y$
such that $\ovr{\{g(nt),\ n\in\Z^{d}\}}\sle Y$ for all $t$
and $=Y$ for \ae\ $t\in\R^{d}$.
But then $Y=\ovr{\{g(t),\ t\in\R^{d}\}}$, 
and so $Y$ is a connected subnilmanifold;
by the second part of \rfr{P-disorb},
the sequence $g(nt)$, $n\in\Z^{d}$, is well distributed in $Y$
for \ae\ $t\in\R^{d}$.
Applying \rfr{P-mulWdis}, we get:
\theorem{P-nilorb}{}{}
Let $X$ be a compact nilmanifold and $g\col\R^{d}\ra X$ be a polynomial mapping.
Then $Y=\ovr{\{g(t),\ t\in\R^{d}\}}$ is a connected subnilmanifold of $X$
and $g(t)$ is well distributed in $Y$ (with respect to the Haar measure in $Y$).
\endtheorem

\remark{}
If we were only interested in proving 
the well distribution of $g$ in a subnilmanifold $Y$,
we could avoid the usage of \rfr{P-genorb};
we need it to show that $g(t)\in Y$ for all $t$.
\endremark
\tsubsection{s-pollim}{Convergence of multiple averages}

Combining \rfr{P-charfac} and \rfr{P-nilorb}, we can now get the following theorem:

\theorem{P-pollim}{}{}
Let $T^{t}$, $t\in\R$, 
be a continuous $1$-parameter group of measure preserving transformations
of a probability measure space $X$,
and let $p_{1},\ld,p_{r}$ be polynomials $\R^{d}\ra\R$.
Then for any $f_{1},\ld,f_{r}\in L^{\infty}(X)$
and any F{\o}lner sequence $(\Phi_{N})$ in $\R^{d}$
the limit 
\equ{
\lim_{N\ras\infty}\frac{1}{w(\Phi_{N})}
\int_{\Phi_{N}}T^{p_{1}(x)}f_{1}\cd\ld\cd T^{p_{r}(x)}f_{r}\,dx
}
exists in $L^{1}$-norm.
\endtheorem
\npar
(In \brfr{Potts}
a version of \rfr{P-pollim} was obtained for ``standard'' Ces\`{a}ro averages
(that is, for the case $\Phi_{N}=\prod_{i=1}^{d}[0,b_{i,N}]$, $N\in\N$,
with $b_{i,N}\ra\infty$ as $N\ra\infty$ for all $i=1,\ld,d$).
In \brfr{Austin-cp}, 
a multidimensional (that is, for $T\col\R^{d}\ra\R^{c}$ with $c\geq 1$) 
version of this result was obtained, again, for the standard Ces\`{a}ro averages.)

\proof{}
We may assume that $T$ is ergodic.
Applying \rfr{P-charfac} we can then replace $(X,T)$ 
by a pro-nilmanifold $Z_{k}(T)$.
Now, given the functions $f_{1},\ld,f_{r}\in L^{\infty}(X)$,
we can approximate them in $L^{1}$-norm by functions
that come from a factor $Y$ of $Z_{k}(T)$ which is a nilmanifold,
and replace $Z_{k}(T)$ by $Y$ and $T$ by a nilrotation $a$ on it.
Next, we note that it is enough to assume that $f_{1},\ld,f_{r}$
are continuous functions on $Y$.
Then an application of \rfr{P-nilorb}
to the polynomial flow $(a^{p_{1}(x)}y,\ld,a_{r}^{p_{r}(x)}y)$, $x\in\R^{d}$,
on the nilmanifold $Y^{r}$
and the function $f_{1}(y_{1})\cd\ld\cd f_{r}(y_{r})\in C(Y^{r})$
proves that the limit $\lim_{N\ras\infty}\frac{1}{w(\Phi_{N})}\int_{\Phi_{N}}
f_{1}(a^{p_{1}(t)}y)\cd\ld\cd f_{r}(a^{p_{r}(t)}y)dt$ exists
for all $y\in Y$, and so, in $L^{1}(Y)$.
\endproof 

Another way to prove \rfr{P-pollim} is to deduce it,
with the help of either \rfr{P-addF} or \rfr{P-mulF},
from the following discrete-time theorem:
\theorem{P-polmdlim}{}{(\brfr{Johnson})}
For any totally ergodic discrete $c$-parameter commutative group $T^{m}$, $m\in\Z^{c}$,
of measure preserving transformations of a probability measure space $X$,
any finite system of polynomials $p_{i}\col\Z^{d}\ra\Z^{c}$, $i=1,\ld,r$,
any $f_{1},\ld,f_{r}\in L^{\infty}(X)$,
and any F{\o}lner sequence $(\Psi_{N})$ in $\Z^{d}$,
the limit 
\equ{
\lim_{N\ras\infty}\frac{1}{|\Psi_{N}|}
\sum_{n\in\Psi_{N}}T^{p_{1}(n)}f_{1}\cd\ld\cd T^{p_{r}(n)}f_{r}
}
exists in $L^{1}$-norm.
\endtheorem

Applying \rfr{P-mulF}, 
we obtain from \rfr{P-polmdlim} the following refinement of \rfr{P-pollim}:
\theorem{P-poleqlim}{}{}
Under the assumptions of \rfr{P-pollim},
\equ{
\lim_{N\ras\infty}\frac{1}{w(\Phi_{N})}
\int_{\Phi_{N}}T^{p_{1}(x)}f_{1}\cd\ld\cd T^{p_{r}(x)}f_{r}\,dx
=\lim_{N\ras\infty}\frac{1}{|\Psi_{N}|}
\sum_{n\in\Psi_{N}}T^{p_{1}(nt)}f_{1}\cd\ld\cd T^{p_{r}(nt)}f_{r}
}
for \ae\ $t\in\R^{d}$ 
and any F{\o}lner sequences $(\Phi_{N})$ in $\R^{d}$ and $(\Psi_{N})$ in $\Z^{d}$.
\endtheorem

As for the actions of several commuting operators,
the following ``linear'' result has been recently obtained:
\theorem{P-dislinlim}{}{(\brfr{Austin-l}; see also \brfr{Host})}
Let $T_{1},\ld,T_{r}$ 
be pairwise commuting measure preserving transformations
of a probability measure space $X$.
Then for any $f_{1},\ld,f_{r}\in L^{\infty}(X)$ 
the limit 
$$
\lim_{b-a\ras\infty}\frac{1}{b-a}\sum_{n=a+1}^{b}
\dsc T_{1}^{n}f_{1}\cd\ld\cd T_{r}^{n}f_{r}
$$
exists in $L^{1}$-norm.
\endtheorem

Applying either \rfr{P-mulF} we obtain:
\theorem{P-linlim}{}{}
Let $T_{1}^{t},\ld,T_{r}^{t}$, $t\in\R$, 
be pairwise commuting continuous $1$-parameter groups of measure preserving transformations
of a probability measure space $X$.
Then for any $f_{1},\ld,f_{r}\in L^{\infty}(X)$
the limit 
$$
\lim_{b-a\ras\infty}\frac{1}{b-a}
\int_{a}^{b}T_{1}^{x}f_{1}\cd\ld\cd T_{r}^{x}f_{r}\,dx
$$
exists in $L^{1}$-norm,
and equals
$\lim_{b-a\ras\infty}\frac{1}{b-a}\sum_{n=a+1}^{b}
\dsc T_{1}^{nt}f_{1}\cd\ld\cd T_{r}^{nt}f_{r}$
for \ae\ $t\in\R$.
\endtheorem
\tsubsection{s-szem}{The polynomial Szemer\'{e}di theorem}

The ``multiparameter multidimensional polynomial ergodic Szemer\'{e}di theorem'' says:
\theorem{P-polDSz}{}{(See \brfr{BR} or \brfr{BLM}.)}
Let $T^{m}$, $m\in\Z^{c}$, be a discrete $c$-parameter commutative group 
of measure preserving transformations of a probability measure space $(X,\mu)$,
let $p_{i}\col\Z^{d}\ra\Z^{c}$, $i=1,\ld,r$, be a system of polynomials 
with $p_{i}(0)=0$ for all $i$,
and let $A\sle X$, $\mu(A)>0$.
Then for any F{\o}lner sequence $(\Psi_{N})$ in $\Z^{d}$,
\equ{
\liminf_{N\ras\infty}\frac{1}{|\Psi_{N}|}\sum_{n\in\Psi_{N}}
\mu\bigl(T^{p_{1}(n)}(A)\cap\ld\cap T^{p_{r}(n)}(A)\bigr)>0.
}
\endtheorem
Since the convergence of the averages
$\lim_{N\ras\infty}\frac{1}{|\Psi_{N}|}\sum_{n\in\Psi_{N}}
\mu\bigl(T^{p_{1}(n)}(A)\cap\ld\cap T^{p_{r}(n)}(A)\bigr)$
is unknown, we cannot apply \rfr{P-addF} or \rfr{P-mulF} 
to get a continuous-time version of \rfr{P-polDSz};
however it can be obtained with the help of either \rfr{P-addSF} or \rfr{P-mulSF}:
\theorem{P-polSz}{}{}
Let $T^{t}$, $t\in\R^{c}$, be a $c$-parameter commutative group 
of measure preserving transformations of a probability measure space $(X,\mu)$,
let $p_{i}\col\R^{d}\ra\R^{c}$, $i=1,\ld,r$, be a system of polynomials 
with $p_{i}(0)=0$ for all $i$,
and let $A\sle X$, $\mu(A)>0$.
Then for any F{\o}lner sequence $(\Phi_{N})$ in $\R^{d}$,
\equ{
\liminf_{N\ras\infty}\frac{1}{w(\Phi_{N})}\int_{\Phi_{N}}
\mu\bigl(T^{p_{1}(x)}(A)\cap\ld\cap T^{p_{r}(x)}(A)\bigr)\,dx>0.
}
\endtheorem
\ignore
It is a more delicate problem -- 
to find the lower bound of the limits (or liminfs) \frfr{f-pSz}
for fixed $p_{1},\ld,p_{r}$ and $\mu(A)$.
An example is {\it the problem of large limits}:
we say that a system of polynomials $p_{1},\ld,p_{r}$, with $p_{i}(0)=0$ for all $i$,
is a system of {\it large limits\/}
if for any measure preserving transformation $T$ of a probability measure space $X$
and any measurable subset $A$ of $X$ one has
$\liminf_{N\ras\infty}\frac{1}{|\Psi_{N}|}\sum_{n\in\Psi_{N}}
\mu\bigl(T^{p_{1}(n)}(A)\cap\ld\cap T^{p_{r}(n)}(A)\bigr)\geq\mu(A)^{r}$.
Dealing with this system, \rfr{P-addSF} is powerless,
but \rfr{P-mulF}, via \rfr{P-poleqlim}, may help.
\endignore

A ($d$-parameter) {\it polynomial sequence\/} in a group $G$
is a sequence of the form $g(n)=\prod_{j=1}^{k}v_{j}^{p_{j}(n)}$,
where $v_{j}$ are elements of $G$ and $p_{j}$ are integer-valued polynomials on $\Z^{d}$.
\rfr{P-polDSz} was extended in \brfr{nsz} to the nilpotent set-up:
\theorem{P-nilDSz}{}{}
Let $G$ be a nilpotent group 
of measure preserving transformations of a probability measure space $(X,\mu)$,
let $g_{i}\col\Z^{d}\ra G$, $i=1,\ld,r$, 
be a system of $d$-parameter polynomial sequences in $G$
with $g_{i}(0)=1_{G}$ for all $i$,
and let $A\sle X$, $\mu(A)>0$.
Then
\equ{
\liminf_{l(b-a)\ras\infty}\frac{1}{w(b-a)}\sum_{n>a}^{b}
\mu\bigl((g_{1}(n))(A)\cap\ld\cap(g_{r}(n))(A)\bigr)>0.
}
\endtheorem

If $G$ a connected nilpotent Lie group, 
then for any $v\in G$ there exists a one-parameter subgroup $v^{t}$, $t\in\R$, of $G$
such that $v^{1}=v$;
this allows one to define $v^{t}$ for all $t\in\R$.
Let us call {\it a polynomial mapping\/} $g\col\R^{d}\ra G$
a mapping of the form $g(x)=\prod_{j=1}^{k}v_{j}^{p_{j}(x)}$,
where $v_{j}$ are elements of $G$ and $p_{j}$ are polynomials on $\R^{d}$.
Applying one of Theorems~\rfrn{P-addSF} or \rfrn{P-mulSF},
we get the following ``continuous-time nilpotent polynomial Szemer\'{e}di theorem'':
\theorem{P-nilSz}{}{}
Let $G$ be a nilpotent Lie group
of measure preserving transformations of a probability measure space $(X,\mu)$,
let $g_{i}\col\R^{d}\ra G$, $i=1,\ld,r$, be a system of polynomial mappings
with $g_{i}(0)=1_{G}$ for all $i$,
and let $A\sle X$, $\mu(A)>0$.
Then for any F{\o}lner sequence $(\Phi_{N})$ in $\R^{d}$,
\equ{
\liminf_{N\ras\infty}\frac{1}{w(\Phi_{N})}\int_{\Phi_{N}}
\mu\bigl((g_{1}(x))(A)\cap\ld\cap(g_{r}(x))(A)\bigr)\,dx>0.
}
\endtheorem
\tsubsection{s-genpol}{Distribution of values of generalized polynomials}

Another application of \rfr{P-nilorb} is a sharpening of the results from \brfr{sko} 
about the distribution of values of bounded generalized polynomials.
Recall that a generalized polynomial 
is a function from $\R^{d}$ or from $\Z^{d}$ to $\R$ 
that is constructed from conventional polynomals
by applying the operaions of addition, multiplication,
and taking the integer part.
We call a function $u\col\R^{d}\ra\R^{c}$
a {\it generalized polynomial mapping}
if all components of $u$ are generalized polynomials.
Under a {\it piecewise polynomial surface\/} $\Sp\sle\R^{c}$
we understand the image $\Sp=S(Q)$ of the cube $Q=[0,1]^{s}$
where $S$ is a {\it piecewise polynomial mapping},
which means that $Q$ can partitioned
into a finite union $Q=\bigcup_{i=1}^{l}Q_{i}$ of subsets
so that for each $i$,
$Q_{i}$ is defined by a system of polynomial inequalities
and $S\rest{Q_{i}}$ is a polynomial mapping.
We endow $\Sp$ with the measure $\mu_{\Sp}=S_{*}(w)$,
the push-forward of the standard Lebesgue measure $w$ on $Q$.
In \brfr{sko}, the following theorem was proved:
\theorem{P-DisGP}{}{(\brfr{sko})}
Let $u\col\Z^{d}\ra\R^{c}$ be a bounded generalized polynomial mapping.
Then the sequence $u(n)$, $n\in\Z^{d}$
is well distributed with respect to $\mu_{\Sp}$
on a piecewise polynomial surface $\Sp\sln\R^{c}$.
\endtheorem
\npar
(Note that it is not claimed in this theorem that $u(n)\in\Sp$ for all $n$;
it follows however that the set $\bigl\{n:u(n)\not\in\Sp\bigr\}$
has zero uniform density in $\Z^{d}$.)

Applying \rfr{P-mulWdis}, 
we may now obtain the $\R$-version of \rfr{P-DisGP}:
\theorem{P-genpol}{}{}
Any bounded generalized polynomial mapping $u\col\R^{d}\ra\R^{c}$
is well distributed on a piecewise polynomial surface $\Sp\sln\R^{c}$.
\endtheorem

An application of the spectral theorem gives, as a corollary,
the following proposition:
\proposition{P-spegenpol}{}{}
Let $U^{t}$, $t\in\R^{c}$, 
be a continuous $c$-parameter group of unitary operators
on a Hilbert space $\cH$,
and let $u\col\R^{d}\ra\R^{c}$ be a generalized polynomial mapping.
Then for any $v\in\cH$ and any F{\o}lner sequence $(\Phi_{N})$ in $\R^{d}$,
the limit
$\lim_{N\ras\infty}\frac{1}{w(\Phi_{N})}\int_{\Phi_{N}}U^{g(x)}v\,dx$ 
exists.
\endproposition
\tsubsection{s-hardywm}{Ergodic theorems along functions from Hardy fields}

We will now deal with a situation
where our ``uniform Ces\'{a}ro theorems'' are not applicable,
but the ``standard Ces\`{a}ro'' \rfr{P-mulC} is;
namely, we will deal with multiple ergodic averages 
along (not necessarily) polynomial functions of polynomial growth.
Such averages for functions of integer argument
were considered in \brfr{Inger} and \brfr{Nikos2}.

To state the results obtained in \brfr{Inger}
we first need to introduce some notation:

$\cT$ is the set of real-valued $C^{\infty}$ functions $g$
defined on intervals $[a,\infty)$, $a\in\R$,
such that a finite $\lim_{x\ras+\infty}xg^{(j+1)}(x)/g^{(j)}(x)$ exists
for all $j\in\N$ and there exists an integer $i\geq 0$ and $\alf\in(i,i+1]$
such that $\lim_{x\ras+\infty}xg'(x)/g(x)=\alf$
and $\lim_{x\ras+\infty}g^{(i+1)}(x)=0$;

$\cP$ is the set of real-valued $C^{\infty}$ functions $g$
defined on intervals $[a,\infty)$, $a\in\R$
such that for some integer $i\geq 0$
a finite nonzero $\lim_{x\ras+\infty}g^{(i+1)}(x)$ exists
and $\lim_{x\ras+\infty}x^{j}g^{(i+j+1)}(x)=0$ for all $j\in\N$;

$\cG=\cT\cup\cP$;

$\cL$ is the Hardy field of logarithmico-exponential functions, 
that is, the minimal field of real-valued functions 
defined on intervals $[a,\infty)$, $a\in\R$,
that contains polynomials
and is closed under the operations of taking exponent and logarithm-of-modulus;

for $\alf>0$, $\cG(\alf)$ is the set of functions $g\in\cG$
with $\lim_{x\ras+\infty}xg'(x)/g(x)=\alf$,
$\cT(\alf)$ is the set of functions $g\in\cT$
with $\lim_{x\ras+\infty}xg'(x)/g(x)=\alf$,
and for any $G\sle\cG$, $G(\alf)=G\cap\cG(\alf)$;

a finite family $G\sln\cG$ with $g_{1}-g_{2}\in\cG$ for all $g_{1},g_{2}\in\cG$
is said to have {\it R-property\/}
if for any $\alf>0$, any $g_{1},g_{2}\in(G(\alf)\cup(G(\alf)-G(\alf)))\sm\{0\}$,
any integer $l\geq 0$ and $\bet\in(0,\alf)$
such that $g_{1}^{(l)},g_{2}\in\cT(\bet)$,
a finite nonzero $\lim_{x\ras+\infty}g_{1}^{([\bet]+l+1)}(x)/g_{2}^{([\bet]+1)}(x)$ exists.

\vbreak{-9000}{3mm}
The following theorem was proved in \brfr{Inger}:
\theorem{P-Inger}{}{(\brfr{Inger})}
Let $g_{1},\ld,g_{r}\in\cG$
be such that $g_{i}-g_{j}\in\cG$ for all $i\neq j$,
and also either $g_{1},\ld,g_{r}\in\cL$
or the family $\{g_{1},\ld,g_{r}\}$ has the R-property.
Then for any invertible weakly mixing transformation $T$
of a probability measure space $(X,\mu)$
and any $f_{1},\ld,f_{r}\in L^{\infty}(X)$,
the sequence $F_{n}=T^{[g_{1}(n)]}f_{1}\cd\ld\cd T^{[g_{r}(n)]}f_{r}$, $n\in\N$,
tends in density in $L^{1}$-norm to $\prod_{i=1}^{r}\int f_{i}\,d\mu$.
\endtheorem

The statement ``$F_{n}$ tends in density in $L^{1}$-norm'' means that
$\lim_{N\ras\infty}
\frac{1}{N}\sum_{n=1}^{N}\bigl\|T^{[g_{1}(n)]}f_{1}\cd\ld\cd T^{[g_{r}(n)]}f_{r}
-\prod_{i=1}^{r}\int f_{i}\,d\mu\bigr\|_{L^{1}(X)}\dsc=0$.
From this and \rfr{P-mulC} we get that,
under the assumptions of \rfr{P-Inger},
$\lim_{b\ras\infty}
\frac{1}{b}\int_{0}^{b}\bigl\|T^{[g_{1}(x)]}f_{1}\cd\ld\cd T^{[g_{r}(x)]}f_{r}\,dx
-\prod_{i=1}^{r}\int f_{i}\,d\mu\bigr\|_{L^{1}(X)}\,dx=0$,
that is, the function
$F_{x}=T^{[g_{1}(x)]}f_{1}\cd\ld\cd T^{[g_{r}(x)]}f_{r}$, $x\in[0,\infty)$,
(whose range is in $L^{1}(X)$)
tends in density in $L^{1}$-norm to $\prod_{i=1}^{r}\int f_{i}\,d\mu$.
Hence, we obtain:
\theorem{P-Ingersk}{}{}
Let $g_{1},\ld,g_{r}\in\cG$
be such that $g_{i}-g_{j}\in\cG$ for all $i\neq j$,
and also either $g_{1},\ld,g_{r}\in\cL$ 
or the family $\{g_{1},\ld,g_{r}\}$ has the R-property.
Then for any invertible weakly mixing transformation $T$
of a probability measure space $(X,\mu)$
and any $f_{1},\ld,f_{r}\in L^{\infty}(X)$,
the  function $F_{x}=T^{[g_{1}(x)]}f_{1}\cd\ld\cd T^{[g_{r}(x)]}f_{r}$, $x\in[0,\infty)$,
tends in density in $L^{1}$-norm to $\prod_{i=1}^{r}\int f_{i}\,d\mu$.
\endtheorem

Actually, one can eliminate the brackets appearing in the exponents
in the expression for $F_{x}$.
Indeed, 
put $G_{x}=T^{g_{1}(x)}f_{1}\cd\ld\cd T^{g_{r}(x)}f_{r}$, $x\in[0,\infty)$,
and let $L=\prod_{i=1}^{r}\int f_{i}\,d\mu$.
Assume that $\|f_{i}\|\leq 1$ for all $i$.
Fix any $\eps>0$
and for each $i=1,\ld,r$ choose functions $g_{i,j}\in L^{\infty}(X)$, $j=1,\ld,k$,
that form an $\eps$-net in the (compact) set $\{T^{t}f_{i},\ t\in[0,1]\}\sln L^{1}(X)$.
For any $J=(j_{1},\ld,j_{r})\in\{1,\ld,k\}^{r}$
the function $(F_{J})_{x}=T^{[g_{1}(x)]}f_{1,j_{1}}\ld T^{[g_{r}(x)]}f_{r,j_{r}}$
tends in density to $L$,
and for any $x\in[0,\infty)$ there exists $J=(j_{1},\ld,j_{r})\in\{1,\ld,k\}^{r}$
such that $\|T^{g_{i}(x)}f_{i}-T^{[g_{i}(x)]}f_{i,j_{i}}\|
=\|T^{\{g_{i}(x)\}}f_{i}-f_{i,j_{i}}\|<\eps$ for all $i$
and so, $\|G_{x}-(F_{J})_{x}\|<2^{2r}\eps$.
This implies that $\limsup_{N\ras\infty}
\frac{1}{N}\sum_{1}^{N}\|G_{x}-L\|<2^{2r}\eps$.
Since this holds for any positive $\eps$,
we see that $G_{x}$ also tends in density to $L$.
So, we have the following result:
\theorem{P-Ingerns}{}{}
Let $g_{1},\ld,g_{r}\in\cG$
be such that $g_{i}-g_{j}\in\cG$ for all $i\neq j$,
and also either $g_{1},\ld,g_{r}\in\cL$ 
or the family $\{g_{1},\ld,g_{r}\}$ has the R-property.
Then for any weakly mixing continuous 1-parameter group $T^{t}$, $t\in\R$,
of measure preserving transformations
of a probability measure space $(X,\mu)$
and any $f_{1},\ld,f_{r}\in L^{\infty}(X)$,
the function $G_{x}=T^{g_{1}(x)}f_{1}\cd\ld\cd T^{g_{r}(x)}f_{r}$, $x\in[0,\infty)$,
tends in density in $L^{1}$-norm to $\prod_{i=1}^{r}\int f_{i}\,d\mu$.
\endtheorem

\vbreak{-9000}{3mm}
Another paper dealing with multiple-ergodic averages
along non-polynomial functions of polynomial growth 
is \brfr{Nikos2}.
Let $\cH$ denote the union of all Hardy fields of real-valued functions.
\theorem{P-Nikos1}{}{(\brfr{Nikos2})}
Let $g\in\cH$ satisfy $\lim_{x\ras+\infty}g(x)/x^{j}=0$ for some $j\in\N$,
and assume that one of the following is true:
either $\lim_{x\ras+\infty}(g(x)-cp(x))/\log x=\infty$ for all $c\in\R$ and $p\in\Z[x]$;
or $\lim_{x\ras+\infty}(g(x)-cp(x))=d$ for some $c,d\in\R$ and $p\in\Z[x]$;
or $(g(x)-x/m)/\log x$ is bounded on $[2,\infty)$ for some $m\in\Z$.
Then for any invertible measure preserving transformation
of a probability measure space $X$,
$\lim_{N\ras\infty}\frac{1}{N}\sum_{n=1}^{N}
T^{[g(n)]}f_{1}\cd T^{2[g(n)]}f_{2}\cd\ld\cd T^{r[g(n)]}f_{r}$
exists in $L^{1}(X)$ for any $r\in\N$ and any $f_{1},\ld,f_{r}\in L^{\infty}(X)$.
\endtheorem

\theorem{P-Nikos2}{}{(\brfr{Nikos2})}
Let $g_{1},\ld,g_{r}\in\cL$ be logarithmico-exponential functions
satisfying $\lim_{x\ras+\infty}g_{i}(x)/x^{k_{i}+1}\dsc
=\lim_{x\ras+\infty}x^{k_{i}+\eps_{i}}/g_{i}(x)=0$
for some integer $k_{i}\geq 0$ and $\eps_{i}>0$, $i=1,\ld,r$,
and $\lim_{x\ras+\infty}g_{i}(x)/g_{j}(x)=0$ or $\infty$
for any $i\neq j$.
Then for any invertible ergodic measure preserving transformation
of a probability measure space $X$,
$\lim_{N\ras\infty}\frac{1}{N}\sum_{n=1}^{N}
T^{[g_{1}(n)]}f_{1}\cd\ld\cd T^{[g_{r}(n)]}f_{r}
=\prod_{i=1}^{r}\int_{X}f_{i}\,d\mu$ in $L^{1}(X)$ 
for any $f_{1},\ld,f_{r}\in L^{\infty}(X)$.
\endtheorem

From this and \rfr{P-mulC} we get
\theorem{P-Nikos1sk}{}{}
Let $g\in\cH$ satisfy $\lim_{x\ras+\infty}g(x)/x^{j}=0$ for some $j\in\N$,
and assume that one of the following is true:
either $\lim_{x\ras+\infty}(g(x)-cp(x))/\log x=\infty$ for all $c\in\R$ and $p\in\Z[x]$;
or $\lim_{x\ras+\infty}(g(x)-cp(x))=d$ for some $c,d\in\R$ and $p\in\Z[x]$;
or $(g(x)-x/m)/\log x$ is bounded on $[2,\infty)$ for some $m\in\Z$.
Then for any invertible measure preserving transformation
of a probability measure space $X$,
$\lim_{b\ras\infty}\frac{1}{b}\int_{0}^{b}
T^{[g(x)]}f_{1}\cd T^{2[g(x)]}f_{2}\cd\ld\cd T^{r[g(x)]}f_{r}\,dx$
exists in $L^{1}(X)$ for any $r\in\N$ and any $f_{1},\ld,f_{r}\in L^{\infty}(X)$.
\endtheorem

\theorem{P-Nikos2sk}{}{}
Let $g_{1},\ld,g_{r}\in\cL$ be logarithmico-exponential functions
satisfying $\lim_{x\ras+\infty}g_{i}(x)/x^{k_{i}+1}\dsc
=\lim_{x\ras+\infty}x^{k_{i}+\eps_{i}}/g_{i}(x)=0$
for some integer $k_{i}\geq 0$ and $\eps_{i}>0$, $i=1,\ld,r$,
and $\lim_{x\ras+\infty}g_{i}(x)/g_{j}(x)=0$ or $\infty$
for any $i\neq j$.
Then for any invertible ergodic measure preserving transformation
of a probability measure space $X$,
$\lim_{b\ras\infty}\frac{1}{b}\int_{0}^{b}
T^{[g_{1}(x)]}f_{1}\cd\ld\cd T^{[g_{r}(x)]}f_{r}\,dx
=\prod_{i=1}^{r}\int_{X}f_{i}\,d\mu$ in $L^{1}(X)$ 
for any $f_{1},\ld,f_{r}\in L^{\infty}(X)$.
\endtheorem

\tsubsection{s-pointwise}{Pointwise ergodic theorems}

Here are two theorems of Bourgain dealing with the pointwise convergence:
\theorem{P-BopntLin}{}{(\brfr{Bour1})}
Let $T$ be a measure preserving transformation 
of a probability measure space $X$.
Then for any $f_{1},f_{2}\in L^{2}(X)$,
the sequence $\frac{1}{N}\sum_{n=1}^{N}T^{n}f_{1}\cd T^{2n}f_{2}$, $N\in\N$,
converges \ae.
\endtheorem

\theorem{P-BopntPol}{}{(\brfr{Bour2})}
Let $T_{1},\ld,T_{r}$ be commuting invertible measure preserving transformations 
of a probability measure space $X$.
Then for any $f\in L^{2}(X)$
and any polynomials $p_{1},\ld,p_{r}\col\Z\ra\Z$,
the sequence
$\frac{1}{N}\sum_{n=1}^{N}\bigl(\prod_{i=1}^{r}T_{i}^{p_{i}(n)}\bigr)f$, $N\in\N$,
converges \ae.
\endtheorem

We now have:
\theorem{P-CpntLin}{}{}
Let $T^{t}$, $t\in\R$, be a continuous action of the semigroup $[0,\infty)$ 
by measure preserving transformations on a probability measure space $X$.
Then for any $f_{1},f_{2}\in L^{2}(X)$,
$\lim_{b\ras\infty}\frac{1}{b}\int_{0}^{b}T^{t}f_{1}\cd T^{2t}f_{2}\,dt$
exists \ae.
\endtheorem
\proof{}
By \rfr{P-BopntLin}, for every $t\in\R$,
the sequence $\frac{1}{N}\sum_{n=1}^{N}T^{nt}f_{1}(\om)\cd T^{2nt}f_{2}(\om)$, $N\in\N$,
converges for \ae\ $\om\in X$;
let $S_{t}\sln X$ be the set of points $\om$ for which this is not so.
Then $\bigl\{(t,\om):\om\in S_{t}\bigr\}$ is a null-subset of $\R\times X$,
thus for \ae\ $\om\in X$,
the limit $\lim_{N\ras\infty}\frac{1}{N}\sum_{n=1}^{N}T^{nt}f_{1}(\om)\cd T^{2nt}f_{2}(\om)$ exists
for \ae\ $t\in\R$.
By (the scalar version of) \rfr{P-mulC},
the limit $\lim_{b\ras\infty}\frac{1}{b}\int_{0}^{b}T^{t}f_{1}(\om)\cd T^{2t}f_{2}(\om)\,dt$ exists
for \ae\ $\om\in X$.
\endproof

In the same way, from \rfr{P-BopntPol} we get:
\theorem{P-CpntPol}{}{}
Let $T^{t}$, $t\in\R^{c}$, be a continuous $c$-parameter group
of measure preserving transformations of a probability measure space $X$.
Then for any $f\in L^{2}(X)$
and any polynomial $p\col\R\ra\R^{c}$,
$\lim_{b\ras\infty}\frac{1}{b}\int_{0}^{b}T^{p(t)}f\,dt$
exists \ae.
\endtheorem

Here are two more pointwise theorems, established by Assani:
\theorem{P-Ass1}{}{(\brfr{Ass})}
Let $T$ be a weakly mixing measure preserving transformation 
of a probability measure space $X$,
let $(P,S)$ be the Pinsker factor of $(X,T)$,
and assume that the spectrum of $S$ is singular.
Then for any $f_{1},\ld,f_{r}\in L^{\infty}(X)$,
the sequence $\frac{1}{N}\sum_{n=1}^{N}T^{n}f_{1}\cd\ld\cd T^{rn}f_{r}$, $N\in\N$,
converges to $\prod_{i=1}^{r}\int f_{i}\,d\mu$ \ae\ on $X$.
\endtheorem

\theorem{P-Ass2}{}{(\brfr{Ass})}
Let $T$ be a weakly mixing measure preserving transformation
of a probability measure space $X$,
let $(P,S)$ be the Pinsker factor of $(X,T)$,
let $L\sle L^{2}(P)$ be the space of functions on $P$
whose spectral measure under the action of $S$
is absolutely continuous with respect to the Lebesgue measure.
For any function $f\in L^{2}(X)$ let $\hf$ denote the projection of $E(f|P)$ to $L$.
Then for any $f_{1},f_{2},f_{3}\in L^{\infty}(X)$,
$$
\lim_{N\ras\infty}\Bigl(\frac{1}{N}\sum_{n=1}^{N}T^{n}f_{1}\cd T^{2n}f_{2}\cd T^{3n}f_{3}
-\frac{1}{N}\sum_{n=1}^{N}T^{n}\hf_{1}\cd T^{2n}\hf_{2}\cd T^{3n}\hf_{3}\Bigr)=0\
\hbox{\ae}
$$
\endtheorem

Let $T^{t}$, $t\in\R$, be a continuous action of $\R$
by measure preserving transformations on a measure space $X$.
Then, with the help of either \rfr{P-addC} or \rfr{P-mulC}, 
repeating (the first two phrases from) the proof of \rfr{P-CpntLin},
and taking into account that
(i) if $T$ is weakly mixing then $T^{t}$ is weakly mixing for all $t\neq0$;
(ii) the Pinsker algebra of $T$ is the Pinsker algebra of $T^{t}$ for all $t\neq0$; and
(iii) if the spectrum of $T$ is singular (respectively, absolutely continuous),
then the spectrum of $T^{t}$ is singular (respectively, absolutely continuous)
for all $t\neq0$;
we obtain:
\theorem{P-xAss1}{}{}
Let $T$ be a continuous action of $\R$
on a probability measure space $X$ by weakly mixing measure preserving transformations,
let $(P,S)$ be the Pinsker factor of $(X,T)$,
and assume that the spectrum of $S$ is singular.
Then for any $f_{1},\ld,f_{r}\in L^{\infty}(X)$ one has
$\lim_{b\ras\infty}\frac{1}{b}\int_{0}^{b}T^{t}f_{1}\cd\ld\cd T^{rt}f_{r}\,dt=
\prod_{i=1}^{r}\int f_{i}\,d\mu$ \ae.
\endtheorem

\theorem{P-xAss2}{}{}
Let $T$ be a continuos action of $\R$
on a probability measure space $X$ by weakly mixing measure preserving transformations
let $(P,S)$ be the Pinsker factor of $(X,T)$,
let $L\sle L^{2}(P)$ be the space of functions on $P$
whose spectral measure under the action of $S$
is absolutely continuous with respect to the Lebesgue measure.
For any function $f\in L^{2}(X)$ let $\hf$ denote the projection of $E(f|P)$ to $L$.
Then for any $f_{1},f_{2},f_{3}\in L^{\infty}(X)$,
$$
\lim_{b\ras\infty}\Bigl(\frac{1}{b}\int_{0}^{b}T^{t}f_{1}\cd T^{2t}f_{2}\cd T^{2t}f_{3}\,dt
-\frac{1}{b}\int_{0}^{b}T^{t}\hf_{1}\cd T^{2t}\hf_{2}\cd T^{3t}\hf_{3}\,dt\Bigr)=0\
\hbox{\ae}
$$
\endtheorem
\bibliography{}
\biblleft=13mm
\bibart Ass/A
a:I. Assani
t:Multiple recurrence and almost sure convegrence for weakly mixing dynamical systems
j:Israel Journal of Math
n:163
y:1998
p:111-124
*
\bibart Austin-l/Au1
a:T. Austin
t:On the norm convergence of nonconventional ergodic averages
j:Ergodic Theory and Dynamical Systems 
n:30
y:2010
p:\no{2}, 321-338
*
\bibartpr Austin-cp/Au2
a:T. Austin
t:Norm convergence of continuous-time polynomial multiple ergodic averages
*
\bibart Inger/BK
a:V. Bergelson and I. H{\aa}land Knutson
t:Weak mixing implies weak mixing of higher orders along tempered functions
j:Ergodic Theory and Dynamical Systems
n:29
y:2009
p:1375-1416
*
\bibart sko/BL
a:V. Bergelson and A. Leibman
t:Distribution of values of bounded generalized polynomials
j:Acta Mathematica
n:198
y:2007
p:155-230
*
\bibart BLM/BLM
a:V. Bergelson, A. Leibman, and R. McCutcheon
t:Polynomial Szemeredi theorem for countable modules over integral domains and finite fields 
j:Journal d'Analyse Mathematique 
n:95
y:2005
p:243-296
*
\bibart isp/BLLe
a:V. Bergelson, A. Leibman, and E. Lesigne
t:Intersective polynomials and the polynomial Szemer\'{e}di theorem
j:Advances in Mathematics
n:219
y:2008
p:\no{1}, 369-388
*
\bibart BR/BM
a:V. Bergelson and R. McCutcheon
t:An ergodic IP polynomial Szemeredi theorem
j:Memoirs of AMS
n:146
y:2000
p:viii+106pp
*
\bibart BirKoo/BiKo
a:G.D. Birkhoff and B.O. Koopman
t:Recent contributions to the ergodic theory
j:Proceedings of the National Academy of Sciences of the USA
n:18
y:1932
p:279–282
*
\bibart Bour1/Bo1
a:J. Bourgain
t:Double recurrence and almost sure convergence
j:Journal f\"{u}r die Reine Angewandte Mathematik
n:404
y:1990
p:140-161
*
\bibart Bour2/Bo2
a:J. Bourgain
t:On the maximal ergodic theorem for certain subsets of the integers
j:Israel Journal of Mathematics
n:61
y:1988
p:\no{1}, 39-72
*
\bibart Nikos2/F
a:N. Frantzikinakis
t:Multiple recurrence and convergence for Hardy sequences of polynomial growth
j:Journal d'Analyse Mathematique
n:112
y:2010
p:79-135 
*
\bibart FrKr/FK
a:N. Frantzikinakis and B. Kra
t:Convergence of multiple ergodic averages for some commuting transformations
j:Ergodic Theory and Dynamical Systems
n:25
y:2005
p:\no{3}, 799-809
*
\bibook Hopf/H
a:E. Hopf 
t:Ergodentheorie
i:Springer, Berlin, 1937Multiple recurrence theorem for measure preserving actions of a nilpotent group, Geometric and Functional Analysis 8 (1998), 853-931 
*
\bibart Host/Ho1
a:B. Host
t:Ergodic seminorms for commuting transformations and applications
j:Studia Mathematica
n:195
y:2009
p:31-49
* 
\bibart HKo/HoK1
a:B. Host and B. Kra 
t:Non-conventional ergodic averages and nilmanifolds
j:Annals of Mathematics
n:161
y:2005
p:\no{1}, 397-488
* 
\bibart HKp/HoK2
a:B. Host and B. Kra 
t:Convergence of polynomial ergodic averages 
j:Israel Journal of Mathematics
n:149
y:2005
p:1-19
* 
\bibart Johnson/J
a:M. Johnson
t:Convergence of polynomial ergodic averages of several variables 
for some commuting transformations
j:Illinois Journal of Mathematics
n:53
y:2009
p:\no{3}, 865-882
*
\bibart Kolmogorov/Ko
a:A.N. Kolmogorov
t:A simplified proof of the Birkhoff-Khinchin ergodic theorem
j:Uspekhi Mat. Nauk
n:
y:1938
p:\no{5}, 52-56
*
\bibart nsz/L1
a:A. Leibman
t:Multiple recurrence theorem for measure preserving actions of a nilpotent group
j:Geometric and Functional Analysis 
n:8 
y:1998
p:853-931 
*
\bibart mpn/L2
a:A. Leibman
t:Pointwise convergence of ergodic averages 
for polynomial actions of $\Z^{d}$ by translations on a nilmanifold
j:Ergodic Theory and Dynamical Systems
n:25
y:2005
p:215-225
*
\bibart cpm/L3
a:A. Leibman
t:Convergence of multiple ergodic averages along polynomials of several variables
j:Israel Journal of Mathematics
n:146
y:2005
p:303--322
*
\bibart ran/L4
a:A. Leibman
t:Rational sub-nilmanifolds of a compact nilmanifold
j:Ergodic Theory and Dynamical Systems
n:26
y:2006
p:787-798
*
\bibart orb/L5
a:A. Leibman
t:Orbits on a nilmanifold 
under the action of a polynomial sequences of translations
j:Ergodic Theory and Dynamical Systems
n:27
y:2007
p:1239-1252
*
\bibart vonNeumann/vN
a:J. von Neumann
t:Physical applications of the ergodic hypothesis
j:Proceedings of the National Academy of Sciences of the USA
n:18
y:1932
p:\no{3} 263-266
*
\bibartpr Potts/P
a:A. Potts
t:Multiple ergodic averages for flows and an application
*
{\small(arXiv:0910.3687)}
\bibart Shah/Sh
a:N. Shah
t:Limit distributions of polynomial trajectories on homogeneous spaces
j:Duke Mathematical Journal
n:75
y:1994
p:no. 3, 711-732
*
\bibart Zo/Z
a:T. Ziegler
t:Universal characteristic factors and Furstenberg averages
j:Journal of AMS
n:20
y:2007
p:\no{1}, 53-97
* 
\endbibliography
\finish
\bye